\documentclass{article}

\usepackage{amsmath}
\usepackage{amsthm}
\usepackage{mathrsfs}
\usepackage{amsfonts}
\begin {document}
\topmargin= -.2in \baselineskip=20pt

\def\bbmu{\mu}

\newtheorem{theorem}{Theorem}[section]
\newtheorem{proposition}[theorem]{Proposition}
\newtheorem{lemma}[theorem]{Lemma}
\newtheorem{corollary}[theorem]{Corollary}
\newtheorem{conjecture}[theorem]{Conjecture}
\theoremstyle{remark}
\newtheorem{remark}[theorem]{Remark}

\title {Calculation of $\ell$-adic Local Fourier Transformations
\thanks{I would like to thank anonymous referees of various versions of
the paper for communicating to me their deep insight on the
stationary phase principle and the Legendre transformation, and for
their many suggestions on improving the exposition of the paper. I
also thank P. Deligne, L. Illusie and N. Katz for their help while
preparing this paper. The research is supported by the NSFC
(10525107).}}

\author {Lei Fu\\
{\small Chern Institute of Mathematics and LPMC, Nankai University,
Tianjin 300071, P. R. China}\\
{\small leifu@nankai.edu.cn}}

\date{}
\maketitle

\begin{abstract}
We calculate the local Fourier transformations for a class of
$\overline{\mathbb Q}_\ell$-sheaves. In particular, we verify a
conjecture of Laumon and Malgrange (\cite{L} 2.6.3). As an
application, we calculate the local monodromy of $\ell$-adic
hypergeometric sheaves introduced by Katz (\cite{K1}). We also
discuss the characteristic $p$ analogue of the Turrittin-Levelt
Theorem for $D$-modules.

\medskip
\noindent {\bf Key words:} Fourier transformation, local Fourier
transformation, hypergeometric sheaf.

\medskip
\noindent {\bf Mathematics Subject Classification:} 14F20, 14G15,
11L05.

\end{abstract}

\section*{Introduction}

The global $\ell$-adic Fourier transformation was first introduced
by Deligne. To study the local behavior of the global Fourier
transformation, Laumon \cite{L} discovered the stationary phase
principle and introduced local Fourier transformations. All these
transformations are defined by cohomological functors and are rarely
computable. However, in \cite{L} 2.6.3, Laumon and Malgrange give
conjectural formulas of local Fourier transformations for a class of
$\overline {\mathbb Q}_\ell$-sheaves. In this paper we calculate
local Fourier transformations for a more general class of
$\overline{\mathbb Q}_\ell$-sheaves. In particular, we prove the
conjecture of Laumon and Malgrange. It turns out that to get the
correct result, the conjectural formulas of Laumon and Malgrange
have to be slightly modified. As an application of our result, we
calculate the local monodromy of $\ell$-adic hypergeometric sheaves
introduced by Katz (\cite{K1}). The method used in this paper can
also be used to show that the Fourier transformation of an analytic
sheaf with meromorphic ramification defined in \cite{R} still has
meromorphic ramification, thus answering positively a question
proposed by Ramero.

\medskip
Throughout this paper, for any ring $A$, we use the notations
$${\mathbb A}_A^1=\mathrm {Spec}\, A[t], \; {\mathbb G}_{m,A}=\mathrm {Spec}\,
A[t,1/t], \; {\mathbb P}_A^1=\mathrm {Proj}\, A[t_0, t_1].$$ Let $k$
be a perfect field of characteristic $p$, $\bar k$ an algebraic
closure of $k$, $q$ any power of $p$, ${\mathbb F}_q$ the finite
field with $q$ elements contained in $\bar k$, and $\ell$ a prime
number distinct from $p$. Fix a nontrivial additive character
$\psi:{\mathbb F}_p\to \overline{\mathbb Q}_\ell^\ast$. The
$k$-morphism
$$\wp: {\mathbb A}_k^1\to {\mathbb A}_k^1,\quad t\mapsto t^p-t$$
is a finite Galois \'etale covering space, and it defines an
${\mathbb F}_p$-torsor
$$0\to {\mathbb F}_p\to {\mathbb A}_k^1\stackrel \wp\to {\mathbb A}_k^1\to
0.$$ Pushing-forward this torsor by $\psi^{-1}$, we get a lisse
$\overline{\mathbb Q}_\ell$-sheaf ${\mathscr L}_\psi$ of rank $1$ on
${\mathbb A}_k^1$, which we call the Artin-Schreier sheaf. Let $X$
be a scheme over $k$ and let $f$ be an element in the ring of global
sections $\Gamma(X,{\mathscr O}_X)$ of the structure sheaf of $X$.
Then $f$ defines a $k$-morphism $X\to {\mathbb A}_k^1$ so that the
induced $k$-algebra homomorphism $k[t]\to \Gamma(X,{\mathscr O}_X)$
maps $t$ to $f$. We often denote this canonical morphism also by
$f$, and denote by ${\mathscr L}_{\psi}(f)$ the inverse image of
${\mathscr L}_\psi$ under this morphism. Let $f_1, f_2\in \Gamma(X,
{\mathscr O}_X)$. We have
$${\mathscr L}_\psi(f_1-f_2)={\mathscr L}_\psi(f_1)\otimes {\mathscr
L}_\psi(f_2)^{-1}.$$ Moreover, we have ${\mathscr L}_\psi(f_1)\cong
{\mathscr L}_\psi(f_2)$ if and only if $f_1-f_2=g^p-g$ for some
$g\in \Gamma(X,{\mathscr O}_X)$. Here for the ``only if" part, we
need the assumption that $\psi:{\mathbb F}_p\to \overline {\mathbb
Q}_l^\ast$ is nontrivial, which implies that $\psi$ is injective.
These facts will be used throughout this paper. (Confer \cite{SGA 4
1/2} Sommes trig. 1.2-1.8.)

\medskip
For any positive integer $N$ prime to $p$, denote by $[N]$ the
$k$-morphism  $$[N]: {\mathbb G}_{m,k}\to {\mathbb G}_{m,k},\quad
t\mapsto t^N.$$ Let $I(k)$ be the set of positive integers $N$ prime
to $p$ such that $k$ contains a primitive $N$-th root of unity. For
any $N\in I(k)$, $[N]$ is a finite Galois etale covering space. Let
$$\bbmu_N(k)=\{\mu\in k|\mu^N=1\}.$$ Then $[N]$ defines a $\bbmu_N(k)$-torsor
$$1\to \bbmu_N(k)\to {\mathbb G}_{m,k}\stackrel {[N]}\to {\mathbb
G}_{m,k}\to 1.$$ When $N$ runs over $I(k)$, we get an inverse system
of extensions of ${\mathbb G}_{m,k}$ by $\bbmu_N(k)$. Passing to the
limit, we get an extension of ${\mathbb G}_{m,k}$ by the profinite
abelian group $\varprojlim_{N\in I(k)} \bbmu_{N}(k).$ Let
$$\rho: \varprojlim_{N\in I(k)} \bbmu_{N}(k)\to \mathrm {GL}(n,
\overline {\mathbb Q}_\ell)$$ be a continuous representation.
Pushing forward the above extension by $\rho^{-1}$, we get a lisse
$\overline{\mathbb Q}_\ell$-sheaf ${\mathscr K}_\rho$ on ${\mathbb
G}_{m,k}$ of rank $n$, which we call the Kummer sheaf associated to
$\rho$. It is tamely ramified at $0$ and at $\infty$. Conversely,
when $k$ is algebraically closed, any lisse $\overline{\mathbb
Q}_\ell$-sheaf on ${\mathbb G}_{m,k}$ tamely ramified at $0$ and at
$\infty$ is obtained this way. This follows from \cite{SGA 1} XIII
2.12. Let $X$ be a scheme over $k$ and let $f$ be a unit in
$\Gamma(X,\mathscr O_X)$. Then $f$ defines a $k$-morphism $X\to
\mathbb G_{m,k}$ so that the induced $k$-algebra homomorphism
$k[t,1/t]\to \Gamma(X,\mathscr O_X)$ maps $t$ to $f$. We denote by
$\mathscr K_{\rho}(f)$ the inverse image of $\mathscr K_{\rho}$
under this morphism.

\medskip
Let $k((t))$ (resp. $k((1/t))$) be the field of formal Laurent
series in the variable $t$ (resp. $1/t$), and let $\eta_0=\mathrm
{Spec}\,k((t))$ (resp. $\eta_\infty=\mathrm {Spec}\,k((1/t))$). We
have canonical morphisms $\eta_0\to {\mathbb G}_m$ and
$\eta_\infty\to {\mathbb G}_m$ defined by the inclusions
$k[t,1/t]\hookrightarrow k((t))$ and $k[t,1/t]\hookrightarrow
k((1/t))$, respectively. Objects on ${\mathbb G}_{m,k}$ can be
restricted to $\eta_0$ and to $\eta_\infty$ through these morphisms.
By \cite{L} 2.2.2.1, any tamely ramified $\overline {\mathbb
Q}_\ell$-sheaf on $\eta_0$ (resp. $\eta_\infty$) is the restriction
of a lisse $\overline{\mathbb Q}_\ell$-sheaf on ${\mathbb G}_{m,k}$
tamely ramified at $0$ and at $\infty$.

\medskip
For any positive integer $r$ prime to $p$, denote also by
$[r]:\eta_0\to \eta_0$ the morphism induced by the $k$-algebra
homomorphism
$$k((t))\to k((t)), \quad t\mapsto t^r.$$ Any formal Laurent series
$$\alpha(t)=\frac{a_{-s}}{t^s}+\frac{a_{-(s-1)}}
{t^{s-1}}+\cdots+ \frac{a_{-1}}{t}+a_0+a_1{t}+\cdots$$ in $k((t))$
can be regarded as a global section of the structure sheaf of
$\eta_0$. In this paper, we calculate the local Fourier
transformation
$${\mathfrak F}^{(0,\infty')}\bigl([r]_\ast\bigl({\mathscr L}_\psi(\alpha(t))\otimes {\mathscr
K}\bigr)\bigr)$$ for any Laurent series $\alpha(t)$ and any tamely
ramified $\overline {\mathbb Q}_\ell$-sheaf ${\mathscr K}$ on
$\eta_0$. Write $$\alpha(t)=\alpha_1(t)+\alpha_2(t)$$ such that
$\alpha_1(t)$ is the polar part of $\alpha(t)$ and $\alpha_2(t)$ is
a formal power series of $t$. Then we have
$$[r]_\ast\bigl({\mathscr L}_\psi(\alpha(t))\otimes {\mathscr
K}\bigr)\cong [r]_\ast\bigl({\mathscr L}_\psi(\alpha_1(t))\otimes
\big({\mathscr L}_\psi(\alpha_2(t))\otimes {\mathscr K}\big)\bigr)$$
and ${\mathscr L}_\psi(\alpha_2(t))\otimes {\mathscr K}$ is again
tamely ramified. So to calculate ${\mathscr
F}^{(0,\infty')}\bigl([r]_\ast\bigl({\mathscr
L}_\psi(\alpha(t))\otimes {\mathscr K}\bigr)\bigr)$ for all Laurent
series $\alpha(t)$ and all tamely ramified sheaf $\mathscr K$, we
may assume $\alpha(t)$ is a polynomial of $1/t$.

Similarly, let $[r]:\eta_\infty\to \eta_\infty$ be the morphism
induced by the $k$-algebra homomorphism
$$k((1/t))\to k((1/t)), \quad t\mapsto t^r.$$ Regard a formal Laurent series
$$\alpha(1/t)={a_{-s}}{t^s}+{a_{-(s-1)}}
t^{s-1}+\cdots+ {a_{-1}}{t}+a_0+\frac{a_1}{t}+\cdots$$ in $k((1/t))$
as a global section of the structure sheaf of $\eta_\infty$. We also
calculate the local Fourier transformations
$${\mathfrak F}^{(\infty,0')}\bigl([r]_\ast\bigl({\mathscr L}_\psi(\alpha(1/t))\otimes {\mathscr
K}\bigr)\bigr) \hbox { and } {\mathscr
F}^{(\infty,\infty')}\bigl([r]_\ast\bigl({\mathscr
L}_\psi(\alpha(1/t))\otimes {\mathscr K}\bigr)\bigr)$$ for any
Laurent series $\alpha(1/t)$ and any tamely ramified $\overline
{\mathbb Q}_\ell$-sheaf ${\mathscr K}$ on $\eta_\infty$. We refer
the reader to \cite{L} for the definitions and properties of local
Fourier transformations.

\medskip
Before presenting the main theorems of this paper, let's recall some
facts from the classical analysis. Consider a smooth concave
function $f(t)$ with $\frac{\mathrm d^2}{\mathrm dt^2}(f(t))<0$. The
Legendre transformation of $f(t)$ is defined as follows: First
consider the function
$$G(t,t')=f(t)+tt'.$$ Require the partial derivative $\frac{\partial
G}{\partial t}$ vanishes:
$$\frac{\mathrm d}{\mathrm dt}(f(t))+t'=0.$$
Solve this equation for $t$ as a function of $t'$ and write the
solution as $t=\mu(t')$. Note that for each fixed $t'$, the function
$G(t,t')$ reaches its greatest value at $t=\mu(t')$. The Legendre
transformation $L_f$ of $f$ is defined to be
$$L_f(t')= G(\mu(t'),t')=f(\mu(t'))+\mu(t')t'.$$
(Classically the Legendre transformation is defined for convex
functions $f(t)$ and $G(t,t')$ is taken to be $-f(t)+tt'$. Confer
\cite{A} Chap. 3, \S 14.) According to Laplace, (confer \cite{C}
Chapter 5, especially \S 17 and \S 20), under suitable conditions on
$\phi(t,t')$, the major contribution to the value of the integral
$$\int_{-\infty}^{\infty}\phi(t,t') e^{f(t)+tt'}\mathrm dt$$ for large real $t'$ comes
from neighborhoods of those points where $G(t,t')=f(t)+tt'$ attains
its greatest values, that is, those points on the curve $t=\mu(t')$.
Moreover, the dominant term in the asymptotic expansion of
$\int_{-\infty}^{\infty}\phi(t,t') e^{f(t)+tt'}\mathrm dt$ is
$$
\phi(\mu(t'),t')e^{L_f(t')}\sqrt{\frac{-\pi}{2\frac{\mathrm
d^2f}{\mathrm dt^2}(\mu(t'))}}.$$

Go back to the characteristic $p$ case, and let
$$f(t)=\frac{a_{-s}}{(\sqrt[r] {t})^s}+\frac{a_{s-1}}{(\sqrt[r]{t})^{s-1}}+\cdots
+\frac{a_{-1}}{\sqrt[r]{t}}+a_0+\cdots$$ be a Laurent series in the
variable $\sqrt[r]{t}$ with $a_{-s}\not=0$. Similar to the above
discussion, we define the Legendre transformation $L_f(t')$ of
$f(t)$ by the system of equations
\[\left\{ \begin{array}{l}
f(t)+tt'=L_f(t'),\cr \frac{\mathrm d}{\mathrm dt}(f(t))+t'=0.
\end{array}\right.\]
Assume $k$ is algebraically closed. Then $L_f(t')$ is a Laurent
series in the variable $\frac{1}{\sqrt[r+s]{t'}}$ of the form
$$L_f(t')=b_{-s}(\sqrt[r+s] {t'})^s +b_{-(s-1)} (\sqrt[r+s]
{t'})^{s-1} +\cdots +b_0+ \frac{b_1}{\sqrt[r+s] {t'}}+\cdots$$ with
$b_{-s}\not=0$. Set
$$
\alpha(t)=\sum_{i=-s}^\infty {a_i}{t^i},\;
\beta(1/t')=\sum_{i=-s}^\infty b_i(1/t')^i,
$$
and let ${\mathscr L}(f)=[r]_\ast {\mathscr L}_\psi(\alpha(t))$ and
${\mathscr L}(L_f)=[r+s]_\ast{\mathscr L}_\psi(\beta(1/t'))$.
Probably motivated by the asymptotic expansion discussed in the
previous paragraph, Laumon and Malgrange conjecture that when $p\gg
r,s$, we have
$${\mathfrak F}^{(0,\infty')}({\mathscr L}(f))\cong {\mathscr
L}(L_f).$$ They also made similar conjectures for ${\mathfrak
F}^{(\infty,0')}$ and ${\mathfrak F}^{(\infty,\infty')}$ (\cite{L}
2.6.3). The main goal of this paper is to prove these conjectures
(which need to be slightly corrected). Our main theorems are as
follows:

\begin{theorem} Let $\gamma(t)=t^r$, let
$\alpha(t)\in k[t,1/t]$, let
\begin{eqnarray*}
&&\delta(t)
=-\frac{1}{rt^{r-1}}\frac{\mathrm d}{dt}(\alpha(t)),\\
&&\beta(t')=\alpha(t')+\gamma(t')\delta(t'),
\end{eqnarray*}
and let $\mathscr K$ be a rank $1$ lisse $\overline{\mathbb
Q}_\ell$-sheaf on ${\mathbb G}_{m,k}$ tamely ramified at $0$ and at
$\infty$. Denote by $\gamma,\delta:\mathbb {\mathbb G}_{m,k}\to
{\mathbb A}_k^1$ the $k$-morphisms defined by $\gamma(t)$ and
$\delta(t')$, respectively.

(i) Suppose $\alpha(t)$ is of the form
$$\alpha(t)=\frac{a_{-s}}{t^s}+\frac{a_{-(s-1)}}{t^{s-1}}+\cdots+\frac{a_{-1}}{t}$$
with $a_{-s}\not=0$, and suppose $r,s\geq 1$, $s<p$, and $p$ is
relatively prime to $2$, $r$, $s$ and $r+s$. We have
$${\mathfrak F}^{(0,\infty')}\bigg(\Big(\gamma_\ast\bigl({\mathscr L}_\psi
(\alpha(t))\otimes {\mathscr K}\bigr)\Big)|_{\eta_0}\bigg)\cong
\bigg(\delta_\ast \Big({\mathscr L}_\psi(\beta(t'))\otimes {\mathscr
K}\otimes {\mathscr
K}_{\chi_2}\Big(\frac{1}{2}s(r+s)a_{-s}t'^s\Big)\otimes
G(\chi_2,\psi)\Big)\bigg)|_{\eta_{\infty'}},$$ where ${\mathscr
K}_{\chi_2}$ is the Kummer sheaf associated to the (unique)
nontrivial character $\chi_2:\bbmu_2(k)\to \overline {\mathbb
Q}_\ell^\ast$ of order $2$, and $G(\chi_2,\psi)$ is the unramified
sheaf associated to the Gauss sum, that is, the $\overline {\mathbb
Q}_\ell$-vector space of dimension $1$ with continuous
$\mathrm{Gal}(\bar k/k)$-action defined by
$$G(\chi_2,\psi)=H_c^1(\mathbb G_{m,k}\otimes_k \bar k,\mathscr
K_{\chi_2}\otimes\mathscr L_\psi).$$

(ii) Suppose $\alpha(t)$ is of the form
$$\alpha(t)=a_st^s+a_{s-1}t^{s-1}+\cdots+a_1t$$ with $a_s\not=0$, and
suppose $r,s\geq 1$ and $p$ is relatively prime to $r$ and $s$. If
$s\leq r$, then
$${\mathfrak F}^{(\infty,\infty')}\bigg(\Big(\gamma_\ast\bigl({\mathscr L}_\psi
(\alpha(t))\otimes {\mathscr K}\bigr)\Big)|_{\eta_\infty}\bigg)=0.$$
Suppose $s>r$ and suppose furthermore that $s<p$ and $p$ is
relatively prime to $2$, $r$, $s$ and $s-r$. We have
$${\mathfrak F}^{(\infty,\infty')}\bigg(\Big(\gamma_\ast\bigl({\mathscr L}_\psi
(\alpha(t))\otimes {\mathscr K}\bigr)\Big)|_{\eta_\infty}\bigg)\cong
\bigg(\delta_\ast \Big({\mathscr L}_\psi(\beta(t'))\otimes {\mathscr
K}\otimes {\mathscr
K}_{\chi_2}\Big(\frac{1}{2}s(s-r)a_{s}t'^s\Big)\otimes
G(\chi_2,\psi)\Big)\bigg)|_{\eta_{\infty'}}.$$

(iii) Suppose $\alpha(t)$ is of the form
$$\alpha(t)=a_st^s+a_{s-1}t^{s-1}+\cdots+a_1t$$ with $a_s\not=0$, and
suppose $r,s\geq 1$ and $p$ is relatively prime to $r$ and $s$. If
$s\geq r$, then
$${\mathfrak F}^{(\infty,0')}\bigg(\Big(\gamma_\ast\bigl({\mathscr L}_\psi
(\alpha(t))\otimes {\mathscr K}\bigr)\Big)|_{\eta_\infty}\bigg)=0.$$
Suppose $s<r$ and suppose furthermore that $s<p$ and $p$ is
relatively prime to $2$, $r$, $s$ and $s-r$. We have
$${\mathfrak F}^{(\infty,0')}\bigg(\Big(\gamma_\ast\bigl({\mathscr L}_\psi
(\alpha(t))\otimes {\mathscr K}\bigr)\Big)|_{\eta_\infty}\bigg)\cong
\bigg(\delta_\ast \Big({\mathscr L}_\psi(\beta(t'))\otimes {\mathscr
K}\otimes {\mathscr
K}_{\chi_2}\Big(\frac{1}{2}s(s-r)a_{s}t'^s\Big)\otimes
G(\chi_2,\psi)\Big)\bigg)|_{\eta_{0'}}.$$
\end{theorem}

In the case where $k$ is algebraically closed, we have the following
theorems:

\begin{theorem}
Suppose $k$ is algebraically closed, $r,s\geq 1$, $s<p$, and $p$ is
relatively prime to $2$, $r$, $s$ and $r+s$. Let
$$\alpha(t)=\sum_{i=-s}^\infty \frac{a_i}{t^i}$$ be a formal Laurent series in
$k((1/t))$ with $a_{-s}\not=0$. Consider the system of equations
\[\left\{ \begin{array}{l}
\alpha(t)+t^rt'^{r+s}=\beta(1/t'),\cr \frac{\mathrm d}{\mathrm
dt}(\alpha(t))+rt^{r-1}t'^{r+s}=0.
\end{array}\right.\]
Using the second equation, we find an expression of $t$ in terms of
$t'$. We then substitute this expression into the first equation to
get $\beta(1/t')$, which is a formal Laurent series in $k((1/t'))$
of the form
$$\beta(1/t')=\sum_{i=-s}^\infty b_i(1/t')^i$$ with
$b_{-s}\not=0.$ For any lisse $\overline{\mathbb Q}_\ell$-sheaf
${\mathscr K}$ on ${\mathbb G}_{m,k}$ tamely ramified at $0$ and at
$\infty$, we have
$${\mathfrak F}^{(0,\infty')}\biggl([r]_\ast\bigl({\mathscr L}_\psi(\alpha(t))\otimes {\mathscr
K}|_{\eta_0}\bigr)\biggr)\cong [r+s]_\ast \bigl({\mathscr
L}_\psi(\beta(1/t'))\otimes (\mathrm {inv}^\ast {\mathscr
K})|_{\eta_{\infty'}}\otimes ([s]^\ast {\mathscr
K}_{\chi_2})|_{\eta_{\infty'}}\bigr),$$ where $\mathrm
{inv}:{\mathbb G}_{m,k}\to {\mathbb G}_{m,k}$ is the morphism
$t'\mapsto \frac{1}{t'}$, and ${\mathscr K}_{\chi_2}$ is the Kummer
sheaf associated to the (unique) nontrivial character
$\chi_2:\bbmu_2(k)\to \overline {\mathbb Q}_\ell^\ast$ of order $2$.
\end{theorem}

\begin{theorem} Suppose $k$ is algebraically closed,
$r,s\geq 1$ and $p$ is relatively prime to $r$ and
$s$. Let
$$\alpha(1/t)=\sum_{i=-s}^\infty {a_i}{(1/t)^i}$$ be a formal Laurent series in
$k((1/t))$ with $a_{-s}\not=0$, and let ${\mathscr K}$ be a lisse
$\overline{\mathbb Q}_\ell$-sheaf on ${\mathbb G}_{m,k}$ tamely
ramified at $0$ and at $\infty$. If $s\leq r$, then
$${\mathfrak F}^{(\infty,\infty')}\biggl([r]_\ast\bigl({\mathscr
L}_\psi(\alpha(1/t))\otimes {\mathscr
K}|_{\eta_\infty}\bigr)\biggr)=0.$$ Suppose $s>r$ and suppose
furthermore that $s<p$ and $p$ is relatively prime to $2$, $r$, $s$
and $s-r$. Define $\beta(1/t')$ by the system of equations
\[\left\{ \begin{array}{l}
\alpha(1/t)+t^rt'^{s-r}=\beta(1/t'),\cr \frac{\mathrm d}{\mathrm
dt}(\alpha(1/t))+rt^{r-1}t'^{s-r}=0.
\end{array}\right.\]
It is a formal Laurent series in $k((1/t'))$ of the form
$$\beta(1/t')=\sum_{i=-s}^\infty b_i(1/t')^i$$ with
$b_{-s}\not=0.$  We have
$${\mathfrak F}^{(\infty,\infty')}\biggl([r]_\ast\bigl({\mathscr L}_\psi(\alpha(1/t))\otimes
{\mathscr K}|_{\eta_\infty}\bigr)\biggr)\cong [s-r]_\ast
\bigl({\mathscr L}_\psi(\beta(1/t'))\otimes {\mathscr
K}|_{\eta_{\infty'}}\otimes ([s]^\ast {\mathscr
K}_{\chi_2})|_{\eta_{\infty'}}\bigr),$$ where ${\mathscr
K}_{\chi_2}$ is the Kummer sheaf associated to the (unique)
nontrivial character $\chi_2:\bbmu_2(k)\to \overline {\mathbb
Q}_\ell^\ast$ of order $2$.
\end{theorem}

\begin{theorem} Suppose $k$ is algebraically closed,
$r,s\geq 1$ and $p$ is relatively prime to
$r$ and $s$. Let
$$\alpha(1/t)=\sum_{i=-s}^\infty {a_i}(1/t)^i$$ be a formal Laurent series in
$k((1/t))$ with $a_{-s}\not=0$, and let ${\mathscr K}$ be a lisse
$\overline{\mathbb Q}_\ell$-sheaf on ${\mathbb G}_{m,k}$ tamely
ramified at $0$ and at $\infty$. If $s\geq r$, then
$${\mathfrak F}^{(\infty,0')}\biggl([r]_\ast\bigl({\mathscr
L}_\psi(\alpha(1/t))\otimes {\mathscr
K}|_{\eta_{\infty}}\bigr)\biggr)=0.$$ Suppose $s<r$ and suppose
furthermore that $s<p$ and $p$ is relatively prime to $2$, $r$, $s$
and $r-s$. Define $\beta(t')$ by the system of equations
\[\left\{ \begin{array}{l}
\alpha(1/t)+t^rt'^{r-s}=\beta(t'),\cr \frac{\mathrm d}{\mathrm
dt}(\alpha(1/t))+rt^{r-1}t'^{r-s}=0.
\end{array}\right.\]
It is a formal Laurent series in $k((t'))$ of the form
$$\beta(t')=\sum_{i=-s}^\infty b_it'^i$$ with
$b_{-s}\not=0.$ We have
$${\mathfrak F}^{(\infty,0')}\biggl([r]_\ast\bigl({\mathscr L}_\psi(\alpha(1/t))\otimes
{\mathscr K}|_{\eta_{\infty}}\bigr)\biggr)\cong [r-s]_\ast
\bigl({\mathscr L}_\psi(\beta(t'))\otimes (\mathrm {inv}^\ast
{\mathscr K})|_{\eta_{0'}}\otimes ([s]^\ast {\mathscr
K}_{\chi_2})|_{\eta_{0'}}\bigr),$$ where $\mathrm {inv}:{\mathbb
G}_{m,k}\to {\mathbb G}_{m,k}$ is the morphism $t'\mapsto
\frac{1}{t'}$, and ${\mathscr K}_{\chi_2}$ is the Kummer sheaf
associated to the (unique) nontrivial character
$\chi_2:\bbmu_2(k)\to \overline {\mathbb Q}_\ell^\ast$ of order $2$.
\end{theorem}

When ${\mathscr K}$ is trivial, Theorems 0.2-0.4 are conjectured by
Laumon and Malgrange except that the term $[s]^\ast \mathscr
K_{\chi_2}$ is missing in their conjecture. Note that in these
theorems, the system of equations relating $\alpha$ and $\beta$ is
equivalent to the system of equations defining the Legendre
transformation. For example, in Theorem 0.2,
$\beta(1/\sqrt[r+s]{t'})$ is the Legendre transformation of
$\alpha(\sqrt[r]{t})$.

Note that we take $\mathscr K$ to be of rank $1$ if $k$ is not
assumed to be algebraically closed, whereas $\mathscr K$ has
arbitrary rank if $k$ is algebraically closed. This is for the
technical reason that when $k$ is algebraically closed, we know the
tame fundamental group of $\mathbb G_{m,k}$ is isomorphic to
$\varprojlim_{(N,p)=1}\bbmu_N(k)$ so that we have a structure
theorem (Lemma 2.5) for lisse $\overline {\mathbb Q}_\ell$-sheaves
on $\mathbb G_{m,k}$ tamely ramified at $0$ and $\infty$. This is
not available when $k$ is not algebraically closed.

In \cite{R}, Ramero constructs the Fourier transformations for a
class of \'etale analytic sheaves on the analytification of the
affine line over a field of characteristic $0$, complete with
respect to a non-Archimedean metric. In order for the Fourier
transformations to be constructible sheaves, the analytic sheaves
are assumed to have meromorphic ramification (\cite{R} 8.3.6).
Ramero asks whether the Fourier transformation of an analytic sheaf
with meromorphic ramification still has meromorphic ramification,
and this problem can be translated into a problem concerning the
local Fourier transformation of a meromorphic representation
(\cite{R} 8.7). An analytic sheaf with meromorphic ramification on
the germ of a punctured disc at $0$ is nothing but an analytic
analogue of the $\overline {\mathbb Q}_\ell$-sheaf
$[r]_\ast\bigl({\mathscr L}_\psi(\alpha(t))\otimes {\mathscr
K}|_{\eta_0}\bigr)$ studied in this paper (confer \cite{R} \S 5 and
8.1.1). The method used in this paper can also be used to calculate
the local Fourier transformations of analytic sheaves with
meromorphic ramification. The resulting formulas are the same as
those in Theorems 0.2-0.4. In particular, we can give a positive
answer to Ramero's question asking whether the Fourier
transformation of an analytic sheaf with meromorphic ramification
still has meromorphic ramification. In fact, it was for the purpose
of solving Ramero's problem that I started working on the
Laumon-Malgrange conjecture.

\medskip About the same time when the main results of this paper was
obtained, analogous formulas for the local Fourier transformations
of $D$-modules were proved independently by J. Fang \cite{Fa} and C.
Sabbah \cite{Sa}. By the refined Turrittin-Levelt Theorem (\cite{Sa}
Corollary 3.3), when $k$ is an algebraically closed field of
characteristic $0$, any finite dimensional $k((t))$-vector space $M$
with a connection is a direct sum of the $D$-module analogues of the
$\overline {\mathbb Q}_\ell$-sheaves $[r]_\ast\bigl({\mathscr
L}_\psi(\alpha(t))\otimes {\mathscr K}|_{\eta_0}\bigr)$. So to
calculate the local Fourier transformation of $M$, it suffices to
calculate the local Fourier transformations of the $D$-module
analogues of $[r]_\ast\bigl({\mathscr L}_\psi(\alpha(t))\otimes
{\mathscr K}|_{\eta_0}\bigr)$. Of course, the Turrittin-Levelt
Theorem doesn't hold for $\overline {\mathbb Q}_\ell$-sheaves on
$\eta_0$ in the characteristic $p$ case. But Deligne (\cite{D1})
suggested to me the following analogue of the Turrittin-Levelt
Theorem in characteristic $p$.

\begin{proposition} Suppose $k$ is algebraically closed.
Let $\rho: I\to \mathrm {GL}(V)$ be an irreducible $\overline
{\mathbb Q}_\ell$-representation of $I={\rm Gal}\Big(\overline
{k((t))}/k((t))\Big)$, where $\overline {k((t))}$ is a separable
closure of $k((t))$. Suppose the following conditions hold:

(a) $\rho(P^p[P,P])=1$, where $P=\mathrm {Gal}\left(\overline
{k((t))}/\bigcup_{(N,p)=1} k(\sqrt[N]{t})\right)$ is the wild
inertia subgroup.

(b) $\rho(I)$ is finite.

(c) Let $s$ be the Swan conductor of $\rho$. We have $s<p$.

\noindent Then there exist a character
$\chi:\varprojlim\limits_{(N,p)=1} \bbmu_N(k) \to \overline{\mathbb
Q}_\ell^\ast$ of finite order, a Laurent series $\alpha(t)\in
k((t))$ of the form
$$\alpha(t)=\frac{a_{-s}}{t^s}+\frac{a_{-(s-1)}}{t^{s-1}}+\cdots+\frac{a_{-1}}{t}$$
with $a_{-s}\not=0$, and a positive integer $r$ prime to $p$ such
that $V$ is isomorphic to the representation corresponding to the
lisse $\overline{\mathbb Q}_\ell$-sheaf $[r]_\ast(\mathscr
L_\psi(\alpha(t))\otimes \mathscr K_\chi)$ on $\eta_0=\mathrm
{Spec}\,k((t))$.
\end{proposition}

In the above proposition, condition (a) is necessary for the
representation $V$ to be isomorphic to the representation
corresponding to a sheaf of the form $[r]_\ast(\mathscr
L_\psi(\alpha(t))\otimes \mathscr K_\chi)$ on $\eta_0=\mathrm
{Spec}\,k((t))$. Condition (b) holds for most interesting
representations. For example, if we know $\rho$ is irreducible and
is quasi-unipotent, that is, there exists a closed subgroup $I_1$ of
$I$ of finite index such that $\rho(\sigma)$ is unipotent for any
$\sigma\in I_1$, then (b) holds. We know representations coming from
arithmetic and geometry are often quasi-unipotent. See for example,
the appendix of \cite{ST}, and \cite{SGA 7} I 1.3. Condition (c)
should not be essential. But we need it in our proof for some
algebraic argument to work.

\medskip In the following, ${\mathbb Z}$-schemes and ${\mathbb
Z}$-morphisms are assumed to be separated and of finite type. Taking
fibers at a prime $p$ of ${\mathbb Z}$-schemes and ${\mathbb
Z}$-morphisms is indicated by putting the subscript $p$. In
\cite{K2} pages 297-301, in order to have a general framework to
study the variation with $p$ of exponential sums on schemes over
${\mathbb Z}$, Katz proposes to study the family
$$\{(Y_p, {\mathscr G}_p)| p \hbox { is prime}\},$$ where $Y$ is a
scheme over ${\mathbb Z}$, and for each $p$ prime to $\ell$,
${\mathscr G}_p$ is a $\overline {\mathbb Q}_\ell$-sheaf on the
fiber $Y_p$ so that there exist a scheme $X$ over $\mathbb Z$, an
object $K$ in the category $D_c^b(X,\overline{\mathbb Q}_\ell)$ of
complexes of $\overline{\mathbb Q}_\ell$-sheaves on $X$ constructed
in \cite{D2} 1.1, and two $\mathbb Z$-morphisms $f:X\to Y$ and
$g:X\to \mathbb A_{\mathbb Z}^1$ with the property
$${\mathscr G}_p\cong R^if_{p!}(K|_{X_p}\otimes g_p^\ast{\mathscr
L}_{\psi_p})$$ for some $i$ and all $p$ prime to $\ell$, where
${\mathscr L}_{\psi_p}$ is the Artin-Schreier sheaf defined by the
standard additive character
$$\psi_p:{\mathbb F}_p\to \overline {\mathbb Q}_l^\ast,\quad \psi_p(t)=e^{\frac{2\pi
i t}{p}}.$$ Let
$$Y\stackrel {\mathrm {pr}_1}\leftarrow Y\times {\mathbb A}_{\mathbb Z}^1
\stackrel{\mathrm {pr}_2} \to {\mathbb A}_{\mathbb Z}^1$$ be the
projections, let $h:X\to Y\times {\mathbb A}_{\mathbb Z}^1$ be the
morphism defined by $$\mathrm {pr}_1\circ h=f, \; \mathrm
{pr}_2\circ h=g,$$ and let $L= Rh_!K.$ Then we have
$${\mathscr G}_p\cong R^i \mathrm {pr}_{1p!}
(L|_{(Y\times {\mathbb A}_{\mathbb Z}^1)_p}\otimes \mathrm
{pr}_{2p}^\ast {\mathscr L}_{\psi_p}).$$ So in the definition of the
family $\{(Y_p,\mathscr G_p)\}$, we may require $X=Y\times {\mathbb
A}_{\mathbb Z}^1$, $f=\mathrm {pr}_1$ and $g=\mathrm {pr}_2$.

\begin{conjecture} With the above notations, suppose $Y$ is an algebraic
curve over ${\mathbb Z}$. Let $\overline {\mathbb F}_p$ be an
algebraic closure of ${\mathbb F}_p$. Then for sufficiently large
$p$, the inverse image of ${\mathscr G}_p$ by any dominant morphism
$$\mathrm {Spec}\, \overline {\mathbb F}_p((t)) \to Y_p$$ can be written as a direct sum
of objects of the form
$$[r]_\ast ({\mathscr L}_{\psi_p}(\alpha(t))\otimes {\mathscr K}),$$
where $r$ is relatively prime to $p$, $\alpha(t)\in \overline
{\mathbb F}_p((t))$, and ${\mathscr K}$ is a tame $\overline
{\mathbb Q}_\ell$-sheaf on $\mathrm {Spec}\, \overline {\mathbb
F}_p((t))$.
\end{conjecture}

We expect the method used in this paper is helpful to prove the
above conjecture. In any case, Theorems 0.1-0.4 can be used to
calculate the local monodromy of sheaves arising from the study of
exponential sums. Indeed, as a direct application of our theorems,
we calculate the local monodromy of hypergeometric $\overline
{\mathbb Q}_\ell$-sheaves introduced by Katz (\cite{K1}).

\medskip Let
$$\lambda_1,\ldots, \lambda_n,\rho_1,\ldots, \rho_m:{\mathbb
F}_q^\ast\to \overline {\mathbb Q}_\ell^\ast$$ be multiplicative
characters. For any extension ${\mathbb F}_{q^k}$ of ${\mathbb F}_q$
of degree $k$, and any $t\not=0$ in ${\mathbb F}_{q^k}$, we define
the hypergeometric sum by
\begin{eqnarray*}
&&\mathrm {Hyp}(\psi; \lambda_1,\ldots, \lambda_n;\rho_1,\ldots,
\rho_m)({\mathbb F}_{q^k},t)\\
&=&\sum\limits_{x_1,\ldots, x_n,y_1,\ldots, y_m\in {\mathbb
F}_{q^k}, x_1\cdots x_n=ty_1\cdots y_m}
\psi\biggl(\mathrm{Tr}_{{\mathbb F}_{q^k}/{\mathbb
F}_p}(\sum_{i=1}^nx_i-\sum_{j=1}^m
y_j)\biggr)\\
&&\qquad\qquad \qquad\qquad\times \prod_{i=1}^n
\lambda_i\biggl({\mathrm N}_{{\mathbb F}_{q^k}/{\mathbb
F}_q}(x_i)\biggr)\prod_{j=1}^m\rho_j\biggl({\mathrm N}_{{\mathbb
F}_{q^k}/{\mathbb F}_q}\left(\frac{1}{y_j}\right)\biggr).
\end{eqnarray*}
In the case where $m=0$, we get the Kloosterman sum
$$\mathrm {Kl}_n(\psi; \lambda_1,\ldots, \lambda_n)({\mathbb F}_{q^k},t)\\
=\sum_{x_1,\ldots, x_n\in {\mathbb F}_{q^k}, x_1\cdots x_n=t}
\psi\biggl(\mathrm{Tr}_{{\mathbb F}_{q^k}/{\mathbb
F}_p}(\sum_{i=1}^nx_i)\biggr)\prod_{i=1}^n \lambda_i\biggl({\mathrm
N}_{{\mathbb F}_{q^k}/{\mathbb F}_q}(x_i)\biggr).$$ In \cite{K1}
8.2, Katz constructs an object $\mathrm {Hyp}(!,\psi;
\lambda_1,\ldots, \lambda_n;\rho_1,\ldots, \rho_m)$ in the
triangulated category $D_c^b(\mathbb G_{m,{\mathbb F}_q},\overline
{\mathbb Q}_\ell)$ so that for any ${\mathbb F}_{q^k}$-rational
point $t$ of ${\mathbb G}_{m,{\mathbb F}_q}$, we have
$$\mathrm {Tr}(F_t, (\mathrm {Hyp}(!,\psi; \lambda_1,\ldots,
\lambda_n;\rho_1,\ldots, \rho_m))_{\bar t})=(-1)^{n+m} \mathrm
{Hyp}(\psi; \lambda_1,\ldots, \lambda_n;\rho_1,\ldots,
\rho_m)({\mathbb F}_{q^k},t),$$ where $F_t$ is the geometric
Frobenius element at $t$. In fact, Katz constructs $\mathrm
{Hyp}(!,\psi; \lambda_1,\ldots, \lambda_n;\rho_1,\ldots, \rho_m)$ on
${\mathbb G}_{m,k}$ for any field $k$ of characteristic $p$ and any
characters
$$\lambda_1,\ldots, \lambda_n,\rho_1,\ldots, \rho_m:
\varprojlim_{N\in I(k)}\bbmu_N(k)\to \overline {\mathbb
Q}_\ell^\ast.$$

Assume $k$ is algebraically closed. Fix an element
$$\zeta=(\zeta_N)\in \varprojlim_{(N,p)=1}\bbmu_N(k)$$ so that
$\zeta_N$ is a primitive $N$-th root of unity for each $N$ prime to
$p$. Note that $\zeta$ is a topological generator of
$\varprojlim\limits_{(N,p)=1}\bbmu_N(k)$. Let $U(n)$ be an
$n$-dimensional vector space $\overline {\mathbb Q}_\ell^n$ on which
$\zeta$ acts through the unipotent $(n\times n)$-matrix with a
single Jordan block
$$\left(\begin{array}{cccc}
1&1&&\\
&1&\ddots& \\
&&\ddots&1\\
&&&1
\end{array}
\right).$$ This action can be extended to a continuous action of
$\varprojlim\limits_{(N,p)=1}\bbmu_N(k)$ on $U(n)$. Denote also by
$U(n)$ the corresponding lisse $\overline{\mathbb Q}_\ell$-sheaf on
${\mathbb G}_{m,k}$ tamely ramified at $0$ and at $\infty$. For
sufficiently large $p$, the following proposition is more precise
than \cite{K1} 8.4.2 (6) and (7).

\begin{proposition} Let $k$ be an algebraically closed field of
characteristic $p$, and let
$$\lambda_1,\ldots, \lambda_n,\rho_1,\ldots, \rho_m:
\varprojlim_{(N,p)=1}\bbmu_N(k)\to \overline {\mathbb Q}_\ell^\ast$$
be characters. For any character $\lambda:
\varprojlim\limits_{(N,p)=1}\bbmu_N(k)\to \overline {\mathbb
Q}_\ell^\ast$, let $\mathrm {mult}_0(\lambda)$ be the number of
times of $\lambda$ occurring in $\lambda_1,\ldots, \lambda_n$, and
let $\mathrm {mult}_\infty(\lambda)$ be the number of times of
$\lambda$ occurring in $\rho_1,\ldots, \rho_m$. Suppose the family
of characters $\lambda_1,\ldots, \lambda_n$ is disjoint from the
family of characters $\rho_1,\ldots, \rho_m$, and suppose $p$ is
relatively prime to $2, 3, \ldots, \mathrm{max}(m,n)$.

(i) If $n>m$, then $\mathrm {Hyp}(!,\psi; \lambda_1,\ldots,
\lambda_n;\rho_1,\ldots, \rho_m)[-1]$ is an irreducible lisse sheaf
of rank $n$ on ${\mathbb G}_{m,k}$, and
\begin{eqnarray*}
&& \bigl(\mathrm {Hyp}(!,\psi; \lambda_1,\ldots,
\lambda_n;\rho_1,\ldots, \rho_m)[-1]\bigr)|_{\eta_0}\cong
\bigoplus_{\lambda} \bigl({\mathscr K}_\lambda\otimes U(\mathrm
{mult}_0(\lambda))\bigr)|_{\eta_0},\\
&&\bigl(\mathrm {Hyp}(!,\psi; \lambda_1,\ldots,
\lambda_n;\rho_1,\ldots, \rho_m)[-1]\bigr)|_{\eta_\infty}\\&\cong
&[n-m]_\ast\bigl({\mathscr L}_\psi((n-m)t)\otimes {\mathscr
K}_{\lambda_1\cdots\lambda_n\rho_1^{-1}\cdots
\rho_m^{-1}}\otimes{\mathscr
K}_{{\chi_2}^{n+m-1}}\bigr)|_{\eta_\infty}\bigoplus \bigoplus_{\rho}
\bigl({\mathscr K}_\rho\otimes U(\mathrm
{mult}_\infty(\rho))\bigr)|_{\eta_\infty}.
\end{eqnarray*}

(ii) If $n<m$, then $\mathrm {Hyp}(!,\psi; \lambda_1,\ldots,
\lambda_n;\rho_1,\ldots, \rho_m)[-1]$ is an irreducible lisse sheaf
of rank $m$ on ${\mathbb G}_{m,k}$, and
\begin{eqnarray*}
&&\bigl(\mathrm {Hyp}(!,\psi; \lambda_1,\ldots,
\lambda_n;\rho_1,\ldots, \rho_m)[-1]\bigr)|_{\eta_\infty}\cong
\bigoplus_{\rho} \bigl({\mathscr K}_\rho\otimes U(\mathrm
{mult}_\infty(\rho))\bigr)|_{\eta_\infty},\\
&&\bigl(\mathrm {Hyp}(!,\psi; \lambda_1,\ldots,
\lambda_n;\rho_1,\ldots, \rho_m)[-1]\bigr)|_{\eta_0}\\&\cong
&[m-n]_\ast\left({\mathscr L}_\psi\left(-\frac{m-n}{t}\right)\otimes
{\mathscr K}_{\lambda_1^{-1}\cdots\lambda_n^{-1}\rho_1\cdots
\rho_m}\otimes{\mathscr K}_{{\chi_2}^{n+m-1}}\right)|_{\eta_0}
\bigoplus \bigoplus_{\lambda}\bigl({\mathscr K}_\lambda\otimes
U(\mathrm {mult}_0(\lambda))\bigr)|_{\eta_0}.
\end{eqnarray*}

(iii) If $n=m$, then $\mathrm {Hyp}(!,\psi; \lambda_1,\ldots,
\lambda_n;\rho_1,\ldots, \rho_m)[-1]$ is an irreducible lisse sheaf
of rank $n$ when restricted to ${\mathbb G}_{m,k}-\{1\}$, and we
have
\begin{eqnarray*}
\bigl(\mathrm {Hyp}(!,\psi; \lambda_1,\ldots,
\lambda_n;\rho_1,\ldots, \rho_m)[-1]\bigr)|_{\eta_0}&\cong&
\bigoplus_{\lambda} \bigl({\mathscr K}_\lambda\otimes U(\mathrm
{mult}_0(\lambda))\bigr)|_{\eta_0},\\
\bigl(\mathrm {Hyp}(!,\psi; \lambda_1,\ldots,
\lambda_n;\rho_1,\ldots, \rho_m)[-1]\bigr)|_{\eta_\infty}&\cong&
\bigoplus_{\rho} \bigl({\mathscr K}_\rho\otimes U(\mathrm
{mult}_\infty(\rho))\bigr)|_{\eta_\infty}.
\end{eqnarray*}
Let $k((t-1))$ be the field of Laurent power series in the variable
$t-1$, let $\eta_1=\mathrm {Spec}\, k((t-1))$, and let $$\mathrm
{tran}:\mathbb A_k^1\to \mathbb A_k^1,\quad x\mapsto x-1$$ be the
translation by $1$. We have
$$(\mathrm {Hyp}(!,\psi;\lambda_1,\ldots,
\lambda_n;\rho_1,\ldots,\rho_n)[-1])_{\bar 1}\cong (\mathrm
{Hyp}(!,\psi;\lambda_1,\ldots,
\lambda_n;\rho_1,\ldots,\rho_n)[-1])_{\bar
\eta_1}^{\mathrm{Gal}(\bar\eta_1/\eta_1)},$$ and we have a short
exact sequence
\begin{eqnarray*}
0&\to& (\mathrm {Hyp}(!,\psi;\lambda_1,\ldots,
\lambda_n;\rho_1,\ldots,\rho_n)[-1])_{\bar
\eta_1}^{\mathrm{Gal}(\bar\eta_1/\eta_1)} \\
&\to& (\mathrm {Hyp}(!,\psi;\lambda_1,\ldots,
\lambda_n;\rho_1,\ldots,\rho_n)[-1])_{\bar \eta_1}\\
&\to& (\mathrm{tran}^\ast{\mathscr
K}_{\lambda_1^{-1}\cdots\lambda_n^{-1}\rho_1\cdots \rho_n})|_{\bar
\eta_1} \to 0.
\end{eqnarray*}
\end{proposition}

Taking $m=0$, we get the following.

\begin{proposition} Suppose $k$ is algebraically closed and
$p=\mathrm{char}\,k$ is relatively prime to $2, 3, \ldots, n$. The
Kloosterman sheaf
$$\mathrm {Kl}_n(\psi; \lambda_1,\ldots, \lambda_n)
= \mathrm {Hyp}(!,\psi; \lambda_1,\ldots, \lambda_n;\emptyset)[-1]$$
is an irreducible lisse sheaf of rank $n$ on ${\mathbb G}_{m,k}$ and
\begin{eqnarray*}
\mathrm {Kl}_n(\psi; \lambda_1,\ldots, \lambda_n)|_{\eta_0}&\cong&
\bigoplus_{\lambda} \bigl({\mathscr K}_\lambda\otimes U(\mathrm
{mult}_0(\lambda))\bigr)|_{\eta_0},\\
\mathrm {Kl}_n(\psi; \lambda_1,\ldots,
\lambda_n)|_{\eta_\infty}&\cong& [n]_\ast({\mathscr
L}_\psi(nt)\otimes {\mathscr
K}_{\lambda_1\cdots\lambda_n}\otimes{\mathscr
K}_{{\chi_2}^{n-1}})|_{\eta_\infty}.
\end{eqnarray*}
\end{proposition}

When $\lambda_1,\ldots,\lambda_n$ are trivial, the above assertion
about $\mathrm {Kl}_n(\psi; \lambda_1,\ldots,
\lambda_n)|_{\eta_\infty}$ is Lemma 1.6 in \cite{FW}. However in
\cite{FW}, we only assume $p$ is relatively prime to $n$. The result
in \cite{FW} is based on the work of Katz (\cite{K3} 10.1 and
5.6.2). Katz deduces his results from an identity of
Hasse-Davenport. In \cite{K1} 8.4.2 (8), Katz shows that Proposition
0.7 (iii) holds without the assumption that $p$ is relatively prime
to $1,\ldots, n$. It is an interesting question to determine
explicitly the local monodromy at $0$ and at $\infty$ for $\mathrm
{Hyp}(!,\psi; \lambda_1,\ldots, \lambda_n;\rho_1,\ldots, \rho_m)$
without assuming that $p$ is relatively prime to $2, 3, \ldots,
\mathrm {max}(n,m)$. In fact, I don't know $\mathrm
{Kl}_2(\psi;1,1)|_{\eta_\infty}$ explicitly if $p=2$.

\medskip
The paper is organized as follows. In \S 1, we calculate a direct
factor of the restriction of the global Fourier transformation at
$\eta_{\infty'}$. In \S 2, we use a group theoretical argument to
deduce Theorems 0.1-0.4 from results in \S 1, and we prove
Proposition 0.5 at the end of this section. We prove Proposition 0.7
in \S 3.

\medskip Shortly after the first version of the paper was finished,
Abbes and Saito \cite{AS} were able to calculate the local Fourier
transformations for a class of monomial Galois representations of
Artin--Schreier-Witt type. Their results are more general than ours,
and they used a blowing-up technique and the ramification theory of
Kato. The method used in this paper is more direct and global. In
the first version of this paper, we work over an algebraically
closed field $k$ and prove Theorems 0.2-0.4. Theorems 0.1 is
inspired by \cite{AS} and is suggested by the referee.

\section{Key lemmas}

The main result of this section is Lemma 1.5, which calculates a
direct factor of
$$\Big(R^1\pi_{2!}(\mathscr L_\psi(f(x,t'))\otimes \pi_1^\ast
\mathscr K)\Big)|_{\eta_\infty'},$$ where $\pi_1,\pi_2:\mathbb
G_{m,k}\times_k\mathbb G_{m,k}\to\mathbb G_{m,k}$ are the
projections, $f(x,t')\in k[x,1/x,t',1/t']$ satisfies certain
conditions and $\mathscr K$ is a rank $1$ lisse $\overline {\mathbb
Q}_\ell$-sheaf on $\mathbb G_{m,k}$ such that $\mathscr
K|_{1}=\overline{\mathbb Q}_\ell$. Here we consider the rational
point $1$ of $\mathbb G_{m,k}$ as a closed subscheme of $\mathbb
G_{m,k}$.

\begin{lemma} Let $R$ be a commutative ring and let
$h(x)=\sum\limits_{i=-m}^n c_ix^i$ be a Laurent polynomial with
$c_i\in R$ and $m,n\geq 0$.

(i) Suppose $m,n\geq 1$ and $c_n$ and $c_{-m}$ are units in $R$.
Then the $R$-morphism
$$H:{\mathbb G}_{m,R}\to {\mathbb A}_R^1,\quad
x\mapsto h(x)$$ corresponding to the $R$-algebra homomorphism
$$R[x]\to R[x,1/x],\quad x\mapsto h(x)$$ is a finite morphism, and it
can be extended to an $R$-morphism
$$\hat H:{\mathbb P}_R^1\to {\mathbb P}_R^1$$ which in terms of homogenous
coordinates can be described by
$$[x_0:x_1]\mapsto [x_0^nx_1^m: \sum_{i=-m}^n
c_ix_0^{n-i}x_1^{m+i}].$$ The canonical diagram
$$\begin{array}{rcl}
{\mathbb G}_{m,R}&\hookrightarrow & {\mathbb P}_R^1 \\
{\scriptstyle H}\downarrow&& \downarrow {\scriptstyle \hat H}\\
{\mathbb A}_R^1&\hookrightarrow &{\mathbb P}_R^1
\end{array}$$ is Cartesian.

(ii) Suppose $m=0$, $n\geq 1$ and $c_n$ is a unit in $R$. Then the
$R$-morphism
$$H:{\mathbb A}_R^1\to {\mathbb A}_R^1,\quad x\mapsto h(x)$$ is a finite morphism,
and it can be extended to an $R$-morphism
$$\hat H:{\mathbb P}_R^1\to {\mathbb P}_R^1$$ which in terms of homogenous
coordinates can be described by
$$[x_0:x_1]\mapsto [x_0^n: \sum_{i=0}^n
c_ix_0^{n-i}x_1^i].$$ The canonical diagram
$$\begin{array}{rcl}
{\mathbb A}_R^1&\hookrightarrow & {\mathbb P}_R^1 \\
{\scriptstyle H}\downarrow&& \downarrow {\scriptstyle \hat H}\\
{\mathbb A}_R^1&\hookrightarrow &{\mathbb P}_R^1
\end{array}$$ is Cartesian.
\end{lemma}

\begin{proof} We sketch a proof of (i). Consider the $R$-algebra homomorphism of
graded rings
\begin{eqnarray*}
&& R[x_0,x_1]\to R[x_0,x_1], \\
&& x_0\mapsto x_0^nx_1^m,\quad x_1\mapsto \sum_{i=-m}^n
c_ix_0^{n-i}x_1^{m+i}.
\end{eqnarray*}
This homomorphism does not preserve degrees, but maps a homogenous
polynomial of degree $d$ to a homogenous polynomial of degree
$(m+n)d$. However by \cite{EGA II} 2.4.7 (i) and 2.8.2, it still
defines an $R$-morphism $\hat H:{\mathbb P}_R^1\to {\mathbb P}_R^1$.
Here we need the fact that $c_n$ and $c_{-m}$ are units in order for
the domain of definition of $\hat H$ to be the whole projective line
${\mathbb P}_R^1$. One can verify that $\hat H^{-1}({\mathbb
A}_R^1)={\mathbb G}_{m,R}.$ This implies the diagram in the lemma is
Cartesian. So $H$ is a proper morphism. By definition, $H$ is an
affine morphism. So $H$ is a finite morphism by \cite{EGA III}
4.4.2.
\end{proof}

\begin{lemma} Let $g(x,z')\in k[x,1/x,z']$, and let $a_0(x)\in k[x,1/x]$ be the
constant term of $g(x,z')$ considered as a polynomial of $z'$. Write
$g(x,z')$ as
$$g(x,z')=\sum_{i=-m}^n c_i(z')x^i,$$ where $m, n\geq 0$ and $c_i(z')\in k[z']$.
Let $\mathscr K$ be a lisse $\overline{\mathbb Q}_\ell$-sheaf on
$\mathbb G_{m,k}$, and let $\bar k$ be an algebraic closure of $k$.

(i) Suppose $m, n\geq 1$ and $c_n(0)c_{-m}(0)\not=0$. Denote by $G$
the $k$-morphism $$G:\mathbb G_{m,k}\times_k\mathbb A_k^1\to \mathbb
A_{k}^1\times_k \mathbb A_k^1,\quad (x,z')\to (g(x,z'),z').$$ Let
$$S=\{a\in\bar k\;|\; a=a_0(b) \hbox{ for some } b\in \bar k-\{0\}
\hbox { satisfying }\frac{\mathrm d a_0}{\mathrm dx}(b)=0\}.$$ For
any $x\in \bar k-S$, there exists a neighborhood of $(x,0)$ in
$\mathbb A_k^1\times_k \mathbb A_k^1$ in which
$G_!\pi_1^\ast{\mathscr K}$ is lisse and $R^iG_!\pi_1^\ast{\mathscr
K}$ $(i\geq 1)$ vanish, where $\pi_1:\mathbb G_{m,k}\times \mathbb
A_k^1\to \mathbb G_{m,k}$ is the projection.

(ii) Suppose $m=0$, $n\geq 1$ and $c_n(0)\not=0$. Denote by $G$ the
$k$-morphism
$$G:\mathbb A_{k}^1\times_k\mathbb A_{k}^1\to \mathbb
A_{k}^1\times_k \mathbb \mathbb A_{k}^1,\quad (x,z')\to
(g(x,z'),z').$$ Let $$S=\{a\in\bar k| a=a_0(b) \hbox{ for some }
b\in \bar k \hbox { satisfying }\frac{\mathrm d
a_0}{dx}(b)=0\}\cup\{a_0(0)\}.$$ For any $x\in \bar k-S$, there
exists a neighborhood of $(x,0)$ in $\mathbb A_k^1\times_k \mathbb
A_k^1$ in which $G_!p_1^\ast j_!{\mathscr K}$ is lisse and
$R^iG_!p_1^\ast j_!{\mathscr K}$ $(i\geq 1)$ vanish, where
$j:\mathbb G_{m,k}\hookrightarrow \mathbb A_k^1$ is the open
immersion, and $p_1:\mathbb A_k^1\times_k \mathbb A_k^1\to \mathbb
A_k^1$ is the projection to the first factor.
\end{lemma}

\begin{proof} The proofs of (i) and (ii) are similar. We give a proof (i).
The Jacobian $G$ is $$\mathrm {det}\left(\begin{array}{cc}
\frac{\partial}{\partial x}(g(x,z'))& \frac{\partial}{\partial
z'}(g(x,z')) \\
\frac{\partial}{\partial x}(z')&\frac{\partial}{\partial z'}(z')
\end{array}
\right) \equiv \frac{\mathrm da_0}{\mathrm dx}  \mod z'.$$ For
$z'=0$, the above Jacobian vanishes if and only if $\frac{\mathrm
da_0}{\mathrm dx} =0.$ So $G:{\mathbb G}_{m,k}\times_k {\mathbb
A}_k^1 \to {\mathbb A}_k^1\times_k {\mathbb A}_k^1$ is \'etale at
$(b,0)$ if $\frac{\mathrm da_0}{\mathrm dx}(b)\not=0$. Hence $G$ is
\'etale at every preimage of $(x,0)$ if $x\not\in S$. Next we prove
$G$ is a finite morphism over some neighborhood of $0$ in $\mathbb
A_k^1$. This will imply (i). Set $A=k[z']$. Then $G$ is the morphism
$$G: {\mathbb G}_{m,A}\to {\mathbb A}_A^1$$ corresponding to the
$A$-algebra homomorphism
$$
A[x]\to A[x,1/x], \quad x\mapsto g(x,z').$$ Let $R=A_{c_nc_{-m}}$ be
the localization of $A$ with respect to $c_n(z')c_{-m}(z')\in A$.
Then $c_n(z')$ and $c_{-m}(z')$ are unit in $R$. By Lemma 1.1 (i),
$G:\mathbb G_{m,R}\to \mathbb A_R^1$ is a finite morphism. Since
$c_n(0)c_{-m}(0)\not=0$, $0$ is contained in the image of the open
immersion
$$\mathrm {Spec}\, R=\mathrm {Spec}\, A_{c_nc_{-m}}
\hookrightarrow \mathrm {Spec}\, A={\mathbb A}_k^1.$$ Our assertion
follows.
\end{proof}

From now on we assume the characteristic $p$ of $k$ is distinct from
$2$.

\begin{lemma} Let ${\mathscr O}_{{\mathbb
G}_{m,k}\times_k {\mathbb A}_{k}^1, (1,0)}^{\mathrm h}$ (resp.
${\mathscr O}_{{\mathbb A}^1_k\times_k {\mathbb A}_{k}^1,
(0,0)}^{\mathrm h}$) be the henselization of ${\mathbb
G}_{m,k}\times_k {\mathbb A}_{k}^1$ (resp. ${\mathbb A}^1_k\times_k
{\mathbb A}^1_k$) at the point $(1,0)$ (resp. $(0,0)$), let
$({\mathbb G}_{m,k}\times_k {\mathbb A}^1_k)_{(1,0)}^{\mathrm {h}}$
(resp. $({\mathbb A}^1_k\times_k {\mathbb A}^1_k)^{\mathrm
{h}}_{(0,0)}$ ) be their spectrum, let
\begin{eqnarray*}
G:\mathbb G_{m,k}\times_k{\mathbb A}^1_k&\to&
\mathbb A_k^1\times_k {\mathbb A}^1_k,\\
(x,z')&\mapsto& (g(x, z'),z')
\end{eqnarray*} be the $k$-morphism defined by some
$g(x,z')\in k[x,1/x,z']$ with the property
\begin{eqnarray*}
&&g(1,z')=0,\\
&&\frac{\partial g}{\partial x}(1,z')=0,\\
&&\frac{1}{2}\frac{\partial^2 g}{\partial x^2}(1,0) \hbox { is a
nonzero square in } k,
\end{eqnarray*}
and let
$$\omega: ({\mathbb
A}^1_k\times_k {\mathbb A}^1_k)^{\mathrm {h}}_{(0,0)}\to ({\mathbb
A}^1_k\times_k {\mathbb A}^1_k)^{\mathrm {h}}_{(0,0)}$$ be the
morphism induced on henselizations by the morphism
$${\mathbb A}_k^1\times_k {\mathbb A}^1_k\to
{\mathbb A}_k^1\times_k{\mathbb A}^1_k, \quad (x,z')\mapsto (x^2,
z').$$ Then there exists an isomorphism $$\phi: ({\mathbb
G}_{m,k}\times_k {\mathbb A}^1_k)^{\mathrm
{h}}_{(1,0)}\stackrel\cong\to ({\mathbb A}^1_k\times_k {\mathbb
A}^1_k)^{\mathrm {h}}_{(0,0)}$$ such that the morphism on
henselizations
$$G^{\mathrm
{h}}: ({\mathbb G}_{m,k}\times_k {\mathbb A}^1_k)^{\mathrm
{h}}_{(1,0)}\to ({\mathbb A}^1_k\times_k {\mathbb A}^1_k)^{\mathrm
{h}}_{(0,0)}$$ induced by $G$ coincides with $\omega\circ\phi$.
\end{lemma}

\begin{proof} Since $g(1,z')=\frac{\partial g}{\partial x}(1,z')=0$,
we have $\frac{g(x,z')}{(x-1)^2}\in k[x,1/x,z']$. Hence
$\frac{g(x,z')}{(x-1)^2}$ lies in ${\mathscr O}^{\mathrm
{h}}_{{\mathbb G}_{m,k}\times_k {\mathbb A}^1_k, (1,0)}$ and its
image in the residue field is $\frac{1}{2}\frac{\partial^2
g}{\partial x^2}(1,0)$, which is a nonzero square in $k$. Since
$\mathrm{char}\,k\not=2$, this element in the residue field has two
distinct square roots. By the Hensel lemma \cite{EGA IV} 18.5.13,
$\frac{g(x,z')}{(x-1)^2}$ has two distinct square roots in
${\mathscr O}^{\mathrm {h}}_{{\mathbb G}_{m,k}\times_k {\mathbb
A}^1_k, (1,0)}$. Let $\delta(x-1,z')$ be one of the square roots.
Its image in the residue field is nonzero, and we have
$$g(x,z')=((x-1)\delta(x-1,z'))^2.$$ Consider the morphism
$$({\mathbb G}_{m,k}\times_k {\mathbb A}^1_k)^{\mathrm {h}}_{(1,0)}\to {\mathbb
A}_k^1\times_k{\mathbb A}^1_k$$ corresponding to the $k$-algebra
homomorphism
\begin{eqnarray*}
&&k[x,z']\to {\mathscr O}_{{\mathbb G}_{m,k}\times_k {\mathbb
A}^1_k,
(1,0)}^{\mathrm {h}},\\
&&x\mapsto (x-1)\delta(x-1,z'),\quad z'\mapsto z'.
\end{eqnarray*}
It maps the Zariski closed point of $({\mathbb G}_{m,k}\times_k
{\mathbb A}^1_k)^{\mathrm h}_{(1,0)}$ to $(0,0)$. Denote the induced
morphism on henselizations by
$$\phi: ({\mathbb G}_{m,k}\times_k {\mathbb A}^1_k)^{\mathrm {h}}_{(1,0)}
\to ({\mathbb A}^1_k\times_k {\mathbb A}^1_k)^{\mathrm
{h}}_{(0,0)}.$$ Then $G^{\mathrm h}=\omega\circ\phi$. Let's prove
$\phi$ is an isomorphism. By \cite{EGA IV} 17.6.3, it suffices to
the homomorphism induced by $\phi$ on completions
$$\hat \phi: \widehat {\mathscr O}_{{\mathbb
A}_k^1\times_k {\mathbb A}^1_k, (0,0)}\to \widehat {\mathscr
O}_{{\mathbb G}_{m,k}\times_k {\mathbb A}^1_k, (1,0)}$$ is an
isomorphism. This homomorphism is given by
\begin{eqnarray*}
&& k[[x,z']] \to k[[x-1,z']], \\
&& x\mapsto (x-1)\delta(x-1,z'), \quad z'\mapsto z'.
\end{eqnarray*}
Here we regard $\delta(x-1,z')$ as a square root of
$\frac{g(x,z')}{(x-1)^2}$ in $k[[x-1,z']]$. Using the fact that the
image of $\delta(x-1,z')$ in the residue field, that is, the
constant term of $\delta(x-1,z')$, is nonzero, one can check the
above homomorphism on formal power series rings is indeed an
isomorphism.\end{proof}

\begin{lemma} Let $g(x,z')\in k[x,1/x,z']$ such that
\begin{eqnarray*}
&&g(1,z')=0,\\
&&\frac{\partial g}{\partial x}(1,z')=0,\\
&&\frac{1}{2}\frac{\partial^2 g}{\partial x^2}(1,0) \hbox { is a
nonzero square in } k,
\end{eqnarray*}
Write $g(x,z')$ as
$$g(x,z')=\sum_{i=-m}^n c_i(z')x^i,$$ where $m, n\geq 0$ and $c_i(z')\in k[z']$.
Let $\mathscr K$ be a rank $1$ lisse $\overline{\mathbb
Q}_\ell$-sheaf on $\mathbb G_{m,k}$ such that $\mathscr
K|_{1}\cong\overline{\mathbb Q}_\ell$.

(i) Suppose $m, n\geq 1$ and $c_n(0)c_{-m}(0)\not=0$. Denote by $G$
the $k$-morphism $$G:\mathbb G_{m,k}\times_k\mathbb G_{m,k}\to
\mathbb A_{k}^1\times_k \mathbb G_{m,k},\quad (x,t')\to
(g(x,1/t'),t').$$ For any $\lambda\in k$,  $\bigl({\mathscr
K}_{\chi_2}(\lambda t'^s)\otimes
G(\chi_2,\psi)\bigr)|_{\eta_\infty'}$ is a direct factor of
$$\biggl(R^1\mathrm {pr}_{2!}\bigl( RG_! (\pi_1^\ast {\mathscr
K}_Q)\otimes {\mathscr L}_\psi(\lambda
xt'^s)\bigr)\biggr)|_{\eta_{\infty'}},$$ where
$$\pi_1:\mathbb G_{m,k}\times_k\mathbb G_{m,k}\to\mathbb G_{m,k},\quad
\mathrm{pr_2}:\mathbb A_k^1\times_k\mathbb G_{m,k}\to\mathbb
G_{m,k}$$ are the projections to the first and second factors,
respectively.

(ii) Suppose $m=0$, $n\geq 1$ and $c_n(0)\not=0$. Denote by $G$ the
$k$-morphism
$$G:\mathbb A_{k}^1\times_k\mathbb G_{m,k}\to \mathbb
A_{k}^1\times_k \mathbb \mathbb G_{m,k},\quad (x,t')\to
(g(x,1/t'),t').$$ For any $\lambda\in k$, $\bigl({\mathscr
K}_{\chi_2}(\lambda t'^s)\otimes
G(\chi_2,\psi)\bigr)|_{\eta_{\infty'}}$ is a direct factor of
$$\biggl(R^1\mathrm {pr}_{2!}\bigl( RG_! (\mathrm{pr}_1^\ast j_! {\mathscr
K})\otimes {\mathscr L}_\psi(\lambda
xt'^s)\bigr)\biggr)|_{\eta_{\infty'}},$$ where $j:\mathbb
G_{m,k}\hookrightarrow \mathbb A_k^1$ is the open immersion, and
$$\mathrm{pr}_1:\mathbb A_k^1\times_k\mathbb G_{m,k}\to\mathbb A_k^1,\quad
\mathrm{pr_2}:\mathbb A_k^1\times_k\mathbb G_{m,k}\to\mathbb
G_{m,k}$$ are the projections.
\end{lemma}

\begin{proof} We give a proof of (i). The morphism $G$ can be
extended to a $k$-morphism
$$\bar G:\mathbb G_{m,k}\times_k(\mathbb G_{m,k}\cup\{\infty'\})\to
\mathbb A_k^1\times_k(\mathbb G_{m,k}\cup\{\infty'\})$$ as follows:
We have
$$\mathbb G_{m,k}\cup\{\infty'\}=\mathrm{Spec}\, k[z'].$$ The
morphism $\bar G$ is
$$\bar G:\mathbb G_{m,k}\times_k(\mathbb G_{m,k}\cup\{\infty'\})\to
\mathbb A_k^1\times_k(\mathbb G_{m,k}\cup\{\infty'\}),\quad
(x,z')\mapsto (g(x,z'),z').$$ Fix notations by the following
commutative diagram:
$$\begin{array}{ccccccc}
{\mathbb G}_{m,k}&=&{\mathbb G}_{m,k}&=&{\mathbb G}_{m,k}&& \\
\uparrow\pi_1&&\uparrow \bar \pi_1&& \uparrow \tilde \pi_1&&\\
{\mathbb G}_{m,k}\times_ k {\mathbb G}_{m,k}&\hookrightarrow&
{\mathbb G}_{m,k}\times_k({\mathbb
G}_{m,k}\cup\{\infty'\})&\hookleftarrow& {\mathbb
G}_{m,R}={\mathbb G}_{m,k}\times_k \mathrm {Spec}\, R&\hookrightarrow & {\mathbb P}_R^1\\
\downarrow G&&\downarrow \bar G&& \downarrow \tilde G&& \downarrow
\hat G \\
{\mathbb A}_k^1\times_k {\mathbb G}_{m,k}&\stackrel \iota
\hookrightarrow& {\mathbb A}_k^1\times_k({\mathbb
G}_{m,k}\cup\{\infty'\})&\hookleftarrow& {\mathbb A}^1_R={\mathbb
A}_k^1\times_k \mathrm {Spec}\, R
&\stackrel \kappa\hookrightarrow & {\mathbb P}_R^1={\mathbb P}_k^1\times_k \mathrm {Spec}\, R\\
\mathrm {pr}_2 \downarrow &&\downarrow \bar{\mathrm
{pr}}_2&&\downarrow \tilde {\mathrm {pr}}_2&& \downarrow \hat
{\mathrm {pr}}_2  \\
{\mathbb G}_{m,k}&\stackrel \delta \hookrightarrow& {\mathbb
G}_{m,k}\cup\{\infty'\}=\mathrm {Spec}\, A&\hookleftarrow& \mathrm
{Spec}\, R&=& \mathrm {Spec}\, R,
\end{array}$$
where $R=k[z']_{c_nc_{-m}}$ is the localization of $k[z']$ with
respect to $c_n(z')c_{-m}(z')$, $\pi_1$, $\bar \pi_1$ and
$\tilde\pi_1$ are projections to the first factors, $\mathrm
{pr}_2$, $\bar {\mathrm {pr}}_2$, $\tilde{\mathrm {pr}}_2$ and
$\hat{\mathrm {pr}}_2$ are projections to the second factors,
$\kappa$, $\iota$ and $\delta$ are the canonical open immersions. We
have
\begin{eqnarray*}
\biggl(R\mathrm {pr}_{2!}\bigl( RG_!(\pi_1^\ast{\mathscr K})\otimes
{\mathscr L}_\psi(\lambda xt'^s) \bigr)\biggr)|_{\eta_{\infty'}}
&\cong& \biggl(R\bar {\mathrm {pr}}_{2!}\bigl( R\bar G_!(\bar
\pi_1^\ast{\mathscr K})\otimes \iota_! {\mathscr L}_\psi(\lambda
xt'^s)
\bigr)\biggr)|_{\eta_{\infty'}}\\
&\cong& \biggl(R\tilde {\mathrm {pr}}_{2!}\bigl( R\tilde G_!(\tilde
\pi_1^\ast{\mathscr K})\otimes (\iota_! {\mathscr L}_\psi(\lambda
xt'^s))|_{{\mathbb A}_R^1}
\bigr)\biggr)|_{\eta_{\infty'}}\\
&\cong& \biggl(R\hat {\mathrm {pr}}_{2\ast}\bigl(\kappa_! R\tilde
G_!(\tilde \pi_1^\ast{\mathscr K})\otimes \kappa_! ((\iota_!
{\mathscr L}_\psi(\lambda xt'^s))|_{{\mathbb A}_R^1})
\bigr)\biggr)|_{\eta_{\infty'}}.
\end{eqnarray*}
Let $a_0(x)\in k[x,1/x]$ be the constant term of $g(x,z')\in
k[x,1/x,z']$ consider as a polynomial of $z'$. Since
$\frac{\partial^2 g}{\partial x^2}(1,0)$ is nonzero, $\frac{\mathrm
da_0}{\mathrm dx}$ is nonzero. Let $S$ be the set of those $a\in
\bar k$ such that $a=a_0(b)$ for some $b\in\bar k$ satisfying
$\frac{\mathrm da_0}{\mathrm dx}(b)=0$. Then $S$ is a finite set.
Since $g(1,z')=\frac{\partial g}{\partial x}(1,z')=0,$ we have $0\in
S$. By Lemma 1.2, $R\tilde G_!(\tilde \pi_1^\ast{\mathscr K})$ is
lisse in a neighborhood of $(x,\infty')$ if $x\not\in S$. By
\cite{L} 1.3.1.2, $\iota_! {\mathscr L}_\psi(xt')$ is universally
strongly locally acyclic relative to $\bar{\mathrm {pr}}_2$. Since
$\iota_! {\mathscr L}_\psi(\lambda xt'^s)$ is obtained from $\iota_!
{\mathscr L}_\psi(xt')$ by base change, $\iota_! {\mathscr
L}_\psi(\lambda xt'^s)$ is also universally strongly locally acyclic
relative to $\bar{\mathrm {pr}}_2$. It follows that the vanishing
cycle
$$R\Phi_{\bar\eta_{\infty'}}\bigl(\kappa_!R\tilde G_!(\tilde
\pi_1^\ast{\mathscr K})\otimes \kappa_! ((\iota_! {\mathscr
L}_\psi(\lambda xt'^s))|_{{\mathbb A}_R^1})\bigr)$$ is supported at
those points $(x,\infty')$ in ${\mathbb P}_k^1\times \infty'$ with
$x=\infty$ or $x\in S$. Moreover, as $\kappa_! R\tilde G_!(\tilde
\pi_1^\ast{\mathscr K})\otimes \kappa_! ((\iota_! {\mathscr
L}_\psi(\lambda xt'^s))|_{{\mathbb A}_R^1})$ vanishes on ${\mathbb
P}_k^1\times \infty'$, its vanishing cycle and nearby cycle
coincide, that is,
$$R\Phi_{\bar\eta_{\infty'}}\bigl(\kappa_! R\tilde G_!(\tilde
\pi_1^\ast{\mathscr K})\otimes \kappa_! ((\iota_! {\mathscr
L}_\psi(\lambda xt'^s))|_{{\mathbb A}_R^1})\bigr) \cong
R\Psi_{\bar\eta_{\infty'}}\bigl(\kappa_! R\tilde G_!(\tilde
\pi_1^\ast{\mathscr K})\otimes \kappa_! ((\iota_! {\mathscr
L}_\psi(\lambda xt'^s))|_{{\mathbb A}_R^1})\bigr).$$ So we have
\begin{eqnarray*}
&&\biggl(R\hat {\mathrm {pr}}_{2\ast}\bigl(\kappa_! R\tilde
G_!(\tilde \pi_1^\ast{\mathscr K})\otimes \kappa_! ((\iota_!
{\mathscr L}_\psi(\lambda xt'^s))|_{{\mathbb
A}_R^1})\bigr)\biggr)_{\bar\eta_{\infty'}}\\
&\cong& R\Gamma\biggl({\mathbb P}_{\bar k}^1\times \infty',
R\Psi_{\bar\eta_{\infty'}}\bigl(\kappa_! R\tilde G_!(\tilde
\pi_1^\ast{\mathscr K})\otimes \kappa_! ((\iota_! {\mathscr
L}_\psi(\lambda xt'^s))|_{{\mathbb A}_R^1})\bigr)\biggr)\\
&\cong & R\Gamma\biggl({\mathbb P}_{\bar k}^1\times \infty',
R\Phi_{\bar\eta_{\infty'}}\bigl(\kappa_!R\tilde G_!(\tilde
\pi_1^\ast{\mathscr K})\otimes \kappa_! ((\iota_! {\mathscr
L}_\psi(\lambda xt'^s))|_{{\mathbb A}_R^1})\bigr)\biggr) \\
&\cong &\bigoplus_{x\in S}
\biggl(R\Phi_{\bar\eta_{\infty'}}\bigl(\kappa_!  R\tilde G_!(\tilde
\pi_1^\ast{\mathscr K})\otimes \kappa_! ((\iota_! {\mathscr
L}_\psi(\lambda xt'^s))|_{{\mathbb
A}_R^1})\bigr)\biggr)_{(x,\infty')}
\\
&& \quad \bigoplus
\biggl(R\Phi_{\bar\eta_{\infty'}}\bigl(\kappa_!R\tilde G_!(\tilde
\pi_1^\ast{\mathscr K})\otimes \kappa_! ((\iota_! {\mathscr
L}_\psi(\lambda xt'^s))|_{{\mathbb
A}_R^1})\bigr)\biggr)_{(\infty,\infty')}.
\end{eqnarray*}
In particular, taking $x=0\in S$, we see
$$\biggl(R\Phi_{\bar\eta_{\infty'}}\bigl(\kappa_!  R\tilde
G_!(\tilde \pi_1^\ast{\mathscr K})\otimes \kappa_! ((\iota_!
{\mathscr L}_\psi(\lambda xt'^s))|_{{\mathbb
A}_R^1})\bigr)\biggr)_{(0,\infty')}$$ is a direct factor of
$$\biggl(R\hat {\mathrm {pr}}_{2\ast}\bigl(\kappa_!  R\tilde G_!(\tilde
\pi_1^\ast{\mathscr K})\otimes \kappa_! ((\iota_! {\mathscr
L}_\psi(\lambda xt'^s))|_{{\mathbb
A}_R^1})\bigr)\biggr)_{\bar\eta_{\infty'}}.$$ We have
$$(\tilde\pi_1^\ast{\mathscr K})|_{({\mathbb G}_{m,k}\times_k ({\mathbb
G}_{m,k}\cup\{\infty'\}))^{\mathrm {h}}_{(1,\infty')}}\cong
\overline {\mathbb Q}_l$$ since ${\mathscr K}$ is lisse on $\mathbb
G_{m,k}$ and $\mathscr K|_{1}\cong \overline {\mathbb Q}_l$, where
$({\mathbb G}_{m,k}\times_k ({\mathbb
G}_{m,k}\cup\{\infty'\}))_{(1,\infty')}^{\mathrm {h}}$ is the
spectrum of the henselization of ${\mathbb G}_{m,k}\times_k
({\mathbb G}_{m,k}\cup\{\infty'\})$ at $(1,\infty')$. Let $\omega$
and $\phi$ be the morphisms defined in Lemma 1.3. Then $(R\tilde
G_!(\tilde \pi_1^\ast{\mathscr K}))|_{({\mathbb A}^1_k\times_k
({\mathbb G}_{m,k}\cup\{\infty'\}))^{\mathrm {h}}_{(0,\infty')}}$
has a direct factor $(\omega\circ\phi)_\ast(\overline{\mathbb
Q}_l)$. Since $\phi$ is an isomorphism, we have
$(\omega\circ\phi)_\ast(\overline{\mathbb Q}_l)\cong
\omega_\ast(\overline{\mathbb Q}_l).$ But $\omega_\ast(\overline
{\mathbb Q}_l)$ has a direct factor $(\tilde {\mathrm {pr}}_1^\ast
j_!{\mathscr K}_{\chi_2})|_{({\mathbb A}^1_k\times_k ({\mathbb
G}_{m,k}\cup\{\infty'\}))^{\mathrm {h}}_{(0,\infty')}},$ where
$\tilde {\mathrm {pr}}_1:\mathbb A_k^1\times_k\mathrm{Spec}\,R\to
\mathbb A_k^1$ is the projection. So
$$\biggl(R^1\Phi_{\bar\eta_{\infty'}}\bigl(\kappa_!
R\tilde G_!(\tilde \pi_1^\ast{\mathscr K})\otimes \kappa_! ((\iota_!
{\mathscr L}_\psi(\lambda xt'^s))|_{{\mathbb
A}_R^1})\bigr)\biggr)_{(0,\infty')}$$ has a direct factor
$$\biggl(R^1\Phi_{\bar\eta_{\infty'}}\bigl(\kappa_! \tilde {\mathrm {pr}}_1^\ast
j_!{\mathscr K}_{\chi_2}\otimes \kappa_! ((\iota_! {\mathscr
L}_\psi(\lambda xt'^s))|_{{\mathbb
A}_R^1})\bigr)\biggr)_{(0,\infty')}.$$ Let
$f:\eta_{\infty'}\to\eta_{\eta'}$ be the $k$-morphism defined by the
$k$-algebra homomorphism $$k((1/t'))\to ((1/t')),\quad t'\mapsto
\lambda t'^s.$$ Then we have
\begin{eqnarray*}
\biggl(R^1\Phi_{\bar\eta_{\infty'}}\bigl(\kappa_! \tilde {\mathrm
{pr}}_1^\ast j_! {\mathscr K}_{\chi_2}\otimes \kappa_! ((\iota_!
{\mathscr L}_\psi(\lambda xt'^s))|_{{\mathbb
A}_R^1})\bigr)\biggr)_{(0,\infty')}&\cong&
f^\ast {\mathfrak F}^{(0,\infty')}({\mathscr K}_{\chi_2}|_{\eta_0})\\
&\cong &f^\ast\bigl({\mathscr K}_{\chi_2}\otimes
G(\chi_2,\psi)\bigr)|_{\eta_{\infty'}}\\
&\cong& \bigl({\mathscr K}_{\chi_2}(\lambda t'^s)\otimes
G(\chi_2,\psi)\bigr)|_{\eta_{\infty'}},
\end{eqnarray*}
where for the first isomorphism, we use \cite{SGA 4 1/2} Th.
finitude 3.7 about the base change of vanishing cycles, and for the
second isomorphism, we use \cite{L} 2.5.3.1 (ii). This proves Lemma
1.4.
\end{proof}

In the rest of this section, we denote by $\pi_1,\pi_2:\mathbb
G_{m,k}\times_k\mathbb G_{m,k} \to \mathbb G_{m,k}$ the two
projections.

\begin{lemma} Let $\mathscr K$ be a rank
$1$ lisse $\overline{\mathbb Q}_\ell$-sheaf on $\mathbb G_{m,k}$
such that $\mathscr K|_1\cong \overline{\mathbb Q}_\ell$, let
$f(x,t')\in k[x,1/x, t',1/t']$, and let $a_s(x)t'^s$ be the highest
degree term of $f(x,t')$ considered as a Laurent polynomial of $t'$,
where $a_s(x)\in k[x,1/x]$. Suppose
\begin{eqnarray*}
&& \frac{\partial f}{\partial x}(1,t')=0,\\
&&\lambda=\frac{1}{2}\frac{\partial^2 a_s}{\partial x^2}\bigg|_{x=1}
\hbox { is nonzero}. \end{eqnarray*} Write $f(x,t')$ as
$$f(x,t')=\sum_{i=-m}^n c_i(t')x^i,$$
where $m,n\geq 0$ and $c_i(t')\in k[t',1/t']$. Suppose one of the
following conditions hold:

(a) $m, n\geq 1$ and $\left(\frac{c_n(t')}{t'^s}\cdot
\frac{c_{-m}(t')}{t'^s}\right)|_{t'=\infty}\not=0$.

(b) $m=0$, $n\geq 1$ and $\frac{c_n(t')}{t'^s}|_{t'=\infty}\not=0$.

(c) $m\geq 1$, $n=0$ and
$\frac{c_{-m}(t')}{t'^s}|_{t'=\infty}\not=0$.

\noindent Then $$\bigl({\mathscr L}_\psi(f(1,t'))\otimes {\mathscr
K}_{\chi_2}(\lambda t'^s)\otimes
G(\chi_2,\psi)\bigr)|_{\eta_{\infty'}}$$ is a direct factor of
$$\biggl(R^1\pi_{2!}\bigl({\mathscr L}_\psi(f(x,t'))\otimes
{\pi_1}^\ast{\mathscr K}\bigr)\biggr)|_{\eta_{\infty'}}.$$
\end{lemma}

\begin{proof} Let
$$g(x,1/t')=\frac{f(x,t')-f(1,t')}{\lambda t'^s}$$ and set $z'=1/t'$.
In the case (a), $g(x,z')$ satisfies the conditions of Lemma 1.4
(i). Let $G$ be the $k$-morphism
$$
G:{\mathbb G}_{m,k}\times_k{\mathbb G}_{m,k}\to {\mathbb
A}_k^1\times_k{\mathbb G}_{m,k},\quad  G(x,t')=(g(x,1/t'),t')$$ and
let $\mathrm {pr}_2: {\mathbb A}_k^1\times_k {\mathbb G}_{m,k}\to
{\mathbb G}_{m,k}$ be the projection. We have $\mathrm {pr}_2\circ
G=\pi_2.$ Combined with the projection formula (\cite{SGA 4} XVII
5.2.9), we get
\begin{eqnarray*}
&&R\pi_{2!}\bigl({\mathscr L}_\psi(f(x,t'))\otimes{\pi_1}^\ast
{\mathscr K}\bigr)\\&\cong& R\pi_{2!}\bigl({\mathscr
L}_\psi(f(x,t')-f(1,t'))\otimes \pi_2^\ast {\mathscr
L}_\psi(f(1,t'))\otimes \pi_1^\ast \mathscr K\bigr)\\
&\cong& {\mathscr L}_\psi(f(1,t'))\otimes R\pi_{2!}\bigl({\mathscr
L}_\psi(f(x,t')-f(1,t'))\otimes
\pi_1^\ast \mathscr K\bigr)\\
&\cong& {\mathscr L}_\psi(f(1,t'))\otimes R(\mathrm {pr}_2\circ G)_!
\bigl(G^\ast {\mathscr L}_\psi(\lambda xt'^s)\otimes \pi_1^\ast
{\mathscr
K}\bigr)\\
&\cong & {\mathscr L}_\psi(f(1,t'))\otimes R\mathrm {pr}_{2!}
\bigl({\mathscr L}_\psi(\lambda xt'^s)\otimes RG_! \pi_1^\ast
{\mathscr K}\bigr).
\end{eqnarray*} Our assertion then follows from Lemma 1.4 (i).

In the case (b), $f(x,t')$ is a polynomial of $x$, and $g(x,z')$
satisfies the conditions of Lemma 1.4 (ii). Let $G$ be the
$k$-morphism
$$
G:{\mathbb A}_k^1\times_k{\mathbb G}_{m,k}\to {\mathbb
A}_k^1\times_k{\mathbb G}_{m,k},\quad  G(x,t')=(g(x,1/t'),t'),$$ let
\begin{eqnarray*}
\mathrm{pr}_1: {\mathbb A}_k^1\times_k {\mathbb G}_{m,k}\to {\mathbb
A}_k^1,&&\quad \mathrm{pr}_2: {\mathbb A}_k^1\times_k {\mathbb
G}_{m,k}\to {\mathbb G}_{m,k}
\end{eqnarray*}
be the projections, and let $j:{\mathbb G}_{m,k}\to{\mathbb A}_k^1$
be the open immersion. The same argument as above shows that
$$R\mathrm{pr}_{2!}\bigl({\mathscr L}_\psi(f(x,t'))\otimes{\mathrm{pr}_1}^\ast
j_!{\mathscr K}\bigr)\cong {\mathscr L}_\psi(f(1,t'))\otimes
R\mathrm {pr}_{2!} \bigl({\mathscr L}_\psi(\lambda xt'^s)\otimes
RG_! \mathrm{pr}_1^\ast j_! {\mathscr K}\bigr).$$ Moreover, we have
\begin{eqnarray*}
&&R\mathrm{pr}_{2!}\bigl({\mathscr
L}_\psi(f(x,t'))\otimes{\mathrm{pr}_1}^\ast j_!{\mathscr
K}\bigr)\\
&\cong& R\mathrm{pr}_{2!}\bigl({\mathscr L}_\psi(f(x,t'))\otimes
(j\times\mathrm {id}_{\mathbb G_{m,k}})_!\pi_1^\ast {\mathscr
K}\bigr)\\
&\cong& R\mathrm{pr}_{2!} (j\times\mathrm {id}_{\mathbb G_{m,k}})_!
\bigl((j\times\mathrm {id}_{\mathbb G_{m,k}})^\ast {\mathscr
L}_\psi(f(x,t'))\otimes
\pi_1^\ast {\mathscr K}\bigr)\\
&\cong& R\pi_{2!}\bigl({\mathscr L}_\psi(f(x,t'))\otimes{\pi_1}^\ast
{\mathscr K}\bigr).
\end{eqnarray*}
So we have
$$R\pi_{2!}\bigl({\mathscr L}_\psi(f(x,t'))\otimes{\pi_1}^\ast
{\mathscr K}\bigr)\cong {\mathscr L}_\psi(f(1,t'))\otimes R\mathrm
{pr}_{2!} \bigl({\mathscr L}_\psi(\lambda xt'^s)\otimes RG_!
\mathrm{pr}_1^\ast j_! {\mathscr K}\bigr).$$  Our assertion then
follows from Lemma 1.4 (ii).

In the case (c), $f(1/x,t')$ satisfies the condition (b). Our
assertion follows from the case (b) and the fact that
$$R\pi_{2!}\bigl({\mathscr L}_\psi(f(x,t'))\otimes{\pi_1}^\ast
{\mathscr K}\bigr)\cong R\pi_{2!}\bigl({\mathscr
L}_\psi(f(1/x,t'))\otimes{\pi_1}^\ast (\mathrm{inv}^\ast {\mathscr
K})\bigr),$$ where $\mathrm {inv}:\mathbb G_{m,k}\to\mathbb G_{m,k}$
is the morphism $x\mapsto 1/x$.
\end{proof}

\begin{lemma} Let $\mathscr K$ be a rank $1$ lisse $\overline {\mathbb
Q}_\ell$-sheaf on $\mathbb G_{m,k}$ tamely ramified at $0$ and
$\infty$, and let $P$ be the morphism
$$
P:\mathbb G_{m,k}\times_k \mathbb G_{m,k}\to\mathbb G_{m,k},\quad
(x,t')\mapsto xt'.$$ Then there exists a lisse $\overline {\mathbb
Q}_\ell$-sheaf $\mathscr K_P$ on $\mathbb G_{m,k}$ tamely ramified
at $0$ and $\infty$ such that $$P^\ast \mathscr K\cong \pi_1^\ast
\mathscr K_P\otimes \pi_2^\ast \mathscr K,\quad \mathscr
K_P|_1\cong\overline{\mathbb Q}_\ell.$$
\end{lemma}

\begin{proof} Let $\mathscr
K^{-1}=\mathscr Hom(\mathscr K,\overline{\mathbb Q}_\ell)$. We claim
that the canonical morphism
$$\pi_1^\ast \pi_{1\ast}(P^\ast\mathscr K\otimes \pi_2^\ast
\mathscr K^{-1})\to P^\ast\mathscr K\otimes \pi_2^\ast
 \mathscr K^{-1}$$ is an isomorphism. To prove
this, we may make base extension from $k$ to its algebraic closure
$\bar k$. There exists a character
$\chi:\varprojlim_{(N,p)=1}\bbmu_N(\bar k)\to\overline{\mathbb
Q}_\ell^\ast$ such that $\mathscr K|_{\mathbb G_{m,\bar k}}\cong
\mathscr K_{\chi}$. By \cite{SGA 4 1/2} Sommes trig. (1.3.1), we
have $P^\ast \mathscr K_\chi\otimes \pi_2^\ast\mathscr
K_\chi^{-1}\cong \pi_1^\ast \mathscr K_\chi.$ We have
\begin{eqnarray*}
\pi_1^\ast\pi_{1\ast} (\pi_1^\ast \mathscr K_\chi)\cong
\pi_1^\ast(\pi_{1\ast}\pi_1^\ast \mathscr K_\chi)\cong
\pi_1^\ast(\mathscr K_\chi\otimes\pi_{1\ast}\pi_1^\ast\overline
{\mathbb Q}_\ell) \cong \pi_1^\ast(\mathscr K_\chi\otimes \overline
{\mathbb Q}_\ell) \cong \pi_1^\ast\mathscr K_\chi.
\end{eqnarray*}
Hence the canonical morphism $\pi_1^\ast
\pi_{1\ast}(\pi_1^\ast\mathscr K_\chi)\to \pi_1^\ast\mathscr K_\chi$
is an isomorphism. It follows that the canonical morphism
$\pi_1^\ast \pi_{1\ast}(P^\ast\mathscr K\otimes \pi_2^\ast \mathscr
K^{-1})\to P^\ast\mathscr K\otimes \pi_2^\ast \mathscr K^{-1}$ is an
isomorphism over $\bar k$, and hence over $k$. Let $\mathscr
K_P=\pi_{1\ast}(P^\ast\mathscr K\otimes \pi_2^\ast \mathscr
K^{-1}).$ Then $P^\ast \mathscr K\cong \pi_1^\ast \mathscr
K_P\otimes \pi_2^\ast\mathscr K.$ Taking the restriction at $(1,1)$
on both sides of this isomorphism, we get $\mathscr K_P|_1\cong
\overline{\mathbb Q}_\ell$.
\end{proof}

\begin{lemma} Let $\alpha(t),\gamma(t),\delta(t)\in k[1,1/t]$,
let $\gamma,\delta:\mathbb {\mathbb G}_{m,k}\to {\mathbb A}_k^1$ be
the $k$-morphisms defined by $\gamma(t),\delta(t)$, respectively,
and let $\mathscr K$ be a $\overline{\mathbb Q}_\ell$-sheaf on
$\mathbb G_{m,k}$. Denote by ${\mathfrak F}:D_c^b(\mathbb
A_k^1,\overline{\mathbb Q}_\ell)\to D_c^b(\mathbb
A_k^1,\overline{\mathbb Q}_\ell)$ the global Fourier transformation.
We have
$$\delta^\ast {\mathfrak F}\Big(\gamma_!\bigl({\mathscr
L}_\psi(\alpha(t))\otimes {\mathscr K}\bigr)\Big)\cong R\pi_{2!}
\bigl({\mathscr L}_\psi(\alpha(t)+\gamma(t)\delta(t'))\otimes
\pi_1^\ast {\mathscr K}\bigr)[1].$$
\end{lemma}

\begin{proof} Fix notations by the following commutative diagram:
$$\begin{array}{crccccclc}
&&&& {\mathbb G}_{m,k}\times_k {\mathbb G}_{m,k} &&&& \\
&&&\swarrow &&\searrow  &&& \\
&&{\mathbb G}_{m,k}\times_k {\mathbb A}_k^1&&&& {\mathbb A}_k^1\times_k {\mathbb G}_{m,k}&& \\
&\swarrow &&\searrow&&\swarrow &&\searrow \scriptstyle \mathrm {pr}_2& \\
{\mathbb G}_{m,k}&&&& {\mathbb A}_k^1\times_k {\mathbb A}_k^1&&&& {\mathbb G}_{m,k}\\
&\scriptstyle \gamma\searrow &&\scriptstyle
p_1\swarrow && \searrow \scriptstyle p_2&& \swarrow \scriptstyle \delta &\\
&&{\mathbb A}_k^1&&&&{\mathbb A}_k^1&&\\
&&&\searrow&&\swarrow&&& \\
&&&&\mathrm {Spec}\, k&&&&
\end{array}$$
Applying the proper base change theorem and the projection formula,
we get
\begin{eqnarray*}
&&\delta^\ast Rp_{2!}\Big(p_1^\ast \gamma_!\bigl({\mathscr
L}_\psi(\alpha(t))\otimes {\mathscr K}\bigr)\otimes {\mathscr
L}_\psi(tt')\Big)\\
&\cong& R\mathrm {pr}_{2!} (\mathrm {id}_{{\mathbb A}_k^1}\times
\delta)^\ast \Big(p_1^\ast \gamma_! \bigl({\mathscr
L}_\psi(\alpha(t))\otimes {\mathscr K}\bigr)\otimes {\mathscr
L}_\psi(tt')\Big) \\
&\cong& R\mathrm {pr}_{2!}\Big(\bigl(p_1\circ (\mathrm
{id}_{{\mathbb A}_k^1}\times \delta)\bigr)^\ast \gamma_!
\bigl({\mathscr L}_\psi(\alpha(t))\otimes {\mathscr K}\bigr)\otimes
(\mathrm {id}_{{\mathbb A}_k^1}\times \delta)^\ast {\mathscr
L}_\psi(tt')\Big) \\
&\cong& R\mathrm {pr}_{2!}\Big((\gamma \times {\mathrm
{id}}_{{\mathbb G}_{m,k}})_!\pi_1^\ast \bigl({\mathscr
L}_\psi(\alpha(t))\otimes {\mathscr K}\bigr) \otimes ({\mathrm
{id}}_{{\mathbb A}_k^1}\times \delta)^\ast {\mathscr
L}_\psi(tt')\Big)
\\
& \cong& R\mathrm {pr}_{2!} (\gamma\times {\mathrm {id}}_{{\mathbb
G}_{m,k}})_!\Big(\pi_1^\ast ({\mathscr L}_\psi(\alpha(t))\otimes
{\mathscr K})\otimes (\gamma\times \mathrm {id}_{{\mathbb
G}_{m,k}})^\ast({\mathrm {id}}_{{\mathbb A}_k^1}\times \delta)^\ast
{\mathscr
L}_\psi(tt') \Big)\\
&\cong&R\pi_{2!} \bigl({\mathscr
L}_\psi(\alpha(t)+\gamma(t)\delta(t'))\otimes \pi_1^\ast {\mathscr
K}\bigr).
\end{eqnarray*} Our assertion follows.
\end{proof}

\begin{proposition} Let $\mathscr K$ be a rank $1$ lisse $\overline{\mathbb
Q}_\ell$-sheaf on ${\mathbb G}_{m,k}$ tamely ramified at $0$ and at
$\infty$, and let $\alpha(t),\gamma(t)\in k[t,1/t]$. Suppose
$\frac{\mathrm d\gamma}{\mathrm dt}$ divides $\frac{\mathrm
d\alpha}{\mathrm dt}$ in $k[t,1/t]$, and let
\begin{eqnarray*}
&&\delta(t)=-\frac{\mathrm d\alpha}{\mathrm d\gamma}
=-\frac{\frac{\mathrm d\alpha}{dt}}{\frac{\mathrm d\gamma}{\mathrm
dt}},\\
&&\beta(t')=\alpha(t')+\gamma(t')\delta(t').
\end{eqnarray*}
Let $a_st^s$ be the highest degree term of $\alpha(t)$, and let $r$
be the highest degree of terms in $\gamma(t)$. Suppose
$p=\mathrm{char}\,k$ is relatively prime to $2, r,s$ and $s-r$.
Write
$$\alpha(xt')+\gamma(xt')\delta(t')=\sum_{i=-m}^nc_i(t')x^i$$ with
$m,n\geq 0$ and $c_i(t')\in k[t',1/t']$. Suppose furthermore that
one of the following conditions holds:

(a) $m, n\geq 1$ and $\left(\frac{c_n(t')}{t'^s}\cdot
\frac{c_{-m}(t')}{t'^s}\right)|_{t'=\infty}\not=0$.

(b) $m=0$, $n\geq 1$ and $\frac{c_n(t')}{t'^s}|_{t'=\infty}\not=0$.

(c) $m\geq 1$, $n=0$ and
$\frac{c_{-m}(t')}{t'^s}|_{t'=\infty}\not=0$.

\noindent Then $$\biggl({\mathscr L}_\psi(\beta(t'))\otimes
{\mathscr K}\otimes {\mathscr
K}_{\chi_2}\left(\frac{1}{2}s(s-r)a_st'^s\right)\otimes
G(\chi_2,\psi)\biggr)|_{\eta_{\infty'}}$$ is a direct factor of
$$\biggl(\delta^\ast \mathscr H^0{\mathfrak F}\Big(\gamma_!\bigl({\mathscr
L}_\psi(\alpha(t))\otimes {\mathscr
K}\bigr)\Big)\biggr)|_{\eta_{\infty'}}.$$
\end{proposition}

\begin{proof} Let
$\lambda=\frac{1}{2}s(s-r)a_s$.  By Lemma 1.7, it suffices to show
$$\biggl({\mathscr L}_\psi(\beta(t'))\otimes
{\mathscr K}\otimes {\mathscr K}_{\chi_2}(\lambda t'^s)\otimes
G(\chi_2,\psi)\biggr)|_{\eta_{\infty'}}$$ is a direct factor of
$$R^1\pi_{2!} \bigl({\mathscr
L}_\psi(\alpha(t)+\gamma(t)\delta(t'))\otimes \pi_1^\ast {\mathscr
K}\bigr)|_{\eta_{\infty'}}.$$ Introduce a change of variables
$$x=\frac{t}{t'}$$ and let
$$f(x,t')=\alpha(t)+\gamma(t)\delta(t')=\alpha(xt')+\gamma(xt')\delta(t').$$
Then we have
$$\beta(t')=f(1,t').$$
Consider the isomorphism $$\tau: {\mathbb G}_{m,k}\times_k {\mathbb
G}_{m,k}\to {\mathbb G}_{m,k}\times_k {\mathbb G}_{m,k},\quad
(x,t')\mapsto\left(xt',t'\right).$$ We have $\pi_2\tau=\pi_2.$ So we
have
\begin{eqnarray*}
R\pi_{2!} \big({\mathscr
L}_\psi(\alpha(t)+\gamma(t)\delta(t'))\otimes {\pi_1}^\ast{\mathscr
K}\big) &\cong& R\pi_{2!}\tau^\ast \big({\mathscr
L}_\psi(\alpha(t)+\gamma(t)\delta(t'))\otimes
{\pi_1}^\ast{\mathscr K}\big)\\
&\cong& R\pi_{2!}\big(\tau^\ast {\mathscr
L}_\psi(\alpha(t)+\gamma(t)\delta(t'))\otimes
\tau^\ast {\pi_1}^\ast{\mathscr K}\big)\\
&\cong& R\pi_{2!}\bigl({\mathscr L}_\psi(f(x,t'))\otimes
({\pi_1}^\ast{\mathscr K}_P\otimes {\pi_2}^\ast {\mathscr
K})\bigr)\\
&\cong& R\pi_{2!}\bigl({\mathscr L}_\psi(f(x,t'))\otimes
{\pi_1}^\ast{\mathscr K}_P\bigr)\otimes {\mathscr K}.
\end{eqnarray*}
Here for third isomorphism, we use Lemma 1.6, and for the fourth
isomorphism, we use the projection formula. We will verify shortly
that the conditions of Lemma 1.5 hold for $f(x,t')$. We then
conclude from lemma 1.5 that
$$\bigl({\mathscr
L}_\psi(f(1,t'))\otimes {\mathscr K}_{\chi_2}(\lambda t'^s)\otimes
G(\chi_2,\psi)\bigr)|_{\eta_{\infty'}}$$ is a direct factor of
$$\biggl(R^1\pi_{2!}\bigl({\mathscr L}_\psi(f(x,t'))\otimes
{\pi_1}^\ast{\mathscr K}_P\bigr)\biggr)|_{\eta_{\infty'}},$$ and
hence
$$\bigl({\mathscr L}_\psi(\beta(t'))\otimes
{\mathscr K}\otimes {\mathscr K}_{\chi_2}(\lambda t'^s)\otimes
G(\chi_2,\psi)\bigr)|_{\eta_{\infty'}}$$ is a direct factor of
$$\biggl(R^1\pi_{2!} \bigl({\mathscr
L}_\psi(\alpha(t)+\gamma(t)\delta(t'))\otimes {\pi_1}^\ast{\mathscr
K}\bigr)\biggr)|_{\eta_{\infty'}}.$$ This proves our assertion.

Let's verify $f(x,t')$ satisfies the conditions of Lemma 1.5. We
have
\begin{eqnarray*}
\frac{\partial f}{\partial x}&=&\frac{\partial}{\partial
x}\Big(\alpha(xt')+\gamma(xt')\delta(t')\Big)\\
&=& t'\frac{\mathrm d\alpha}{\mathrm dt}(xt')+t'\frac{\mathrm
d\gamma}{\mathrm dt}(xt')\delta(t').
\end{eqnarray*}
Since $\delta(t)=-\frac{\frac{\mathrm d\alpha}{dt}}{\frac{\mathrm
d\gamma}{dt}}$, we have
$$\frac{\partial f}{\partial x}(1,t')=t'\frac{\mathrm
d}{\mathrm dt'}(\alpha(t'))+t'\frac{\mathrm d}{\mathrm
dt'}(\gamma(t'))\delta(t')=0.$$ Write
\begin{eqnarray*}
\alpha(t)&=&a_st^s+a_{s-1}t^{s-1}+\cdots,\\
\gamma(t)&=&b_r t^r+b_{r-1}t^{r-1}+\cdots.
\end{eqnarray*}
We have
\begin{eqnarray*}
f(x,t')&=&\alpha(xt')+\gamma(xt')\delta(t')\\
&=&\alpha(xt')-\gamma(xt')\frac{\frac{\mathrm d}{\mathrm
dt'}(\alpha(t'))}
{\frac{\mathrm d}{\mathrm dt'}(\gamma(t'))}\\
&=& a_s(xt')^s +a_{s-1}(xt')^{s-1}+\cdots \\
&&\quad -\Big(b_r(xt')^r
+b_{r-1}(xt')^{r-1}+\cdots\Big)\left(\frac{sa_s
t'^{s-1}+(s-1)a_{s-1}t'^{s-2}+\cdots}{rb_r
t'^{r-1}+(r-1)b_{r-1}t'^{r-2}+\cdots}\right).
\end{eqnarray*}
From this expression, we see that the highest degree term of
$f(x,t')$ considered as a Laurent polynomial of $t'$ is $a_s(x)t'^s$
with
$$a_s(x)=a_sx^s-\frac{s}{r}a_s x^r.$$
Note that we have
$$\frac{1}{2}\frac{\partial^2}{\partial
x^2}(a_s(x))|_{x=1}=\frac{1}{2}s(s-r)a_s=\lambda.$$ Since $p$ is
relatively prime to $2,r,s$ and $s-r$, we have
$\frac{1}{2}\frac{\partial^2}{\partial x^2}(a_s(x))|_{x=1}\not=0$.
\end{proof}

\begin{remark} This remark gives a philosophical reason why we
introduce the change of variable $x=\frac{t}{t'}$. By the facts on
classical analysis recalled in the introduction, we expect a
dominant part of
$$\biggl(R^1\pi_{2!} \bigl({\mathscr
L}_\psi(\alpha(t)+\gamma(t)\delta(t'))\otimes {\pi_1}^\ast{\mathscr
K}\bigr)\biggr)|_{\eta_{\infty'}}$$ comes from the contribution on
the curve
$$\frac{\partial}{\partial t}\Big(\alpha(t)+\gamma(t)\delta(t')\Big)=0.$$
By assumption, we have $\frac{\mathrm d\alpha}{\mathrm
dt}+\frac{\mathrm d\gamma}{\mathrm dt}\delta(t)=0.$ So $t=t'$ is a
branch of the above curve. By the change of variable
$x=\frac{t}{t'}$, we see $x=1$ is a branch of the curve. So a
dominant part of the above cohomology should come from the
contribution at $x=1$. This is indeed the case as we have seen.
\end{remark}

In the next section, we will deduce Theorems 0.1-0.4 from the
following lemma by a group theoretical argument:

\begin{lemma} Let $\gamma(t)=t^r$, let
$\alpha(t)\in k[t,1/t]$, let
\begin{eqnarray*}
&&\delta(t)=-\frac{\mathrm d\alpha}{\mathrm d\gamma}
=-\frac{1}{rt^{r-1}}\frac{\mathrm d\alpha}{dt},\\
&&\beta(t')=\alpha(t')+\gamma(t')\delta(t'),
\end{eqnarray*}
and let $\mathscr K$ be a rank $1$ lisse $\overline{\mathbb
Q}_\ell$-sheaf on ${\mathbb G}_{m,k}$ tamely ramified at $0$ and at
$\infty$. We have $\gamma(0)=0$ and $\gamma(\infty)=\infty$.

(i) Suppose $\alpha(t)$ is of the form
$$\alpha(t)=\frac{a_{-s}}{t^s}+\frac{a_{-(s-1)}}{t^{s-1}}+\cdots+\frac{a_{-1}}{t}$$
with $a_{-s}\not=0$, and suppose $r,s\geq 1$ and
$p=\mathrm{char}\,k$ is relatively prime to $2$, $r$, $s$ and $r+s$.
We have $\delta(0')=\infty'$ and
$$\biggl({\mathscr L}_\psi(\beta(t'))\otimes {\mathscr K}\otimes
{\mathscr K}_{\chi_2}\left(\frac{1}{2}s(r+s)a_{-s}t'^s\right)\otimes
G(\chi_2,\psi)\biggr)|_{\eta_{0'}}$$ is a direct factor of
$\delta^\ast {\mathfrak
F}^{(0,\infty')}\biggl(\Big(\gamma_\ast\bigl({\mathscr L}_\psi
(\alpha(t))\otimes {\mathscr K}\bigr)\Big)|_{\eta_0}\bigg)$.

(ii) Suppose $\alpha(t)$ is of the form
$$\alpha(t)=a_st^s+a_{s-1}t^{s-1}+\cdots+a_1t$$ with $a_s\not=0$, and
suppose $s>r\geq 1$ and $p$ is relatively prime $2$, $r$, $s$ and
$s-r$. We have $\delta(\infty')=\infty'$ and
$$\biggl({\mathscr L}_\psi(\beta(t'))\otimes {\mathscr
K}\otimes {\mathscr
K}_{\chi_2}\left(\frac{1}{2}s(s-r)a_st'^s\right)\otimes
G(\chi_2,\psi)\bigg)|_{\eta_{\infty'}}$$ is a direct factor of
$\delta^\ast {\mathfrak
F}^{(\infty,\infty')}\bigg(\Big(\gamma_\ast\bigl({\mathscr
L}_\psi(\alpha(t))\otimes {\mathscr
K}\bigr)\Big)|_{\eta_\infty}\bigg).$

(iii) Suppose $\alpha(t)$ is of the form
$$\alpha(t)=a_st^s+a_{s-1}t^{s-1}+\cdots+a_1t$$ with $a_s\not=0$, and
suppose $r>s\geq 1$ and $p$ is relatively prime to $2$, $r$, $s$ and
$s-r$. We have $\delta(\infty')=0'$ and
$$\biggl({\mathscr L}_\psi(\beta(t'))\otimes {\mathscr K}\otimes
{\mathscr K}_{\chi_2}\left(\frac{1}{2}s(s-r)a_st'^s\right)\otimes
G(\chi_2,\psi)\bigg)|_{\eta_{\infty'}}$$ is a direct factor of
$\delta^\ast {\mathfrak
F}^{(\infty,0')}\bigg(\Big(\gamma_\ast\bigl({\mathscr
L}_\psi(\alpha(t))\otimes {\mathscr
K}\bigr)\Big)|_{\eta_\infty}\bigg).$
\end{lemma}

\begin{proof} ${}$

(i) Let
\begin{eqnarray*}
&&\gamma_1(t)=\gamma(1/t)=t^{-r}, \\
&&\alpha_1(t)=\alpha(1/t)=a_{-s}t^s+a_{-(s-1)}t^{s-1}+\cdots +a_{-1}t,\\
&&\delta_1(t)=-\frac{\frac{\mathrm d}{\mathrm
dt}(\alpha_1(t))}{\frac{\mathrm d}{\mathrm
dt}(\gamma_1(t))}=\delta(1/t)=\frac{s}{r}a_{-s}t^{r+s}+\frac{s-1}{r}a_{-(s-1)}t^{r+s-1}+\cdots
+\frac{1}{r}a_{-1}t^{r+1},\\
&&\beta_1(t')=\alpha_1(t')+\gamma_1(t')\delta_1(t')=\beta(1/t'),
\end{eqnarray*} and let $\lambda=\frac{1}{2}s(r+s)a_{-s}$.
From the above expression for $\delta_1(t')=\delta(1/t')$, we get
$\delta(0')=\infty'$. Note that as a Laurent polynomial of $x$, the
highest and the lowest degree terms of
$\alpha_1(xt')+\gamma_1(xt')\delta_1(t')$ are $(a_{-s}t'^s)x^s$ and
$\Big(\frac{s}{r}a_{-s}t'^s+\frac{s-1}{r}a_{-(s-1)}t'^{s-1}+\cdots
+\frac{1}{r}a_{-1}t'\Big)x^{-r}$, respectively, and we have
$$\frac{a_{-s}t'^s}{t'^s}|_{t'=\infty'}\not=0, \quad
\frac{\frac{s}{r}a_{-s}t'^s+\frac{s-1}{r}a_{-(s-1)}t'^{s-1}+\cdots
+\frac{1}{r}a_{-1}t'}{t'^s}|_{t'=\infty}\not=0.$$ So we can apply
Proposition 1.8 (a) to $\alpha_1(t)$, $\gamma_1(t)$ and
$\mathrm{inv}^\ast \mathscr K$, and we see
$$\Big({\mathscr
L}_\psi(\beta_1(t'))\otimes \mathrm{inv}^\ast {\mathscr K}\otimes
{\mathscr K}_{\chi_2}(\lambda t'^s)\otimes
G(\chi_2,\psi)\Big)|_{\eta_{\infty'}}$$ is a direct factor of
$$\biggl(\delta_1^\ast\mathscr H^0 {\mathfrak
F}\Big(\gamma_{1!}\bigl({\mathscr L}_\psi (\alpha_1(t))\otimes
\mathrm{inv}^\ast {\mathscr
K}\bigr)\Big)\biggr)|_{\eta_{\infty'}}.$$ We have
\begin{eqnarray*}
{\mathscr L}_\psi(\beta_1(t'))\otimes \mathrm{inv}^\ast {\mathscr
K}\otimes {\mathscr K}_{\chi_2}(\lambda t'^s)\otimes
G(\chi_2,\psi)&\cong&\mathrm{inv}^\ast \Big({\mathscr
L}_\psi(\beta(t'))\otimes {\mathscr K}\otimes {\mathscr
K}_{\chi_2}(\lambda t'^s)\otimes
G(\chi_2,\psi)\Big)\\
\gamma_{1!}\bigl({\mathscr L}_\psi (\alpha_1(t))\otimes
\mathrm{inv}^\ast {\mathscr K}\bigr)&\cong& \gamma_!\bigl({\mathscr
L}_\psi (\alpha(t))\otimes {\mathscr K}\bigr).
\end{eqnarray*} So
$$\Big({\mathscr
L}_\psi(\beta(t'))\otimes  {\mathscr K}\otimes {\mathscr
K}_{\chi_2}(\lambda t'^s)\otimes G(\chi_2,\psi)\Big)|_{\eta_{0'}}$$
is a direct factor of
$$\biggl(\delta^\ast \mathscr H^0{\mathfrak
F}\Big(\gamma_!\bigl({\mathscr L}_\psi (\alpha(t))\otimes {\mathscr
K}\bigr)\Big)\biggr)|_{\eta_{0'}}.$$ On the other hand,
$\gamma_!\bigl({\mathscr L}_\psi(\alpha(t))\otimes {\mathscr
K}\bigr)$ is tamely ramified at $\infty$ since $\alpha(t)$ is a
polynomial of $1/t$ and $(r,p)=1$. By \cite{L} 2.4.3 (iii) b) and
the stationary phase principle \cite{L} 2.3.3.1 (iii), as sheaves on
$\eta_{\infty'}$, we have
\begin{eqnarray*}
\biggl({\mathscr H}^0{\mathfrak F}\Big(\gamma_!\bigl({\mathscr
L}_\psi( \alpha(t))\otimes {\mathscr K} \bigr)\Big)\biggr)|_{
\eta_{\infty'}}\cong {\mathfrak
F}^{(0,\infty')}\bigg(\Big(\gamma_!\bigl({\mathscr
L}_\psi(\alpha(t))\otimes {\mathscr K}\bigr)\Big)|_{\eta_0}\bigg).
\end{eqnarray*}
So
$$\Big({\mathscr
L}_\psi(\beta(t'))\otimes  {\mathscr K}\otimes {\mathscr
K}_{\chi_2}(\lambda t'^s)\otimes G(\chi_2,\psi)\Big)|_{\eta_{0'}}$$
is a direct factor of $\delta^\ast {\mathfrak
F}^{(0,\infty')}\bigg(\Big(\gamma_!\bigl({\mathscr L}_\psi
(\alpha(t))\otimes {\mathscr K}\bigr)\Big)|_{\eta_0}\bigg)$. Finally
note that $\gamma$ induces a finite morphism $\gamma:\mathbb
G_{m,k}\to\mathbb G_{m,k}$. So we can replace $\gamma_!$ by
$\gamma_\ast$ in the above expression. Our assertion follows.

(ii) We have
$$\delta(t)=-\frac{1}{rt^{r-1}}\frac{\mathrm d\alpha}{\mathrm dt}
=-\frac{s}{r}a_st^{s-r}-\frac{s-1}{r}a_{s-1}t^{s-r-1}-\cdots-
\frac{1}{r}a_1t^{-r+1}.$$ Since $s>r$, we have
$\delta(\infty')=\infty'$. Note that
$\alpha(xt')+\gamma(xt')\delta(t')$ is a polynomial of $x$, and the
highest degree term as a polynomial of $x$ is $(a_st'^s)x^s$, and we
have $\frac{a_st'^s}{t'^s}|_{t'=\infty'}\not=0$. So we can apply
Proposition 1.8 (b). Let $\lambda=\frac{1}{2}s(s-r)a_s$. Then
$$\Big({\mathscr L}_\psi(\beta(t'))\otimes {\mathscr K}\otimes
{\mathscr K}_{\chi_2}(\lambda t'^s)\otimes
G(\chi_2,\psi)\Big)|_{\eta_{\infty'}}$$ is a direct factor of
$\bigg(\delta^\ast {\mathscr H}^0{\mathfrak
F}\Big(\gamma_!\bigl({\mathscr L}_\psi( \alpha(t))\otimes {\mathscr
K} \bigr)\Big)\bigg)|_{ \eta_{\infty'}}$. By the stationary phase
principle \cite{L} 2.3.3.1 (iii), as sheaves on $\eta_{\infty'}$, we
have
\begin{eqnarray*}
&&\biggl({\mathscr H}^0{\mathfrak F}\Big(\gamma_!\bigl({\mathscr
L}_\psi( \alpha(t))\otimes {\mathscr K} \bigr)\Big)\biggr)|_{
\eta_{\infty'}}\\
&\cong & {\mathscr
F}^{(0,\infty')}\bigg(\Big(\gamma_!\bigl({\mathscr
L}_\psi(\alpha(t))\otimes{\mathscr K}\bigr)\Big)|_{\eta_0}\bigg)
\bigoplus {\mathscr
F}^{(\infty,\infty')}\bigg(\Big(\gamma_!\bigl({\mathscr
L}_\psi(\alpha(t))\otimes  {\mathscr
K}\bigr)\Big)|_{\eta_\infty}\bigg).
\end{eqnarray*}
Since $\alpha(t)$ is a polynomial of $t$,
$\Big(\gamma_!\bigl({\mathscr L}_\psi(\alpha(t))\otimes {\mathscr
K}\bigr)\Big)|_{\eta_0}$ is tame and hence ${\mathscr
F}^{(0,\infty')}\bigg(\Big(\gamma_!\bigl({\mathscr
L}_\psi(\alpha(t))\otimes {\mathscr K}\bigr)\Big)|_{\eta_0}\bigg)$
is also tame by \cite{L} 2.5.3.1. But the slope of $\bigl({\mathscr
L}_\psi(\beta(t'))\otimes {\mathscr K}\otimes {\mathscr
K}_{\chi_2}(\lambda t'^s)\otimes
G(\chi_2,\psi)\bigr)|_{\eta_{\infty'}}$ is nonzero. So
$\bigl({\mathscr L}_\psi(\beta(t'))\otimes {\mathscr K}\otimes
{\mathscr K}_{\chi_2}(\lambda t'^s)\otimes
G(\chi_2,\psi)\bigr)|_{\eta_{\infty'}}$ is a direct factor of
$\delta^\ast{\mathfrak
F}^{(\infty,\infty')}\bigg(\Big(\gamma_!\bigl({\mathscr
L}_\psi(\alpha(t))\otimes {\mathscr
K}\bigr)\Big)|_{\eta_\infty}\bigg).$

(iii) Again we have
$$\delta(t)=-\frac{1}{rt^{r-1}}\frac{\mathrm d}{\mathrm dt}(\alpha(t))
=-\frac{s}{r}a_st^{s-r}-\frac{s-1}{r}a_{s-1}t^{s-r-1}-\cdots-
\frac{1}{r}a_1t^{-r+1}.$$ Since $s<r$, we have $\delta(\infty')=0'$.
Note that $\alpha(xt')+\gamma(xt')\delta_1(t')$ is a polynomial of
$x$, and the highest degree term as a polynomial of $x$ is
$$\Big(-\frac{s}{r}a_st'^s-\frac{s-1}{r}a_{s-1}t'^{s-1}-\cdots-
\frac{1}{r}a_1t'\Big)x^r,$$ and we have
$$
\frac{-\frac{s}{r}a_st'^s-\frac{s-1}{r}a_{s-1}t'^{s-1}-\cdots-\frac{1}{r}a_1t'}
{t'^s}|_{t'=\infty}\not=0.$$ So we can apply Proposition 1.8 (b).
Let $\lambda=\frac{1}{2}s(s-r)a_s$. Then $$\bigl({\mathscr
L}_\psi(\beta(t'))\otimes {\mathscr K}\otimes {\mathscr
K}_{\chi_2}(\lambda t'^s)\otimes
G(\chi_2,\psi)\bigr)|_{\eta_{\infty'}}$$ is a direct factor of
$\biggl(\delta^\ast \mathscr H^0{\mathfrak
F}\Big(\gamma_!\bigl({\mathscr L}_\psi( \alpha(t))\otimes {\mathscr
K} \bigr)\Big)\biggr)|_{\eta_{\infty'}}.$ We have
\begin{eqnarray*}
{\mathfrak F}\Big(\gamma_!\bigl({\mathscr L}_\psi( \alpha(t))\otimes
{\mathscr K}\bigr)[1]\Big)_{\bar 0} &\cong& R\Gamma_c\bigl(\mathbb
A^1_{\bar k},\gamma_!\bigl({\mathscr L}_\psi(
\alpha(t))\otimes {\mathscr K}\bigr)\bigr)[2] \\
&\cong& R\Gamma_c\bigl(\mathbb G_{m,\bar k},{\mathscr L}_\psi(
\alpha(t))\otimes {\mathscr K} \bigr)[2].
\end{eqnarray*}
Since the sheaf ${\mathscr L}_\psi(\alpha(t))\otimes {\mathscr K}$
on $\mathbb G_{m,k}$ is nonconstant and lisse of rank $1$, and thus
irreducible, we have
$$H^2_c\bigl(\mathbb G_{m,\bar k},{\mathscr L}_\psi(
\alpha(t))\otimes {\mathscr K} \bigr)=0$$ by \cite{D2} 1.4.1.3. So
by \cite{L} 2.3.2, we have a short exact sequence
\begin{eqnarray*}
0&\to& H^1_c\bigl(\mathbb G_{m,\bar k},{\mathscr L}_\psi(
\alpha(t))\otimes {\mathscr K} \bigr)\\
&\to&\biggl( \mathscr H^0{\mathfrak F}\Big(\gamma_!\bigl({\mathscr
L}_\psi( \alpha(t))\otimes
{\mathscr K} \bigr)\Big)\biggr)|_{\bar \eta_{0'}}\\
&\to& {\mathscr F}^{(\infty,0')}\Bigg(\Big(\gamma_!\bigl({\mathscr
L}_\psi(\alpha(t'))\otimes {\mathscr
K}\bigr)\Big)|_{\eta_{\infty}}\bigg)\to 0.
\end{eqnarray*}
Note that $H^1_c\bigl(\mathbb G_{m,\bar k},{\mathscr L}_\psi(
\alpha(t))\otimes {\mathscr K} \bigr)$ defines an unramified
representation of $\mathrm{Gal}(\bar\eta_{0'}/\eta_{0'})$, whereas
the representation defined by $\bigl({\mathscr
L}_\psi(\beta(t'))\otimes {\mathscr K} \otimes {\mathscr
K}_{\chi_2}(\lambda t'^s)\otimes
G(\chi_2,\psi)\bigr)|_{\eta_{\infty'}}$ has no unramified
subquotient. So $\bigl({\mathscr L}_\psi(\beta(t'))\otimes {\mathscr
K}\otimes {\mathscr K}_{\chi_2}(\lambda t'^s)\otimes
G(\chi_2,\psi)\bigr)|_{\eta_{\infty'}}$ is a direct factor of
$\delta^\ast {\mathfrak
F}^{(\infty,0')}\bigg(\Big(\gamma_!\bigl({\mathscr
L}_\psi(\alpha(t))\otimes {\mathscr
K}\bigr)\Big)|_{\eta_\infty}\bigg).$
\end{proof}

\section{Proof of Theorems}

In this section, we use a group theoretical argument to deduce
Theorems 0.1-0.4 from Lemma 1.10. We assume $k$ is algebraically
closed throughout this section except when we prove Theorem 0.1 in
the end.

\begin{lemma} Let $G$ be a group, let $H$ be a subgroup of $G$, and let
$V$ be a finite dimensional representation of $G$ over a field of
characteristic $0$.

(i) If the restriction of this representation to $H$ is semisimple
and if $H$ has finite index in $G$, then $V$ is a semisimple
representation of $G$.

(ii) If $V$ is a semisimple representation of $G$ and if $H$ is a
normal subgroup of $G$, then the restriction of this representation
to $H$ is semisimple.
\end{lemma}

\begin{proof}

(i) Let $U$ be a subspace of $V$ invariant under the action of $G$.
We need to  show that there exists a homomorphism $\pi:V\to U$ such
that $\pi|_U$ is the identity and $g^{-1}\pi g=\pi$ for all $g\in
G$. Since $V$ is semisimple as a representation of $H$, we can find
a homomorphism $\pi_0:V\to U$ such that $\pi_0|_U$ is the identity
and $h^{-1}\pi_0h=\pi_0$ for all $h\in H$. Let $Hg_1,\ldots, Hg_n$
be a family of representatives of right cosets of $H$, where
$n=[G:H]$. Define $\pi=\frac{1}{n}\sum_{i=1}^n g_i^{-1}\pi_0 g_i.$
One can check $\pi$ has the required property.

(ii) We may assume $V$ is an irreducible representation of $G$. If
$U$ is an irreducible sub-representation of $H$ contained in $V$,
then we have $V=\sum_{g\in G} gU$ since the righthand side is
invariant under $G$. Since $H$ is normal in $G$, each $gU$ is
invariant under the action of $H$. So $gU$ is a representation of
$H$, and it is irreducible. As a sum of irreducible representation,
$V$ is semisimple as a representation of $H$.
\end{proof}

Fix a separable closure $\overline {k((t))}$ of the field $k((t))$
of Laurent series. Let $p=\mathrm{char}\,k$, let $I=\mathrm
{Gal}\biggl(\overline {k((t))}/k((t))\biggr)$ and let $P$ be the
wild inertia subgroup. Then $P$ is a profinite $p$-group, and we
have a canonical isomorphism
$$I/P\stackrel\cong \to \varprojlim_{(N,p)=1}\bbmu_N(k),\quad \sigma\mapsto
\left(\frac{\sigma(\sqrt[N]{t})}{\sqrt[N]{t}}\right).$$ A $\overline
{\mathbb Q}_\ell$-representation of a topological group is always
assumed to be finite dimensional and continuous. Recall that a
representation is called indecomposable if it is not a direct sum of
two proper sub-representations. Any finite dimensional
representation is a direct sum of indecomposable representations.
Any $\overline {\mathbb Q}_\ell$-representation of $I$ defines a
$\overline {\mathbb Q}_\ell$-sheaf on $\eta_0={\mathrm {Spec}}\,
k((t))$. Fix an element $$\zeta=(\zeta_N)\in
\varprojlim_{(N,p)=1}\bbmu_N(k)$$ so that $\zeta_N$ is a primitive
$N$-th root of unity for each $N$ prime to $p$. It is a topological
generator of $\varprojlim\limits_{(N,p)=1}\bbmu_N(k)$. Let $U(n)$ be
the $n$-dimensional representation
$\varprojlim\limits_{(N,p)=1}\bbmu_N(k)\to \mathrm{GL}(n,\overline
{\mathbb Q}_\ell)$ so that $\zeta$ acts through the unipotent
$(n\times n)$-matrix with a single Jordan block
$$\left(\begin{array}{cccc}
1&1&&\\
&1&\ddots& \\
&&\ddots&1\\
&&&1
\end{array}
\right).$$ Composed with the canonical projection
$$I\to I/P\cong \varprojlim_{(N,p)=1}\bbmu_N(k),$$ we get a
$\overline {\mathbb Q}_\ell$-representation of $I$, which we still
denote by $U(n)$. It is indecomposable. We denote the $\overline
{\mathbb Q}_\ell$-sheaf on $\eta_0=\mathrm {Spec}\, k((t))$
corresponding to this representation also by $U(n)$.

Fix an element $g\in I$ so that its image in $I/P\cong
\varprojlim\limits_{(N,p)=1}\bbmu_N(k)$ is $\zeta$. Let $J$ the
inverse image of $\zeta^{\mathbb Z}$ in $I$. Then $J$ is dense in
$I$, and it is the subgroup of $I$ generated by $P$ and $g$. Define
a topology on $J$ so that open subgroups of $P$ form a base of
neighborhoods of the identity element in $J$. Then $P$ is an open
subgroup of $J$, and the restriction to $J$ of the homomorphism $I
\to \varprojlim\limits_{(N,p)=1}\bbmu_N(k)$ induces an isomorphism
of discrete groups
$$J/P \stackrel {\cong}\to \{\zeta^m |m\in \mathbb Z\}\cong \mathbb Z.$$
Let $\chi:\mathbb Z\to \overline{\mathbb Q}_\ell^\ast$ be a
character. Denote by $\mathscr K_\chi$ the representation
$$J\to \overline{\mathbb Q}_\ell^\ast,\quad \sigma g^m\mapsto \chi(m)
\hbox { for any } \sigma\in P \hbox { and } m\in \mathbb Z.$$ For
any integer $n$, the restriction to $J$ of the representation $U(n)$
of $I$ is also denoted by $U(n)$.

\begin{lemma}

(i) Let $V$ be an indecomposable $\overline {\mathbb
Q}_\ell$-representation of $J$. Then there exists an irreducible
representation $W$ of $J$ and a positive integer $n$ such that
$$V\cong W \otimes U(n).$$ This factorization is unique up to
isomorphism, and $W$ is isomorphic to any irreducible
sub-representation of $V$.

(ii) Let $V$ be an irreducible $\overline{\mathbb
Q}_\ell$-representation of $J$. Then there exists a character
$\chi:{\mathbb Z}\to \overline{\mathbb Q}_\ell^\ast$ such that the
representation $V\otimes {\mathscr K}_{\chi^{-1}}$ factors through a
finite discrete quotient group of $J$.

(iii) If $W$ is an irreducible $\overline {\mathbb
Q}_\ell$-representation of $J$, then for any positive integer $n$,
$W\otimes U(n)$ is indecomposable.
\end{lemma}

\begin{proof}

(i) Since $V$ is a $\overline {\mathbb Q}_\ell$-representation and
$P$ is a profinite $p$-group, the restriction of this representation
to $P$ factors through a finite quotient of $P$. In particular, as a
representation of $P$, $V$ is semisimple. Let ${\mathscr S}$ be the
set of isomorphic classes of irreducible representations of $P$
which are sub-representations of $V$. Consider the isotypical
decomposition
$$V=\bigoplus_{W\in {\mathscr S}} V_W$$
of $V$ as a representation of $P$, where for each $W\in \mathscr S$,
$V_W$ is a direct sum of several copies of $W$. Since $P$ is normal
in $J$, for any $h\in J$, $hV_W$ is invariant under the action of
$P$. So as a representation of $P$, we also have
$$V=\bigoplus_{W\in {\mathscr S}}hV_W,$$
and each $hV_W$ is isotypical (i.e. a direct sum of several copies
of an irreducible representation). By the uniqueness of the
isotypical decomposition, $h$ induces a permutation of the set $\{
V_W|W\in{\mathscr S}\}$. So $J$ acts on this set. Since $V$ is
indecomposable as a representation of $J$, $J$ acts transitively.
Fix a representation $W\in \mathscr S$, and let $$J'=\{h\in J |
hV_W=V_W\}$$ be the stabilizer of $V_W$. Then $J'$ has finite index
in $J$ (since the set $\{V_W|W\in \mathscr S\}$ is finite) and
$P\subset J'$. Let $N=[J:J']$. Then $J'$ is generated by $P$ and
$g^N$. The space $V_W$ is invariant under the action of $J'$ and
hence defines a representation of $J'$. We have
$$V\cong \mathrm {Ind}_{J'}^J (V_W).$$
Let $\rho: J'\to \mathrm {GL}(V_W)$ be the homomorphism defined by
the representation $V_W$ of $J'$, and let $\rho^{g^N}:J'\to \mathrm
{GL}(V_W)$ be the homomorphism
$$\rho^{g^N}(\sigma)=\rho(g^{-N}\sigma g^N)$$ for any $\sigma\in J'$.
Then $g^N:V_W\to V_W$ defines an isomorphism from the representation
$\rho^{g^N}$ to the representation $\rho$. But as a representation
of $P$, $V_W$ is a direct sum of several copies of the irreducible
representation $W$. By the uniqueness (up to isomorphism) of such
decomposition, if $\rho_W: P\to \mathrm {GL}(W)$ is the homomorphism
defined by the representation $W$ of $P$, there exists an
isomorphism $F:W\to W$ from the representation $\rho_W^{g^N}$ to the
representation $\rho_W$, where $\rho_W^{g^N}$ is defined by
$$\rho_W^{g^N}(\sigma)=\rho_W(g^{-N}\sigma g^N)$$ for any $\sigma\in
P$. Consider the map
$$\rho':J'\to \mathrm {GL}(W),\quad \sigma g^{mN}\mapsto \rho_W(\sigma)F^m
\hbox{ for any } \sigma\in P \hbox { and any integer } m.$$ One can
verify $\rho'$ defines a representation of $J'$. Since $W$ is an
irreducible representation of $P$, $\rho'$ defines an irreducible
representation of $J'$. We have an isomorphism
$$W\otimes \mathrm {Hom}_{P}(W, V_W)\stackrel\cong \to V_W, \quad w\otimes f\mapsto
f(w)$$ of representations of $P$, where $P$ acts trivially on
$\mathrm {Hom}_P(W,V_W)$. Define the action of $h\in J'$ on $f\in
\mathrm {Hom}_P(W, V_W)$ to be $\rho(h)\circ f\circ \rho'(h^{-1})$.
(Note that $\rho(h)\circ f\circ \rho'(h^{-1})$ is $P$-equivariant
since $P$ is normal in $J'$.) Then the above isomorphism is an
isomorphism of representations of $J'$. Since $P$ acts trivially on
$\mathrm {Hom}_P(W,V_W)$, the action of $J'$ on $\mathrm
{Hom}_P(W,V_W)$ factors through $J'/P.$ Since $V=\mathrm
{Ind}_{J'}^J(V_W)$ is indecomposable, $V_W$ is also indecomposable
as a representation of $J'$. Hence $\mathrm {Hom}_P(W,V_W)$ is
indecomposable as a representation of $J'/P$. But $J'/P$ is a cyclic
group of infinite order generated by the image of $g^N$. Therefore
the Jordan form of $g^N$ on $\mathrm {Hom}_P(W,V_W)$ has a single
Jordan block. Let $\lambda$ be the eigenvalue of $g^N$ on $\mathrm
{Hom}_P(W,V_W)$, and let $K$ be the one dimensional representation
of $J'/P$ so that the image of $g^N$ acts by the scalar
multiplication by $\lambda$. Then as representations of $J'$, we
have
$$W\otimes \mathrm {Hom}_P(W,V_W)\cong (W\otimes K)\otimes
(K^{-1}\otimes \mathrm {Hom}_P(W,V_W)).$$ Now $W\otimes K$ is still
an irreducible representation of $J'$, and $K^{-1}\otimes \mathrm
{Hom}_P(W,V_W)$ is a representation of $J'/P$, and the Jordan form
of $g^N$ on $K^{-1}\otimes \mathrm {Hom}_P(W,V_W)$ has a single
Jordan block with eigenvalue $1$.  Let $n$ be the dimension of
$\mathrm {Hom}_P(W,V_W)$ and let $U$ be the $n\times n$ matrix
$$U=\left(\begin{array}{cccc}
1&1&&\\
&1&\ddots&\\
&&\ddots&1\\
&&&1\end{array} \right).$$ One can verify
$$U^N=\left(\begin{array}{cccc}
1&N&\ast&\ast\\
&1&\ddots&\ast\\
&&\ddots&N\\
&&&1\end{array} \right).$$ In particular, all eigenvalues of $U^N$
are $1$ and $U^N-I$ has rank $n-1$. So $U^N$ has a single Jordan
block. Hence as representations of $J'$, we have
$$K^{-1}\otimes \mathrm {Hom}_P(W,V_W)\cong \mathrm {Res}_{J'}^J(U(n)).$$
We thus have
\begin{eqnarray*}
V&\cong& \mathrm {Ind}_{J'}^J(V_W)\\
&\cong & \mathrm {Ind}_{J'}^J\biggl((W\otimes K)\otimes
\mathrm {Res}_{J'}^J(U(n))\biggr)\\
&\cong& \mathrm {Ind}_{J'}^J(W\otimes K)\otimes U(n).
\end{eqnarray*}
We have seen that $W$ is an irreducible representation of $J'$. Note
that $J'$ is a normal subgroup of $J$ since $P\subset J'\subset J$
and $J/P$ is abelian. Using this fact, the fact that $W\otimes K$ is
an irreducible representation of $J'$, and Proposition 22 in
\cite{Se} \S 7.3, one can check $\mathrm {Res}_{J'}^J \mathrm
{Ind}_{J'}^J(W\otimes K)$ is semisimple. By Lemma 2.1 (i), $\mathrm
{Ind}_{J'}^J(W\otimes K)$ is a semisimple representation of $J$. As
$V\cong \mathrm {Ind}_{J'}^J(W\otimes K)\otimes U(n)$ and $V$ is
indecomposable, $\mathrm {Ind}_{J'}^J(W\otimes K)$ is an irreducible
representation of $J$. This finishes the proof of (i).

(ii) Keep the notations in the proof of (i). We have shown $V\cong
\mathrm {Ind}_{J'}^J(V_W)$. Suppose $V$ is irreducible. Then $V_W$
is irreducible. We have shown $V_W\cong W\otimes \mathrm
{Hom}_P(W,V_W)$. As $V_W$ is irreducible and the action of $g^N$ on
$\mathrm {Hom}_P(W,V_W)$ has a single Jordan block, $\mathrm
{Hom}_P(W,V_W)$ must have dimension 1. So after twisting $W$ by a
one dimensional representation, we may assume $V\cong \mathrm
{Ind}_{J'}^J(W)$, where $W$ is a representation of $J'$ which is
irreducible as a representation of $P$. Let $\rho:J'\to \mathrm
{GL}(W)$ be the homomorphism defined by the representation $W$.
$\rho(P)$ is a finite group. For any $h\in J'$, let $\rho(h)$ act on
this finite group by conjugation. There exists a positive integer
$r$ such that the conjugation by $\rho(h^r)$ induces the identity
map on $\rho(P)$. For example, we can take $r$ to be the order of
the automorphism group of $\rho(P)$. Thus we have
$$\rho(g^{rN})^{-1}\rho(\sigma)\rho(g^{rN})=\rho(\sigma)$$ for any
$\sigma\in P$, that is, $\rho(g^{rN})$ is a $P$-equivariant
isomorphism of $W$. But $W$ is an irreducible representation of $P$.
So $\rho(g^{rN})$ is a scalar multiplication by Schur's Lemma.
Choose $c\in \overline {\mathbb Q}_\ell$ such that $\rho(g^{rN})$ is
the scalar multiplication by $c^{rN}$. Define a character
$\chi:\mathbb Z\to \overline {\mathbb Q}_\ell^\ast$ by $\chi(1)=c.$
We claim that the representation $V\otimes \mathscr K_{\chi^{-1}}$
factors through a finite discrete quotient group of $J$. We have
$$
V\otimes \mathscr K_{\chi^{-1}}\cong \mathrm {Ind}_{J'}^J(W)\otimes
\mathscr K_{\chi^{-1}}\cong\mathrm {Ind}_{J'}^J\bigl(W\otimes
\mathrm {Res}_{J'}^J (\mathscr K_{\chi^{-1}})\bigr).
$$
It suffices to show $W\otimes \mathrm {Res}_{J'}^J \mathscr
K_{\chi^{-1}}$ factors through a finite discrete quotient group of
$J'$. By our construction, $g^{rN}$ acts trivially on $W\otimes
\mathrm {Res}_{J'}^J \mathscr K_{\chi^{-1}}$. Moreover, there exists
an open subgroup $P'$ of $P$ such that $P'$ acts trivially on
$W\otimes \mathrm {Res}_{J'}^J \mathscr K_{\chi^{-1}}$. Since the
subgroup of $J'$ generated by $P'$ and $g^{rN}$ has finite index,
our assertion follows.

(iii) Let $V=W\otimes U(n)$. We can find a filtration
$$V=F_0\supset F_1\supset\cdots\supset F_n=0$$
of $V$ by $J$-invariant subspaces such that $F_i/F_{i+1}\cong W$ for
any $0\leq i\leq n-1$. By the Jordan-H\"{o}lder theorem, any
irreducible subquotient of $V$ is isomorphic to $W$. Write
$$V=V_1 \oplus \cdots \oplus V_m$$ so that each $V_i$ is
indecomposable. By (i), we must have
$$V_i\cong W\otimes U(n_i)$$ for some positive integer $n_i$.
We thus have
$$V=W\otimes U(n)\cong W\otimes (U(n_1)\oplus\cdots \oplus
U(n_m)).$$ Since $U(n)$ is trivial as a representation of $P$, we
have a canonical isomorphism
\begin{eqnarray*}
&&\mathrm {Hom}_P(W,W)\otimes U(n)\to \mathrm {Hom}_P(W,W\otimes
U(n)),\\
&& f\otimes u\mapsto (w\mapsto f(w)\otimes u) \hbox { for any }
f\in\mathrm {Hom}_P(W,W),\; u\in U(n),\; w\in W.
\end{eqnarray*}
This isomorphism is $J$-equivariant. Through this isomorphism, the
action of $g$ on $\mathrm {Hom}_P(W,W\otimes U(n))$ is identified
with $$g|_{\mathrm {Hom}_P(W,W)}\otimes g|_{U(n)}.$$ Similarly, the
action of $g$ on $\mathrm {Hom}_P\biggl(W,W\otimes
(U(n_1)\oplus\cdots\oplus U(n_m))\biggr)$ is identified with
$$g|_{\mathrm {Hom}_P(W,W)}\otimes (g|_{U(n_1)}\oplus \cdots \oplus
g|_{U(n_m)}).$$ By (ii), there exists a character $\chi:\mathbb
Z\to\overline{\mathbb Q}_\ell^\ast$ such that the representation
$W\otimes \mathscr K_{\chi^{-1}}$ factors through a finite quotient
of $J$. As representations of $J$, we have
$$\mathrm{Hom}_P(W,W)\cong \mathrm{Hom}_P(W\otimes\mathscr
K_{\chi^{-1}},W\otimes\mathscr K_{\chi^{-1}}).$$ It follows that $J$
acts on $\mathrm{Hom}_P(W,W)$ through a finite quotient. So the
matrix of $g$ on $\mathrm{Hom}_P(W,W)$ is diagonalizable. Let
$r=\mathrm{dim}(\mathrm{Hom}_P(W,W))$. Then $g|_{\mathrm
{Hom}_P(W,W)}\otimes g|_{U(n)}$ has $r$ Jordan blocks and each of
them is of size $n\times n$. Similarly, $g|_{\mathrm
{Hom}_P(W,W)}\otimes (g|_{U(n_1)}\oplus \cdots \oplus g|_{U(n_m)})$
has $rm$ Jordan blocks. Among them, $r$ blocks are of size
$n_i\times n_i$ for each $i$. So we must have $m=1$ and $V$ is
indecomposable.
\end{proof}

\begin{corollary}

(i) Let $V$ be an indecomposable $\overline {\mathbb
Q}_\ell$-representation of $I$. Then there exist an irreducible
representation $W$ of $I$ and a positive integer $n$ such that
$$V\cong W \otimes U(n).$$ This factorization is unique up to
isomorphism, and $W$ is isomorphic to any irreducible
sub-representation of $V$.

(ii) Conversely, if $W$ is an irreducible $\overline {\mathbb
Q}_\ell$-representation of $I$, then for any positive integer $n$,
$W\otimes U(n)$ is indecomposable.
\end{corollary}

\begin{proof}

(i) Let $W$ be an irreducible representation of $I$ contained in
$V$. Since $J$ is dense in $I$, $V$ is also an indecomposable
representation of $J$ and $W$ is an irreducible representation of
$J$ contained in $V$. By Lemma 2.2 (i), there exists a positive
integer $n$ such that we have $V\cong W\otimes U(n)$ as
representations of $J$. Again because $J$ is dense in $I$, we have
$V\cong W\otimes U(n)$ as representations of $I$.

(ii) $W$ is an irreducible representation of $J$. By Lemma 2.2
(iii), $W\otimes U(n)$ is indecomposable as a representation of $J$.
It follows that it is indecomposable as a representation of $I$.
\end{proof}

\begin{remark} I don't know whether the analogue of
Lemma 2.2 (ii) holds for irreducible representations of $I$, that
is, whether the following statement is correct: Let $W$ be an
irreducible $\overline{\mathbb Q}_\ell$-representation of $I$. Then
there exists a character
$\chi:\varprojlim\limits_{(N,p)=1}\bbmu_N(k) \to \overline{\mathbb
Q}_\ell^\ast$ such that the representation $W\otimes {\mathcal
K}_{\chi^{-1}}$ factors through a finite discrete quotient of $I$,
where ${\mathcal K}_{\chi^{-1}}$ is the Kummer representation of $I$
associated to the character $\chi^{-1}$. This is true if
$\mathrm{dim}\,W=1$. Indeed, let $\rho:I\to \overline{\mathbb
Q}_\ell^\ast$ be the homomorphism defined by this representation.
Then $\rho(g)$ is an $\ell$-adic unit. So there exists a character
$\chi:\varprojlim\limits_{(N,p)=1}\bbmu_N(k) \to \overline{\mathbb
Q}_\ell^\ast$ such that $\chi(\xi)=\rho(g)$. Then $\rho\chi^{-1}$ is
trivial on $g$ and on an open subgroup of $P$. So $\rho\chi^{-1}$
factors through a finite quotient of $I$.
\end{remark}

\begin{lemma} Any indecomposable representation of
$\varprojlim\limits_{(N,p)=1}\bbmu_N(k)$ is of the form $\mathscr
K_\chi\otimes U(n)$ for some character $\chi:
\varprojlim\limits_{(N,p)=1}\bbmu_N(k)\to \overline {\mathbb
Q}_\ell^\ast$ and some integer $n$.
\end{lemma}

\begin{proof}
Suppose $\rho:\varprojlim\limits_{(N,p)=1}\bbmu_N(k)\to \mathrm
{GL}(n,\overline {\mathbb Q}_\ell)$ is an indecomposable
$n$-dimensional representation. Then $\rho(\zeta)$ is a matrix with
a single Jordan block. Let $\lambda$ be the eigenvalue of
$\rho(\zeta)$ and let $\chi:
\varprojlim\limits_{(N,p)=1}\bbmu_N(k)\to \overline {\mathbb
Q}_\ell^\ast$ be the character defined by $\chi(\zeta)=\lambda$.
Then $\rho$ is isomorphic to $\mathscr K_\chi\otimes U(n)$.
\end{proof}

\begin{lemma} Let $\chi: \varprojlim\limits_{(N,p)=1}\bbmu_N(k) \to
\overline {\mathbb Q}_l^\ast$ be a character, let $\alpha(t)\in
k((t))$ be a formal Laurent series, and let $r$ be a positive
integer prime to $p$. Suppose for any $r$-th root of unity
$\mu\not=1$, the difference $\alpha(\mu t)-\alpha (t)$ is not of the
form $\gamma^p-\gamma$ $\big(\gamma\in k((t))\big)$. Then
$[r]_\ast\bigl({\mathscr L}_\psi(\alpha(t))\otimes {\mathscr
K}_\chi|_{\eta_0}\bigr)$ is an irreducible $\overline {\mathbb
Q}_\ell$-sheaf on $\eta_0$, and for any positive integer $n$,
$[r]_\ast\bigl({\mathscr L}_\psi(\alpha(t))\otimes {\mathscr
K}_\chi|_{\eta_0}\otimes U(n)\bigr)$ is an indecomposable $\overline
{\mathbb Q}_\ell$-sheaf on $\eta_0$.
\end{lemma}

\begin{proof} Choose $T\in \overline {k((t))}$ so that
$T^r=t$. Let ${\mathscr K}'$ be the inverse image of ${\mathscr
K}_\chi$ under the morphism $\mathrm {Spec} k((T))\to {\mathbb
G}_{m,k}$ corresponding to the $k$-algebra homomorphism
$$k[t,1/t]\to k((T)),\quad t\mapsto T,$$ and
let $f: \mathrm {Spec}\, k((T))\to \mathrm {Spec}\, k((t))$ be
the morphism induced by the inclusion $k((t))\hookrightarrow
k((T))$. Then the sheaf $[r]_\ast\bigl({\mathscr
L}_\psi(\alpha(t))\otimes {\mathscr K}_\chi|_{\eta_0}\bigr)$ can be
identified with $f_\ast\bigl({\mathscr L}_\psi(\alpha(T))\otimes
{\mathscr K}'\bigr)$. Let's prove $f_\ast\bigl({\mathscr
L}_\psi(\alpha(T))\otimes {\mathscr K}'\bigr)$ is irreducible. Let
$I=\mathrm {Gal}\bigl(\overline {k((t))}/k((t))\bigr)$ and let
$I'=\mathrm {Gal}\bigl(\overline {k((t))}/k((T))\bigr)$. Note that
$I'$ is a normal subgroup of $I$, and we have a canonical
isomorphism
$$I/I'\stackrel\cong \to \bbmu_r(k), \quad gI' \mapsto
\frac{g(T)}{T}.$$ Let $\rho: I'\to \mathrm {GL}(1, \overline
{\mathbb Q}_\ell)$ be the representation corresponding to the rank
$1$ sheaf ${\mathscr L}_\psi(\alpha(T))\otimes {\mathscr K}'$ on
$\mathrm {Spec}\; k((T))$. Then the representation of $I$
corresponding to the sheaf $f_\ast\bigl({\mathscr
L}_\psi(\alpha(T))\otimes {\mathscr K}'\bigr)$ is $\mathrm
{Ind}_{I'}^I(\rho)$. By the Mackey's criterion (Proposition 23 in
\cite{Se} \S 7.4), to prove $\mathrm {Ind}_{I'}^I(\rho)$ is
irreducible, it suffices to show that for any $g\in I-I'$, the
representation $\rho$ and $\rho^g$ are disjoint, where $\rho^g:
I'\to \mathrm {GL}(1,\overline {\mathbb Q}_\ell)$ is the
representation defined by
$$\rho^g(\sigma)=\rho(g^{-1}\sigma g).$$ Let $\mu=\frac{g(T)}{T}$ and let
$$\tilde g: \mathrm {Spec}\, k((T))\to \mathrm {Spec}\, k((T))$$ be the isomorphism
induced by $g$. Then the sheaf on $\mathrm {Spec}\,k((T))$
corresponding to the representation $\rho^g$ is $$\tilde
g^\ast({\mathscr L}_\psi(\alpha(T))\otimes {\mathscr K}')\cong
{\mathscr L}_\psi(\alpha(\mu T))\otimes \tilde g^\ast {\mathscr
K}'.$$ By our assumption, $\alpha(\mu T)-\alpha(T)$ is not of the
form $\gamma^p-\gamma$ $\big(\gamma\in k((T))\big)$. So ${\mathscr
L}_\psi(\alpha(\mu T))\otimes ({\mathscr L}_\psi(\alpha(T)))^{-1}$
is nontrivial. This representation is wildly ramified. As ${\mathscr
K}'\otimes (\tilde g^\ast {\mathscr K}')^{-1}$ is tamely ramified,
we have
$${\mathscr
L}_\psi(\alpha(\mu T))\otimes ({\mathscr L}_\psi(\alpha(T)))^{-1}
\not\cong {\mathscr K}'\otimes (\tilde g^\ast {\mathscr K}')^{-1}.$$
Therefore,
$${\mathscr L}_\psi(\alpha(\mu T))\otimes \tilde g^\ast {\mathscr K}'\not\cong
{\mathscr L}_\psi(\alpha( T))\otimes {\mathscr K}'.$$ So the two
representations $\rho$ and $\rho^g$ of degree $1$ are disjoint. This
proves $[r]_\ast\bigl({\mathscr L}_\psi(\alpha(t))\otimes {\mathscr
K}_\chi|_{\eta_0}\bigr)$ is irreducible. We have
\begin{eqnarray*}
[r]_\ast\bigl({\mathscr L}_\psi(\alpha(t))\otimes {\mathscr
K}_\chi|_{\eta_0}\otimes U(n)\bigr)&\cong& [r]_\ast\bigl({\mathscr
L}_\psi(\alpha(t))\otimes {\mathscr K}_\chi|_{\eta_0}\otimes
[r]^\ast
U(n)\bigr)\\
&\cong& [r]_\ast\bigl({\mathscr L}_\psi(\alpha(t))\otimes {\mathscr
K}_\chi|_{\eta_0}\bigr)\otimes U(n),
\end{eqnarray*}
and by Corollary 2.3 (ii), $[r]_\ast\bigl({\mathscr
L}_\psi(\alpha(t))\otimes {\mathscr K}_\chi|_{\eta_0}\bigr)\otimes
U(n)$ is indecomposable. So $[r]_\ast\bigl({\mathscr
L}_\psi(\alpha(t))\otimes {\mathscr K}_\chi|_{\eta_0}\otimes
U(n)\bigr)$ is indecomposable.
\end{proof}

\begin{lemma} Let $\alpha(t)=\sum\limits_{i=-s}^\infty a_it^i\in k((t))$ be
a formal Laurent series. For any positive integer $r$, let
$$G(\alpha(t), r)=\{\mu\in \bbmu_r(k)|\alpha(\mu t)-\alpha(t)=\gamma^p-\gamma
\hbox { for some } \gamma\in k((t))\}.$$ Then $G(\alpha(t),r)$ is a
group. Let $\alpha^{(1)}(t)=\sum\limits_{i=-s}^{-1}a_it^i$. Then we
have $G(\alpha(t),r)=G(\alpha^{(1)}(t),r).$
\end{lemma}

\begin{proof} Suppose $\mu_1,\mu_2\in G(\alpha(t),r)$.
Choose $\gamma_1(t),\gamma_2(t)\in k((t))$ such that
\begin{eqnarray*}
\alpha(\mu_1t)-\alpha(t)&=&\gamma_1(t)^p-\gamma_1(t),\\
\alpha(\mu_2t)-\alpha(t)&=& \gamma_2(t)^p-\gamma_2(t).
\end{eqnarray*}
We have
\begin{eqnarray*}
\alpha(\mu_1\mu_2^{-1}t)-\alpha(t)
&=&(\alpha(\mu_1\mu_2^{-1}t)-\alpha(\mu_2^{-1}t))+
(\alpha(\mu_2^{-1}t)-\alpha(t))\\
&=&
(\gamma_1(\mu_2^{-1}t)^p-\gamma_1(\mu_2^{-1}t))-
(\gamma_2(\mu_2^{-1}t)^p-\gamma_2(\mu_2^{-1}t))\\
&=&(\gamma_1(\mu_2^{-1}t)-\gamma_2(\mu_2^{-1}t))^p-
(\gamma_1(\mu_2^{-1}t)-\gamma_2(\mu_2^{-1}t)).
\end{eqnarray*}
Therefore $\mu_1\mu_2^{-1}\in G(\alpha(t),r)$. So $G(\alpha(t),r)$
is a group.

Let $\alpha^{(2)}(t)=\sum\limits_{i=0}^{\infty }a_it^i$. Then for
any $\mu\in\bbmu_r(k)$,  $\alpha^{(2)}(\mu t)-\alpha^{(2)}(t)$ is a
formal power series. By the Hensel Lemma, we can find a formal power
series $\delta\in k[[t]]$ such that
$$\alpha^{(2)}(\mu t)-\alpha^{(2)}(t)=\delta^p-\delta.$$
Then $\alpha(\mu t)-\alpha(t)=\gamma^p-\gamma$ if and only if
$\alpha^{(1)}(\mu
t)-\alpha^{(1)}(t)=(\gamma-\delta)^p-(\gamma-\delta).$ It follows
that $G(\alpha(t),r)=G(\alpha^{(1)}(t),r).$
\end{proof}

\begin{lemma} Let ${\mathscr K}$ be a tamely ramified
$\overline{\mathbb Q}_\ell$-sheaf on $\eta_0=\mathrm {Spec}\,
k((t))$ and let $\theta(t)\in k[[t]]$ be a formal power series of
the form
$$\theta(t)=a_1t+a_2t^2+\cdots$$ with $a_1\not=0$. Denote
by $\theta:\mathrm{Spec}\, k((t))\to \mathrm{Spec}\, k((t))$ the
morphism corresponding to the $k$-algebra homomorphism
$$k((t))\to k((t)), \; t\mapsto \theta(t).$$
Then we have an isomorphism $\theta^\ast {\mathscr K}\cong {\mathscr
K}$.
\end{lemma}

\begin{proof} Set
$$\theta_1(t)=\frac{\theta(t)}{t}=a_1+a_2t+\cdots.$$ For each
integer $N$ with $(N,p)=1$, choose an $N$-th root $\sqrt[N]{t}$ of
$t$ and an $N$-th root $\sqrt[N]{\theta_1(t)}$ of $\theta_1(t)$ in
$\overline {k((t))}$. By the Hensel Lemma, we have
$\sqrt[N]{\theta_1(t)}\in k[[t]]$ for all $N$. We can find a
continuous representation
$\rho:\varprojlim\limits_{(N,p)=1}\bbmu_N(k)\to \mathrm
{GL}(n,\overline {\mathbb Q}_\ell^\ast)$ so that ${\mathscr K}$
corresponds to the Galois representation
$$\mathrm {Gal}\biggl(\overline {k((t))}/k((t))\biggr)\to \mathrm {GL}(n,\overline
{\mathbb Q}_\ell^\ast),\quad \sigma\mapsto
\rho\left(\left(\frac{\sigma(\sqrt[N]{t})}{\sqrt[N]{t}}
\right)_{(N,p)=1}\right).$$ Then $\theta^\ast{\mathscr K}$
corresponds to the Galois representation
$$\mathrm {Gal}\biggl(\overline {k((t))}/k((t))\biggr)\to \mathrm {GL}(n,\overline
{\mathbb Q}_\ell^\ast),\quad \sigma\mapsto
\rho\left(\left(\frac{\sigma(\sqrt[N]{t}\sqrt[N]{\theta_1(t)})}
{\sqrt[N]{t}\sqrt[N]{\theta_1(t)}} \right)_{(N,p)=1}\right).$$ Since
$\sqrt[N]{\theta_1(t)}\in k[[t]]$, we have
$\sigma(\sqrt[N]{\theta_1(t)})=\sqrt[N]{\theta_1(t)}.$ Hence these
two representations are the same.
\end{proof}

\begin{lemma} Let $\chi:\varprojlim\limits_{(N,p)=1}\bbmu_n(k)\to \overline
{\mathbb Q}_\ell^\ast$ be a character. In the notation of Theorem
0.2, as sheaves on $\eta_{\infty'}$,
$${\mathscr L}_\psi(\beta(1/t'))\otimes {\mathscr
K}_{\chi^{-1}}|_{\eta_{\infty'}}\otimes ([s]^\ast{\mathscr
K}_{\chi_2})|_{\eta_{\infty'}}$$ is a direct factor of $[r+s]^\ast
{\mathfrak F}^{(0,\infty')}\biggl([r]_\ast\bigl({\mathscr
L}_\psi(\alpha(t))\otimes {\mathscr K}_\chi|_{\eta_0}
\bigr)\biggr)$.
\end{lemma}

\begin{proof} The equation
$$\frac{\mathrm d}{\mathrm dt}(\alpha(t))+rt^{r-1}t'^{r+s}=0$$ can be written as
$$t'^{r+s}=\frac{1}{t^{r+s}}\left(\frac{sa_{-s}}{r}+
\frac{(s-1)a_{-(s-1)}}{r}t+\cdots\right).$$ Solve this equation, and
write the solution as
$$
t'=\frac{1}{t}(\lambda_0+\lambda_1t+\cdots)=\lambda(t)
$$
with $\lambda_0=\sqrt[r+s]{\frac{sa_{-s}}{r}}.$ The solution is not
unique and different solutions differ by an $(r+s)$-th root of
unity. As long as $\lambda_0$ is fixed to be an $(r+s)$-th root of
$\frac{sa_{-s}}{r}$, for each $i$, $\lambda_i$ depends only on
$a_{-s},a_{-(s-1)},\ldots, a_{-s+i}$. We can write
$$
1/t'=\frac{t}{\lambda_0+\lambda_1t+\cdots}\\
= \nu_1 t+\nu_2t^2+\cdots$$ with $\nu_1=\frac{1}{\lambda_0}$. For
for each $i$, $\nu_i$ depends only on $\lambda_0,\ldots,
\lambda_{i-1}$. Express $t$ in terms of $1/t'$ and write the
solution as
$$
t=\mu_1(1/t')+\mu_2(1/t')^2+\cdots=\mu(1/t')
$$
with $\mu_1=\frac{1}{\nu_1}=\sqrt[r+s]{\frac{sa_{-s}}{r}}.$ For each
$i$, $\mu_i$ depends only on $\nu_1,\ldots, \nu_i$. Substituting
this expression into the equation
$$\beta(1/t')=\alpha(t)+t^rt'^{r+s},$$
we get
\begin{eqnarray*}
\beta(1/t')&=&
\alpha(\mu(1/t'))+(\mu(1/t'))^rt'^{r+s}\\
&=&b_{-s}t'^s+ b_{-(s-1)}t'^{s-1}+\cdots
\end{eqnarray*}
with
$$b_{-s}=\frac{a_{-s}}{\mu_1^s}+\mu_1^r=\frac{a_{-s}(1+\frac{s}{r})}
{\left(\sqrt[r+s]{\frac{sa_{-s}}{r}}\right)^s}.$$ In particular, we
have $b_{-s}\not=0.$ Note that in the above discussion, we need the
fact that $p$ is relatively prime to $r$, $s$ and $r+s$. Moreover,
as long as we fix an $(r+s)$-th root $\sqrt[r+s]{\frac{sa_{-s}}{r}}$
of $\frac{sa_{-s}}{r}$ and let $b_{-s}=\frac{a_{-s}(1+\frac{s}{r})}
{\left(\sqrt[r+s]{\frac{sa_{-s}}{r}}\right)^s}$, then for each $i$,
$b_i$ depends only on $a_{-s},a_{-(s-1)}, \ldots, a_i$. If $\zeta$
is an $(r+s)$-th root of unity, then
$$t=\mu(1/\zeta t'),\; \beta=\beta(1/\zeta t')$$ is also a
solution of the system of equations
\[\left\{ \begin{array}{l}
\alpha(t)+t^rt'^{r+s}=\beta(1/t'),\cr \frac{\mathrm d}{\mathrm
dt}(\alpha(t))+rt^{r-1}t'^{r+s}=0,
\end{array}\right.\]
and all the solutions of this system of equations are of this form.
Note that
$$[r+s]_\ast \bigl({\mathscr
L}_\psi(\beta(1/t'))\otimes {\mathscr
K}_{\chi^{-1}}|_{\eta_{\infty'}}\otimes ([s]^\ast {\mathscr
K}_{\chi_2})|_{\eta_{\infty'}}\bigr)$$ does not depend on the choice
of different solutions of $\beta$.

Write $\alpha(t)=\alpha^{(1)}(t)+\alpha^{(2)}(t),$ where
$$\alpha^{(1)}(t)=\sum_{i=-s}^{-1}a_it^i,\quad
\alpha^{(2)}(t)=\sum_{i=0}^{\infty}a_it^i.$$ We have ${\mathscr
L}_\psi(\alpha(t))\cong {\mathscr L}_\psi(\alpha^{(1)}(t))\otimes
{\mathscr L}_\psi(\alpha^{(2)}(t)),$ and${\mathscr
L}_\psi(\alpha^{(2)}(t))\cong \overline{\mathbb Q}_\ell$ since $k$
is algebraically closed. So we have
$${\mathscr L}_\psi(\alpha(t))\cong {\mathscr
L}_\psi(\alpha^{(1)}(t)).$$ Write
$\beta(1/t')=\beta^{(1)}(1/t')+\beta^{(2)}(1/t'),$ where
$$\beta^{(1)}(1/t')=\sum_{i=-s}^{-1}b_i(1/t')^i,\quad
\beta^{(2)}(1/t')=\sum_{i=0}^{\infty}b_i(1/t')^i.$$ We have
$${\mathscr L}_\psi(\beta(1/t'))\cong {\mathscr
L}_\psi(\beta^{(1)}(1/t')).$$ Note that $\beta^{(1)}(1/t')$ only
depends on $\alpha^{(1)}(t)$. So to prove lemma, we may assume
$\alpha(t)$ is of the form
$$\alpha(t)=\frac{a_{-s}}{t^s}+\frac{a_{-(s-1)}}{t^{s-1}}+\cdots+\frac{a_{-1}}{t}.$$

Let
\begin{eqnarray*}
&&\delta(t')=-\frac{1}{rt'^{r-1}}\frac{\mathrm d}{\mathrm
dt'}(\alpha(t'))=\lambda(t')^{r+s},\\
&&\beta_0(t')=\alpha(t')+t'^r\lambda(t')^{r+s}.
\end{eqnarray*}
Note that $\beta_0(t')$ is the $\beta(t')$ defined Theorem 0.1. We
have
\begin{eqnarray*}
\beta_0(\mu(1/t'))&=&\alpha(\mu(1/t'))+(\mu(1/t'))^r
\lambda(\mu(1/t'))^{r+s}\\
&=& \alpha(\mu(1/t'))+(\mu(1/t'))^r
t'^{r+s}\\
&=&\beta(1/t'),\\
\delta(\mu(1/t'))&=&\lambda(\mu(1/t'))^{r+s}\\
&=& t'^{r+s}.
\end{eqnarray*}
Let $\mu:\eta_{\infty'}\to\eta_{0'}$ and $\delta: \eta_{0'}\to
\eta_{\infty'}$ be the $k$-morphisms corresponding to the
$k$-algebra homomorphisms
\begin{eqnarray*}
k((t'))\to k((1/t')),&& t'\mapsto \mu(1/t'),\\
k((1/t'))\to k((t')),&& t'\mapsto \delta(t'),
\end{eqnarray*}
respectively. Then $\mu$ is an isomorphism and $\delta\circ
\mu=[r+s]$. By Lemma 1.10 (i),
$$\biggl({\mathscr L}_\psi(\beta_0(t'))\otimes {\mathscr
K}_{\chi}\otimes {\mathscr
K}_{\chi_2}\left(\frac{1}{2}s(r+s)a_{-s}t'^s\right)\otimes
G(\chi_2,\psi)\biggr)|_{\eta_{0'}}$$ is a direct factor of
$\delta^\ast {\mathfrak
F}^{(0,\infty')}\biggl([r]_\ast\bigl({\mathscr
L}_\psi(\alpha(t))\otimes {\mathscr K}_\chi\bigr)|_{\eta_0}
\biggr)$. So $$\mu^\ast \biggl(\Big({\mathscr
L}_\psi(\beta_0(t'))\otimes {\mathscr K}_{\chi}\otimes {\mathscr
K}_{\chi_2}\Big(\frac{1}{2}s(r+s)a_{-s}t'^s\Big)\otimes
G(\chi_2,\psi)\Big)|_{\eta_{0'}}\biggr)$$ is a direct factor of
$\mu^\ast \delta^\ast {\mathfrak
F}^{(0,\infty')}\biggl([r]_\ast\bigl({\mathscr
L}_\psi(\alpha(t))\otimes {\mathscr K}_\chi\bigr)|_{\eta_0}
\biggr)$. By Lemma 2.8, we have
\begin{eqnarray*}
\mu^\ast \mathscr K_{\chi}&\cong&\mathscr K_{\chi^{-1}},\\
\mu^\ast {\mathscr
K}_{\chi_2}\left(\frac{1}{2}s(r+s)a_{-s}t'^s\right)&\cong&
[s]^\ast\mathscr K_{\chi_2}.
\end{eqnarray*}
Moreover, we have
\begin{eqnarray*}
\mu^\ast {\mathscr L}_\psi(\beta_0(t'))&\cong& {\mathscr
L}_\psi(\beta_0(\mu(1/t'))) \\&\cong& {\mathscr L}_\psi(\beta(1/t')),\\
\mu^\ast \delta^\ast {\mathfrak
F}^{(0,\infty')}\biggl([r]_\ast\bigl({\mathscr
L}_\psi(\alpha(t))\otimes {\mathscr K}_\chi\bigr)|_{\eta_0}
\biggr)&\cong& (\delta\circ \mu)^\ast{\mathfrak
F}^{(0,\infty')}\biggl([r]_\ast\bigl({\mathscr
L}_\psi(\alpha(t))\otimes {\mathscr K}_\chi\bigr)|_{\eta_0}
\biggr)\\
&\cong& [r+s]^\ast {\mathfrak
F}^{(0,\infty')}\biggl([r]_\ast\bigl({\mathscr
L}_\psi(\alpha(t))\otimes {\mathscr K}_\chi\bigr)|_{\eta_0} \biggr)
\end{eqnarray*}
Since $k$ is algebraically closed, we have
$G(\chi_2,\psi)\cong\overline {\mathbb Q}_\ell$. It follows that
$${\mathscr L}_\psi(\beta(1/t'))\otimes {\mathscr
K}_{\chi^{-1}}|_{\eta_{\infty'}}\otimes ([s]^\ast{\mathscr
K}_{\chi_2})|_{\eta_{\infty'}}$$ is a direct factor of $[r+s]^\ast
{\mathfrak F}^{(0,\infty')}\biggl([r]_\ast\bigl({\mathscr
L}_\psi(\alpha(t))\otimes {\mathscr K}_\chi|_{\eta_0}
\bigr)\biggr)$.
\end{proof}

\begin{lemma} Let $\alpha(t)=\sum_{i=-s}^\infty a_it^i$ be a formal Laurent series
in $k((t))$ and let $\alpha^{(1)}(t)=\sum_{i=-s}^{-1}a_it^i.$ Define
$\beta(1/t')=\sum_{i=-s}^\infty b_i(1/t')^i$ by the system of
equations in Theorem 0.2.

(i) If $\bbmu_d(k)\subset G(\alpha(t),r)$, then $d|r$, $d|s$ and
$a_i=0$ for any $i$ such that $-s\leq i\leq -1$ and $d\not| i$.

(ii) Under the assumption of (i), let $r_0=\frac{r}{d}$,
$s_0=\frac{s}{d}$, $\alpha_0(t)=\sum_{d|i,\; i<0}
a_it^{\frac{i}{d}},$ and let $\beta_0(1/t')$ be defined by the
system of equations
\[\left\{ \begin{array}{l}
\alpha_0(t)+t^{r_0}t'^{r_0+s_0}=\beta_0(1/t'),\\
\frac{\mathrm d}{\mathrm
dt}(\alpha_0(t))+r_0t^{r_0-1}t'^{r_0+s_0}=0.
\end{array}
\right.\] Then $\beta(1/t')-\beta_0(1/t'^d)\in k[[t]].$

(iii) $\bbmu_d(k)\subset G(\alpha(t),r)$ if and only if
$\bbmu_d(k)\subset G(\beta(1/t'),r+s)$.
\end{lemma}

\begin{proof} (i) Since $G(\alpha(t),r)$ is a subgroup of
$\bbmu_r(k)$, we must have $d|r$. Fix a primitive $d$-th root of
unity $\mu_0$ and choose $\gamma\in k((t))$ so that
$$\alpha(\mu_0 t)-\alpha(t)=\gamma^p-\gamma.$$
If the order $\mathrm {ord}_t(\gamma)$ of $\gamma$ with respect to
$t$ is negative, that is, the formal Laurent series $\gamma$
involves negative powers of $t$, then we have
$$\mathrm {ord}_t(\gamma^p-\gamma)=p\,\mathrm {ord}_t(\gamma)\leq -p.$$
But $$\mathrm {ord}_t(\alpha(\mu_0 t)-\alpha(t))\geq -s.$$ So we
have $-s\leq -p$. This contradicts our assumption that $s<p.$ So
$\gamma$ must be a formal power series, and hence the polar part
$$\sum_{i=-s}^{-1}a_i(\mu_0^i-1) t^i$$ of $\alpha(\mu_0
t)-\alpha(t)$ vanishes. Thus $a_i=0$ whenever $i\leq -1$ and $d\not
| i$. As $a_{-s}\not=0$, we have $d|s$.

(ii) Suppose we obtain from the equation
$$\frac{\mathrm d}{\mathrm dt}(\alpha_0(t))+r_0t^{r_0-1}t'^{r_0+s_0}=0$$ the
expressions
$$t'=\lambda(t),\; t=\mu(1/t').$$ Then we have
\begin{eqnarray*}
&&\frac{\mathrm d}{\mathrm dt}(\alpha_0(t))+r_0t^{r_0-1}\lambda(t)^{r_0+s_0}=0,\\
&&\beta_0(1/t')=\alpha_0(\mu(1/t'))+(\mu(1/t'))^{r_0}t'^{r_0+s_0}.
\end{eqnarray*}
So we have
\begin{eqnarray*}
&&\frac{\mathrm d}{\mathrm dt}(\alpha^{(1)}(t))= \frac{\mathrm
d}{\mathrm dt}(\alpha_0(t^d))=\frac{\mathrm d\alpha_0}{\mathrm
dt}(t^d) dt^{d-1}\\
&=& -r_0(t^d)^{r_0-1}\lambda(t^d)^{r_0+s_0}dt^{d-1} = -r t^{r-1}
((\lambda(t^d))^{1/d})^{r+s}.
\end{eqnarray*}
Hence from the equation
$$\frac{\mathrm d}{\mathrm dt}(\alpha^{(1)}(t))+rt^{r-1}t'^{r+s}=0,$$
we get the expressions
$$t'=(\lambda(t^d))^{1/d},\; t=(\mu(1/t'^d))^{1/d},$$
and we have
\begin{eqnarray*}
\beta_0(1/t'^d)&=& \alpha_0(\mu(1/t'^d))+(\mu(1/t'^d))^{r_0}
(t'^d)^{r_0+s_0}\\&=&\alpha^{(1)}((\mu(1/t'^d))^{1/d})+((\mu(1/t'^d))^{1/d})^{r}
t'^{r+s}.
\end{eqnarray*}
So $\alpha^{(1)}(t)$ and $\beta_0(1/t'^d)$ satisfy the same system
of equations for $\alpha(t)$ and $\beta(1/t')$:
\[\left\{ \begin{array}{l}
\alpha^{(1)}(t)+t^rt'^{r+s}=\beta_0(1/t'^d),\\
\frac{\mathrm d}{\mathrm dt}(\alpha^{(1)}(t))+rt^{r-1}t'^{r+s}=0.
\end{array}
\right.\] Since $\alpha^{(1)}(t)$ is the polar part of $\alpha(t)$,
$\beta(1/t)$ and $\beta_0(1/t'^d)$ have the same polar part by the
discussion in the proof of Lemma 2.9. Hence
$\beta(1/t)-\beta_0(1/t'^d)\in k[[t]]$.

(iii) Suppose $\bbmu_d(k)\subset G(\alpha(t),r)$. Then the
conclusion of (ii) holds. It is clear that $\bbmu_d(k)\subset
G(\beta_0(1/t'^d),r+s)$. Since $\beta(1/t')-\beta_0(1/t'^d)\in
k[[t]]$, we have $$G(\beta(1/t'),r+s)=G(\beta_0(1/t'^d),r+s)$$ by
Lemma 2.7. So $\bbmu_d(k)\subset G(\beta(1/t'),r+s)$.

We have
\begin{eqnarray*}
\frac{\mathrm d}{\mathrm dt'}(\beta(1/t'))&=&
\frac{\mathrm d}{\mathrm dt'}(\alpha(t)+t^rt'^{r+s})\\
&=& \frac{\mathrm d}{\mathrm dt}(\alpha(t))\frac{\mathrm dt}{\mathrm
dt'}+rt^{r-1}t'^{r+s}
\frac{\mathrm dt}{\mathrm dt'}+(r+s) t^rt'^{r+s-1}\\
&=&\left(\frac{\mathrm d}{\mathrm
dt}(\alpha(t))+rt^{r-1}t'^{r+s}\right)
\frac{\mathrm dt}{\mathrm dt'}+(r+s) t^rt'^{r+s-1}\\
&=& (r+s) t^rt'^{r+s-1}.
\end{eqnarray*}
So $\alpha(t)$ is a solution of the system of equations
\[\left\{ \begin{array}{l}
\beta(1/t')-t^rt'^{r+s}=\alpha(t),\\
\frac{\mathrm d}{\mathrm dt'}(\beta(1/t'))-(r+s) t^rt'^{r+s-1}=0.
\end{array}
\right.\] (The above discussion is nothing but a re-proof of the
involutivity of the Legendre transformation. See \cite{A} Chapter 3,
\S 14 C.) Reversing the role of $\alpha$ and $\beta$ in the above
discussion, we see that the condition $\bbmu_d(k)\subset
G(\beta(1/t'),r+s)$ implies the condition $\bbmu_d(k)\subset
G(\alpha(t),r)$.
\end{proof}

We are now ready to prove Theorems 0.1-0.4.

\begin{proof}[Proof of Theorem 0.2] Notation as in Theorem 0.2.
We may assume $\alpha(t)$ is of the form
$$\alpha(t)=\frac{a_{-s}}{t^s}+\frac{a_{-(s-1)}}{t^{s-1}}+\cdots+\frac{a_{-1}}{t}$$
with $a_{-s}\not =0$. First we prove that the general case of
Theorem 0.2 follows from the special case where $\alpha(t)$ is
assumed to have the property $G(\alpha(t),r)=1$. Let
$G(\alpha(t),r)=\bbmu_d(k)$ and let $r_0,s_0,
\alpha_0(t),\beta_0(1/t')$ be defined as in Lemma 2.10. We have
$\alpha(t)=\alpha_0(t^d)$, $G(\alpha_0(t),r_0)=\{1\},$ and
$\beta(1/t')=\beta_0(1/t'^d)$. Note that $[d]_\ast{\mathscr K}$ is a
lisse $\overline {\mathbb Q}_\ell$-sheaf on ${\mathbb G}_m$ tamely
ramified at $0$ and at $\infty$. One can verify $\mathrm{inv}^\ast
[d]_\ast {\mathscr K} \cong  [d]_\ast \mathrm {inv}^\ast {\mathscr
K}.$ Applying the special case of Theorem 0.2 to $\alpha_0(t)$, we
get
\begin{eqnarray*}
&&{\mathfrak F}^{(0,\infty')}\biggl([r_0]_\ast\bigl({\mathscr
L}_\psi(\alpha_0(t))\otimes ([d]_\ast {\mathscr
K})|_{\eta_0}\bigr)\biggr)\\&\cong& [r_0+s_0]_\ast \bigl({\mathscr
L}_\psi(\beta_0(1/t'))\otimes (\mathrm {inv}^\ast [d]_\ast {\mathscr
K})|_{\eta_{\infty'}}\otimes ([s_0]^\ast {\mathscr
K}_{\chi_2})|_{\eta_{\infty'}}\bigr)\\
&\cong& [r_0+s_0]_\ast \bigl({\mathscr L}_\psi(\beta_0(1/t'))\otimes
([d]_\ast\mathrm {inv}^\ast  {\mathscr K})|_{\eta_{\infty'}}\otimes
([s_0]^\ast {\mathscr K}_{\chi_2})|_{\eta_{\infty'}}\bigr).
\end{eqnarray*}
So we have
\begin{eqnarray*}
&&{\mathfrak F}^{(0,\infty')}\biggl([r]_\ast\bigl({\mathscr
L}_\psi(\alpha(t))\otimes {\mathscr
K}|_{\eta_0}\bigr)\biggr)\\
&\cong& {\mathscr
F}^{(0,\infty')}\biggl([r_0]_\ast[d]_\ast\bigl([d]^\ast {\mathscr
L}_\psi(\alpha_0(t))\otimes {\mathscr
K}|_{\eta_0}\bigr)\biggr)\\
&\cong& {\mathfrak F}^{(0,\infty')}\biggl([r_0]_\ast\bigl({\mathscr
L}_\psi(\alpha_0(t))\otimes ([d]_\ast{\mathscr
K})|_{\eta_0}\bigr)\biggr) \\
&\cong &[r_0+s_0]_\ast \bigl({\mathscr L}_\psi(\beta_0(1/t'))\otimes
([d]_\ast \mathrm {inv}^\ast {\mathscr K})|_{\eta_{\infty'}}\otimes
([s_0]^\ast {\mathscr K}_{\chi_2})|_{\eta_{\infty'}}\bigr) \\
&\cong& [r_0+s_0]_\ast [d]_\ast \bigl([d]^\ast {\mathscr
L}_\psi(\beta_0(1/t'))\otimes (\mathrm{inv}^\ast {\mathscr
K})|_{\eta_{\infty'}}\otimes ([d]^\ast [s_0]^\ast {\mathscr
K}_{\chi_2})|_{\eta_{\infty'}}\bigr)\\
&\cong & [r+s]_\ast \bigl({\mathscr L}_\psi(\beta(1/t'))\otimes
(\mathrm {inv}^\ast {\mathscr K})|_{\eta_{\infty'}}\otimes ([s]^\ast
{\mathscr K}_{\chi_2})|_{\eta_{\infty'}}\bigr)
\end{eqnarray*}
This prove the general case of Theorem 0.2.

Frow now on, we assume $G(\alpha(t),r)=\{1\}$. By Lemma 2.10 (iii),
we have $G(\beta(1/t'),r+s)=\{1\}$. If ${\mathscr K}={\mathscr
K}_\chi$ for some character
$$\chi:\varprojlim\limits_{(N,p)=1}\bbmu_N(k)\to \overline{\mathbb
Q}_\ell^\ast,$$ then $[r]_\ast\bigl({\mathscr
L}_\psi(\alpha(t))\otimes {\mathscr K}_\chi|_{\eta_0} \bigr)$ is
irreducible by Lemma 2.6, and hence ${\mathscr
F}^{(0,\infty')}\biggl([r]_\ast\bigl({\mathscr
L}_\psi(\alpha(t))\otimes {\mathscr K}_\chi|_{\eta_0} \bigr)\biggr)$
is also irreducible by \cite{L} 2.4.3 (i) c) and (ii) a). By Lemma
2.9, ${\mathscr L}_\psi(\beta(1/t'))\otimes {\mathscr
K}_{\chi^{-1}}|_{\eta_{\infty'}} \otimes ([s]^\ast {\mathscr
K}_{\chi_2})|_{\eta_{\infty'}}$ is a direct factor of $[r+s]^\ast
{\mathfrak F}^{(0,\infty')}\biggl([r]_\ast\bigl({\mathscr
L}_\psi(\alpha(t))\otimes {\mathscr K}_\chi|_{\eta_0}
\bigr)\biggr)$. For each $\mu\in\bbmu_{r+s}(k)$, $\beta(1/\mu t')$
is a solution of the system of equations in Theorem 0.2. So
${\mathscr L}_\psi(\beta(1/\mu t'))\otimes {\mathscr
K}_{\chi^{-1}}|_{\eta_{\infty'}} \otimes ([s]^\ast {\mathscr
K}_{\chi_2})|_{\eta_{\infty'}}$ is also a direct factor of
$[r+s]^\ast {\mathfrak
F}^{(0,\infty')}\biggl([r]_\ast\bigl({\mathscr
L}_\psi(\alpha(t))\otimes {\mathscr K}_\chi|_{\eta_0}
\bigr)\biggr)$. As $G(\beta(1/t),r+s)=\{1\}$, the rank $1$
representations of $\mathrm
{Gal}(\bar\eta_{\infty'}/\eta_{\infty'})$ corresponding to
${\mathscr L}_\psi(\beta(1/\mu t'))\otimes {\mathscr
K}_{\chi^{-1}}|_{\eta_{\infty'}} \otimes ([s]^\ast {\mathscr
K}_{\chi_2})|_{\eta_{\infty'}}$ are disjoint as $\mu$ go over
$\bbmu_{r+s}(k)$. By Lemma 2.1 (ii), $[r+s]^\ast {\mathfrak
F}^{(0,\infty')}\biggl([r]_\ast\bigl({\mathscr
L}_\psi(\alpha(t))\otimes {\mathscr K}_\chi|_{\eta_0} \bigr)\biggr)$
is semisimple. Moreover, the rank of $[r+s]^\ast {\mathfrak
F}^{(0,\infty')}\biggl([r]_\ast\bigl({\mathscr
L}_\psi(\alpha(t))\otimes {\mathscr K}_\chi|_{\eta_0} \bigr)\biggr)$
is $r+s$ by \cite{L} 2.4.3 (i) b). So we have
$$[r+s]^\ast {\mathfrak
F}^{(0,\infty')}\biggl([r]_\ast\bigl({\mathscr
L}_\psi(\alpha(t))\otimes {\mathscr K}_\chi|_{\eta_0}
\bigr)\biggr)\cong \bigoplus_{\mu\in\bbmu_{r+s}(k)} {\mathscr
L}_\psi(\beta(1/\mu t'))\otimes {\mathscr
K}_{\chi^{-1}}|_{\eta_{\infty'}} \otimes ([s]^\ast {\mathscr
K}_{\chi_2})|_{\eta_{\infty'}}.$$ This is the isotypical
decomposition of $[r+s]^\ast {\mathfrak
F}^{(0,\infty')}\biggl([r]_\ast\bigl({\mathscr
L}_\psi(\alpha(t))\otimes {\mathscr K}_\chi|_{\eta_0}
\bigr)\biggr)$. The projection
$$[r+s]^\ast {\mathfrak
F}^{(0,\infty')}\biggl([r]_\ast\bigl({\mathscr
L}_\psi(\alpha(t))\otimes {\mathscr K}_\chi|_{\eta_0}
\bigr)\biggr)\to {\mathscr L}_\psi(\beta(1/t'))\otimes {\mathscr
K}_{\chi^{-1}}|_{\eta_{\infty'}} \otimes ([s]^\ast {\mathscr
K}_{\chi_2})|_{\eta_{\infty'}}$$ defines a homomorphism
$${\mathfrak
F}^{(0,\infty')}\biggl([r]_\ast\bigl({\mathscr
L}_\psi(\alpha(t))\otimes {\mathscr K}_\chi|_{\eta_0}
\bigr)\biggr)\to [r+s]_\ast\bigl({\mathscr
L}_\psi(\beta(1/t'))\otimes {\mathscr
K}_{\chi^{-1}}|_{\eta_{\infty'}} \otimes ([s]^\ast {\mathscr
K}_{\chi_2})|_{\eta_{\infty'}}\bigr).$$ If we apply $[r+s]^\ast$ to
this homomorphism, we get an isomorphism. So this homomorphism
itself is an isomorphism. This prove Theorem 0.2 in the case where
$\mathscr K\cong \mathscr K_\chi$.

By Lemma 2.5, to prove Theorem 0.2 in general, we may assume
$${\mathscr K}= {\mathscr K}_\chi\otimes U(n)$$ for some character
$\chi:\varprojlim\limits_{(N,p)=1}\bbmu_N(k)\to \overline{\mathbb
Q}_\ell^\ast$ and some integer $n$. Since $G(\alpha(t),r)=\{1\}$,
the sheaf $[r]_\ast\bigl({\mathscr L}_\psi(\alpha(t))\otimes
{\mathscr K}|_{\eta_0}\bigr)$ is indecomposable by Lemma 2.6, and it
contains $[r]_\ast\bigl({\mathscr L}_\psi(\alpha(t))\otimes
{\mathscr K}_\chi|_{\eta_0}\bigr)$ as an irreducible subsheaf. By
\cite{L} 2.4.3 (i) c) and (ii) a), ${\mathfrak
F}^{(0,\infty')}\biggl([r]_\ast\bigl({\mathscr
L}_\psi(\alpha(t))\otimes {\mathscr K}|_{\eta_0} \bigr)\biggr)$ is
indecomposable, and it contains  ${\mathfrak
F}^{(0,\infty')}\biggl([r]_\ast\bigl({\mathscr
L}_\psi(\alpha(t))\otimes {\mathscr K}_\chi|_{\eta_0} \bigr)\biggr)$
as an irreducible subsheaf. So $[r+s]_\ast\bigl({\mathscr
L}_\psi(\beta(1/t'))\otimes (\mathrm {inv}^\ast {\mathscr
K}_\chi)|_{\eta_{\infty'}}\otimes ([s]^\ast{\mathscr
K}_{\chi_2})|_{\eta_{\infty'}}\bigr)$ is an irreducible subsheaf of
the indecomposable sheaf ${\mathfrak
F}^{(0,\infty')}\biggl([r]_\ast\bigl({\mathscr
L}_\psi(\alpha(t))\otimes {\mathscr K}|_{\eta_0} \bigr)\biggr)$. By
Corollary 2.3 (i), we have $${\mathfrak
F}^{(0,\infty')}\biggl([r]_\ast\bigl({\mathscr
L}_\psi(\alpha(t))\otimes {\mathscr K}|_{\eta_0} \bigr)\biggr) \cong
[r+s]_\ast\bigl({\mathscr L}_\psi(\beta(1/t'))\otimes (\mathrm
{inv}^\ast {\mathscr K}_\chi)|_{\eta_{\infty'}}\otimes
([s]^\ast{\mathscr K}_{\chi_2})|_{\eta_{\infty'}}\bigr)\otimes
U(n')$$ for some integer $n'$. By \cite{L} 2.4.3 (i) b), the rank of
${\mathfrak F}^{(0,\infty')}\biggl([r]_\ast\bigl({\mathscr
L}_\psi(\alpha(t))\otimes {\mathscr K}|_{\eta_0} \bigr)\biggr)$ is
$$r\,\mathrm{rank}({\mathscr K})+s\,\mathrm {rank}(\mathscr K)=n(r+s),$$
whereas the rank of $[r+s]_\ast\bigl({\mathscr
L}_\psi(\beta(1/t'))\otimes (\mathrm {inv}^\ast {\mathscr
K}_\chi)|_{\eta_{\infty'}}\otimes ([s]^\ast{\mathscr
K}_{\chi_2})|_{\eta_{\infty'}}\bigr)\otimes U(n')$ is $n'(r+s)$. So
we must have $n'=n$, and
\begin{eqnarray*}
&&{\mathfrak F}^{(0,\infty')}\biggl([r]_\ast\bigl({\mathscr
L}_\psi(\alpha(t))\otimes {\mathscr K}|_{\eta_0} \bigr)\biggr)\\
&\cong & [r+s]_\ast\bigl({\mathscr L}_\psi(\beta(1/t'))\otimes
(\mathrm {inv}^\ast {\mathscr K}_\chi)|_{\eta_{\infty'}}\otimes
([s]^\ast{\mathscr K}_{\chi_2})|_{\eta_{\infty'}}\bigr)\otimes U(n)\\
&\cong&  [r+s]_\ast\bigl({\mathscr L}_\psi(\beta(1/t'))\otimes
((\mathrm {inv}^\ast {\mathscr K}_\chi)\otimes [r+s]^\ast
U(n))|_{\eta_{\infty'}}\otimes ([s]^\ast{\mathscr
K}_{\chi_2})|_{\eta_{\infty'}}\bigr)\\
&\cong & [r+s]_\ast\bigl({\mathscr L}_\psi(\beta(1/t'))\otimes
(\mathrm {inv}^\ast ({\mathscr K}_\chi\otimes
U(n)))|_{\eta_{\infty'}}\otimes
([s]^\ast{\mathscr K}_{\chi_2})|_{\eta_{\infty'}}\bigr) \\
&\cong& [r+s]_\ast\bigl({\mathscr L}_\psi(\beta(1/t'))\otimes
(\mathrm {inv}^\ast {\mathscr K})|_{\eta_{\infty'}}\otimes
([s]^\ast{\mathscr K}_{\chi_2})|_{\eta_{\infty'}}\bigr).
\end{eqnarray*}
This verifies the formula in Theorem 0.2.\end{proof}

\begin{proof}[Proof of Theorem 0.1 (i)] Here $k$ is not necessarily algebraically
closed. Notation as in Theorem 0.1 (i). By Lemma 1.10 (i), we have a
projection
$$\delta^\ast {\mathfrak
F}^{(0,\infty')}\biggl([r]_\ast\bigl({\mathscr
L}_\psi(\alpha(t))\otimes {\mathscr K}\bigr)|_{\eta_0} \biggr)\to
\Big({\mathscr L}_\psi(\beta(t'))\otimes {\mathscr K}\otimes
{\mathscr K}_{\chi_2}\Big(\frac{1}{2}s(r+s)a_{-s}t'^s\Big)\otimes
G(\chi_2,\psi)\Big)|_{\eta_{0'}}.$$ It defines a homomorphism
$${\mathfrak
F}^{(0,\infty')}\biggl([r]_\ast\bigl({\mathscr
L}_\psi(\alpha(t))\otimes {\mathscr K}\bigr)|_{\eta_0} \biggr)\to
\biggl(\delta_\ast\Big({\mathscr L}_\psi(\beta(t'))\otimes {\mathscr
K}\otimes {\mathscr
K}_{\chi_2}\Big(\frac{1}{2}s(r+s)a_{-s}t'^s\Big)\otimes
G(\chi_2,\psi)\Big)\biggr)|_{\eta_{\infty'}}.$$ To prove it is an
isomorphism, we may make base extension from $k$ to its algebraic
closure $\bar k$. Over algebraically closed field, it is an
isomorphism by the above proof of the $\mathscr K=\mathscr K_\chi$
case of Theorem 0.2.
\end{proof}

\begin{proof}[Proof of Theorems 0.1 (ii)-(iii), 0.3 and 0.4.]
The proof of Theorems 0.1 (ii) and 0.3 for the case $s>r$, and the
proof of Theorems 0.1 (iii) and 0.4 for the case $s<r$ are similar
to the above proof of Theorems 0.1 (i) and 0.2 by using the same
group theoretical argument and Lemma 1.10 (ii)-(iii).

In the notation of Theorem 0.3, if $s<r$, then the slopes of
$[r]_\ast\bigl({\mathscr L}_\psi(\alpha(1/t))\otimes{\mathscr
K}\bigr)|_{\eta_\infty}$ are $\frac{s}{r}\leq 1$. (Confer \cite{K3}
1.13.) By \cite{L} 2.4.3 (iii) b), we then have
$${\mathfrak F}^{(\infty,\infty')}\biggl([r]_\ast\bigl({\mathscr
L}_\psi(\alpha(1/t))\otimes {\mathscr
K}\bigr)|_{\eta_\infty}\biggr)=0.$$ In the notation of Theorem 0.4,
if $s>r$, then the slopes of $[r]_\ast\bigl({\mathscr
L}_\psi(\alpha(1/t))\otimes{\mathscr K}\bigr)|_{\eta_\infty}$ are
$\frac{s}{r}\geq 1$. By \cite{L} 2.4.3 (iii) b), we then have
$${\mathfrak F}^{(\infty,0')}\biggl([r]_\ast\bigl({\mathscr
L}_\psi(\alpha(1/t))\otimes {\mathscr
K}\bigr)|_{\eta_\infty}\biggr)=0.$$ Similarly we can prove Theorem
0.1 (ii) for the case $s<r$ and Theorem 0.1 (iii) for the case
$s>r$.
\end{proof}

\begin{proof}[Proof of Proposition 0.5] By (b), there exists a finite Galois
extension $L$ of $k((t))$ such that $\rho$ factors through
$I_0=\mathrm {Gal}(L/k((t)))$. The maximal tamely ramified extension
of $k((t))$ contained in $L$ is of the form $k((T))$, where
$T=\sqrt[N]{t}$ is an $N$-th root of $t$ for some positive integer
$N$ prime to $p$. Let $P_0=\mathrm {Gal}(L/k((T)))$. $\mathrm
{Res}_{P_0}^{I_0}(V)$ is semisimple. By (a), $\mathrm
{Res}_{P_0}^{I_0}(V)$ factors through $P_0^{\mathrm
{ab}}=P_0/[P_0,P_0]$. So $\mathrm {Res}_{P_0}^{I_0}(V)$ is a direct
sum of one dimensional representations. Also by (a), the image of
$P_0$ under each of these one dimensional representations is
contained in
$$\bbmu_p(\overline {\mathbb Q}_\ell)=\{\mu\in \overline {\mathbb
Q}_\ell^\ast|\mu^p=1\}.$$ The additive character $\psi$ is an
isomorphism from $\mathbb Z/p$ to $\bbmu_p(\overline {\mathbb
Q}_\ell)$. By the Artin-Schreier theory, each one dimensional direct
summand of $\mathrm {Res}_{P_0}^{I_0}(V)$ is an Artin-Schreier
representation, that is, a representation of $P_0$ corresponding to
a sheaf $\mathscr L_\psi(\alpha(T))$ on $\mathrm {Spec}\, k((T))$
for some Laurent series $\alpha(T)\in k((T))$. Let
$\{\alpha_1,\ldots, \alpha_m\}$ be a set of $\alpha_i(T)\in k((T))$
so that the representations of $P_0$ corresponding to the sheaves
$\mathscr L_\psi(\alpha_i(T))$ on $\mathrm {Spec}\, k((T))$ give all
one dimensional sub-representations of $\mathrm
{Res}_{P_0}^{I_0}(V)$, and so that $\mathscr
L_\psi(\alpha_i(T))\not\cong \mathscr L_\psi(\alpha_j(T))$ for
$i\not =j$. Moreover, we can choose $\alpha_i(T)$ so that in their
formal power series expansion, only non-positive powers of $T$ are
involved. Let
$$V=\bigoplus_{i=1}^m V_i$$ be the isotypical
decomposition of $\mathrm {Res}_{P_0}^{I_0}(V)$ so that for each
$i$, $V_i$ is a direct sum of several copies of the representation
corresponding to the sheaf $\mathscr L_\psi(\alpha_i(T))$. As in the
beginning of the proof of Lemma 2.2 (i), $I_0$ acts on the set
$\{V_1,\ldots, V_m\}$, and if $I_1$ is the stabilizer of $V_1$ under
this action, then $P_0\subset I_1$, $V_1$ is a representation of
$I_1$, and $V=\mathrm {Ind}_{I_1}^{I_0}(V_1).$ The field $L^{I_1}$
is a subfield of $L^{P_0}=k((T))$ and we have $L^{I_1}=k((S))$,
where $S=T^{\frac{N}{r}}$ is an $r$-th root of $t$ for some positive
integer $r$ dividing $N$. Choose $g\in I_1$ so that
$\frac{g(T)}{T}=\mu$ is a primitive $\frac{N}{r}$-th root of unity.
Let $\rho_{V_1}:P_0\to \mathrm {GL}(V_1)$ be the homomorphism
defined by the representation $V_1$ and let $\rho_{V_1}^g:P_0\to
\mathrm {GL}(V_1)$ be defined by
$\rho^g_{V_1}(\sigma)=\rho_{V_1}(g^{-1}\sigma g)$ for any $\sigma\in
P_0$. Let $\rho_{\alpha_1}:P_0\to \overline {\mathbb Q}_\ell^\ast$
be the character corresponding to the sheaf $\mathscr
L_\psi(\alpha_1(T))$ on $\mathrm {Spec}\, k((T))$, and let
$\rho_{\alpha_1}^g:P_0\to \overline {\mathbb Q}_\ell^\ast$ be the
character defined by
$\rho^g_{\alpha_1}(\sigma)=\rho_{\alpha_1}(g^{-1}\sigma g)$ for any
$\sigma\in P_0$. Then the representation defined by $\rho_{V_1}$
(resp. $\rho_{V_1}^g$) is a direct sum of several copies of the
irreducible representation defined by $\rho_{\alpha_1}$ (resp.
$\rho_{\alpha_1}^g$). But $g:V_1\to V_1$ is an isomorphism from the
representation $\rho_{V_1}^g$ to the representation $\rho_{V_1}$. It
follows that the representations defined by $\rho_{\alpha_1}$ and
$\rho_{\alpha_1}^g$ are isomorphic. As $\rho_{\alpha_1}$ (resp.
$\rho_{\alpha_1}^g$) corresponds to the sheaf $\mathscr
L_\psi(\alpha_1(T))$ (resp. $\mathscr L_\psi(\alpha_1(\mu T))$) on
$\mathrm {Spec}\, k((T))$, it follows that $\mathscr
L_\psi(\alpha_1(T))\cong \mathscr L_\psi(\alpha_1(\mu T))$. So there
exists $\gamma(T)\in k((T))$ such that
$$\alpha_1(\mu T)-\alpha_1(T)=\gamma(T)^p-\gamma(T).$$ Write
$$\alpha_1(T)=\frac{a_{-s_1}}{T^{s_1}}+\frac{a_{-(s_1-1)}}{T^{s_1-1}}
+\cdots+\frac{a_{-1}}{T}$$ with $a_{-s_1}\not=0$. Denote by $s(-)$
the Swan conductor of a Galois representation. We have
$$s_1=s(\mathscr
L_\psi(\alpha_1(T)))\leq s(V_1)=s(\mathrm
{Ind}_{I_1}^{I_0}(V_1))=s(V)=s.$$ By (c), we have $s_1<p$. The same
argument as in the proof of Lemma 2.10 (i) shows that $\gamma(T)$
must be a formal power series, and $a_i=0$ whenever $\frac{N}{r}\not
| i$. Let
$$\alpha(S)=\sum_{\frac{N}{r}|i}a_iS^{\frac{ir}{N}}.$$ Then we have
$\alpha_1(T)=\alpha(T^{\frac{N}{r}}).$ Denote the rank one
representation of $I_1$ corresponding to the lisse $\overline
{\mathbb Q}_\ell$-sheaf $\mathscr L_\psi(\alpha(S))$ on $\mathrm
{Spec}\, k((S))$ also by $\mathscr L_\psi(\alpha(S))$. Then $\mathrm
{Res}_{P_0}^{I_1}(\mathscr L_\psi(\alpha(S)))$ is the representation
of $P_0$ corresponding to the sheaf $\mathscr L_\psi(\alpha_1(T))$
on $\mathrm {Spec}\, k((T))$. It follows that the representation
$\mathrm {Res}_{P_0}^{I_1} (V_1\otimes \mathscr
L_\psi(\alpha(S))^{-1})$ is a direct sum of the trivial
representation. Let $\mathscr K=V_1\otimes \mathscr
L_\psi(\alpha(S))^{-1}.$ Then $\mathscr K$ is a tamely ramified
representation of $I_1$. We have
$$V=\mathrm {Ind}_{I_1}^{I_0}(V_1)=\mathrm {Ind}_{I_1}^{I_0}(
\mathscr L_\psi(\alpha(S))\otimes \mathscr K).$$ As $V$ is
irreducible, $\mathscr L_\psi(\alpha(S))\otimes \mathscr K$ is
irreducible, and hence $\mathscr K$ is irreducible. So $\mathscr
K=\mathscr K_\chi$ for some character
$\chi:\varprojlim\limits_{(N,p)=1} \bbmu_N(k) \to \overline{\mathbb
Q}_\ell^\ast$ of finite order.
\end{proof}

\begin{remark} If the statement in Remark 2.4
holds, then by Corollary 2.3 (i) and Proposition 0.5, the following
generalization of Proposition 0.5 holds:  Let $\rho: I\to \mathrm
{GL}(V)$ be an indecomposable $\overline {\mathbb
Q}_\ell$-representation. Suppose $\rho(P^p[P,P])=1$ and suppose the
Swan conductor of $\rho$ is less than $p$. Then $V$ is isomorphic to
the representation corresponding to the lisse $\overline{\mathbb
Q}_\ell$-sheaf $[r]_\ast(\mathscr L_\psi(\alpha(t))\otimes \mathscr
K)$ on $\eta_0=\mathrm {Spec}\,k((t))$ for some positive integer $r$
prime to $p$, some Laurent series $\alpha(t)$ and some tamely
ramified sheaf $\mathscr K$ on $\eta_0$.
\end{remark}

\section{Local monodromy of hypergeometric sheaves}

Let $k$ be an algebraically closed field of characteristic $p$, and
let
$$\lambda_1,\ldots,\lambda_n,\rho_1,\ldots,
\rho_m:\varprojlim_{(N,p)=1}\bbmu_N(k)\to\overline {\mathbb
Q}_\ell^\ast$$ be characters. In \cite{K1} 8.2, Katz defines the
hypergeometric complex $\mathrm
{Hyp}(!,\psi;\lambda_1,\ldots,\lambda_n;\rho_1,\ldots, \rho_m)$ in
$D_c^b({\mathbb G}_{m,k},\overline {\mathbb Q}_\ell)$ using the
multiplicative convolution. Let $j:{\mathbb G}_{m,k}\to {\mathbb
A}_k^1$ be the canonical open immersion, and let $\mathrm {inv}:
{\mathbb G}_{m,k}\to {\mathbb G}_{m,k}$ be the morphism defined by
$t\mapsto \frac{1}{t}$. Recall the following properties of the
hypergeometric complex (\cite{K1} 8.2.2 (2), 8.2.4, 8.2.5, 8.1.12).

(a) If $n=1$ and $m=0$, we have
$$\mathrm {Hyp} (!,\psi;\lambda_1;\emptyset)\cong j^\ast{\mathscr L}_\psi
\otimes {\mathscr K}_{\lambda_1}[1].$$

(b) For any character
$\lambda:\varprojlim\limits_{(N,p)=1}\bbmu_N(k)\to\overline {\mathbb
Q}_\ell^\ast$, we have
$$\mathrm
{Hyp}(!,\psi;\lambda\lambda_1,\ldots,\lambda\lambda_n;\lambda\rho_1,\ldots,
\lambda\rho_m)\cong \mathscr K_\lambda \otimes \mathrm
{Hyp}(!,\psi;\lambda_1,\ldots,\lambda_n;\rho_1,\ldots, \rho_m).$$

(c) In the case where $\lambda_n=1$, we have
$$\mathrm
{Hyp}(!,\psi;\lambda_1,\ldots,\lambda_{n-1},1 ;\rho_1,\ldots,
\rho_m)\cong j^\ast{\mathfrak F}(j_!\mathrm {inv}^\ast \mathrm
{Hyp}(!,\psi;\lambda_1,\ldots,\lambda_{n-1};\rho_1,\ldots,
\rho_m)).$$

(d) We have $$\mathrm {inv}^\ast \mathrm
{Hyp}(!,\psi;\lambda_1,\ldots,\lambda_n;\rho_1,\ldots, \rho_m) \cong
\mathrm {Hyp}(!, \psi^{-1}; \rho_1^{-1},\ldots,
\rho_m^{-1};\lambda_1^{-1},\ldots,\lambda_n^{-1}).$$

The above properties can also be used to define $\mathrm
{Hyp}(!,\psi;\lambda_1,\ldots,\lambda_n;\rho_1,\ldots, \rho_m)$.
Indeed, by properties (b) and (c), we have
\begin{eqnarray*}
&&\mathrm {Hyp}(!,\psi;\lambda_1,\ldots,\lambda_n;\rho_1,\ldots,
\rho_m)\\
&\cong& \mathscr K_{\lambda_n}\otimes \mathrm
{Hyp}(!,\psi;\lambda_1\lambda_n^{-1},\ldots,\lambda_{n-1}\lambda_n^{-1},1
;\rho_1\lambda_n^{-1},\ldots,
\rho_m\lambda_n^{-1})\\
&\cong &\mathscr K_{\lambda_n}\otimes j^\ast{\mathfrak F}(j_!\mathrm
{inv}^\ast \mathrm
{Hyp}(!,\psi;\lambda_1\lambda_n^{-1},\ldots,\lambda_{n-1}\lambda_n^{-1};
\rho_1\lambda_n^{-1},\ldots,\rho_m\lambda_n^{-1})).
\end{eqnarray*}
This reduces the definition of $\mathrm
{Hyp}(!,\psi;\lambda_1,\ldots,\lambda_n;\rho_1,\ldots, \rho_m)$ to
the case where $n=0$. By properties (d), (b) and (c), for the case
where $n=0$, we have
\begin{eqnarray*}
&&\mathrm {Hyp}(!,\psi;\emptyset;\rho_1,\ldots,
\rho_m)\\
&\cong& \mathrm {inv}^\ast \mathrm
{Hyp}(!,\psi^{-1};\rho_1^{-1},\ldots, \rho_m^{-1};\emptyset)\\
&\cong& \mathrm {inv}^\ast(\mathscr K_{\rho_m^{-1}}\otimes \mathrm
{Hyp}(!,\psi^{-1};\rho_1^{-1}\rho_m,\ldots,
\rho_{m-1}^{-1}\rho_m, 1;\emptyset))\\
&\cong & \mathrm {inv}^\ast(\mathscr K_{\rho_m^{-1}}\otimes
j^\ast\mathfrak F_{\psi^{-1}}(j_!\mathrm {inv}^\ast \mathrm
{Hyp}(!,\psi^{-1};\rho_1^{-1}\rho_m,\ldots,
\rho_{m-1}^{-1}\rho_m;\emptyset))),
\end{eqnarray*}
where $\mathfrak F_{\psi^{-1}}$ denotes the Fourier transformation
defined by the additive character $\psi^{-1}$. This reduces the
definition of $\mathrm
{Hyp}(!,\psi;\lambda_1,\ldots,\lambda_n;\rho_1,\ldots, \rho_m)$ to
that in the case where $n=1$ and $m=0$, for which we can use the
property (a).

\medskip
To prove Proposition 0.7 in the Introduction, we first prove
Proposition 0.8.

\begin{proof}[Proof of Proposition 0.8.] When $n=1$, our assertion
follows from the fact that
$$\mathrm {Hyp} (!,\psi;\lambda_1;\emptyset)[-1]\cong
j^\ast{\mathscr L}_\psi \otimes {\mathscr K}_{\lambda_1}.$$ Suppose
the lemma holds if there are $n-1$ characters $\lambda$'s. To prove
the lemma for $n$ characters $\lambda$'s, we use the fact that
$$
\mathrm {Hyp}(!,\psi;\lambda_1,\ldots, \lambda_n;\emptyset)\cong
{\mathscr K}_{\lambda_n}\otimes \mathrm
{Hyp}(!,\psi;\lambda_1\lambda_n^{-1},\ldots,
\lambda_{n-1}\lambda_n^{-1}, 1;\emptyset)$$ to reduce to the case
where $\lambda_n=1$. We have
$$\mathrm {Hyp}(!,\psi;\lambda_1,\ldots, \lambda_{n-1}, 1;\emptyset)
=j^\ast{\mathfrak F}(j_!\mathrm {inv}^\ast \mathrm
{Hyp}(!,\psi;\lambda_1,\ldots,\lambda_{n-1};\emptyset)).$$ By the
induction hypothesis, we have
$$\bigl(\mathrm {Hyp}(!,\psi; \lambda_1,\ldots,
\lambda_{n-1};\emptyset)[-1]\bigr)|_{\eta_\infty}\cong
[n-1]_\ast\bigl({\mathscr L}_\psi((n-1)t)\otimes {\mathscr
K}_{\lambda_1\cdots\lambda_{n-1}}\otimes{\mathscr
K}_{{\chi_2}^{n-2}}\bigr)|_{\eta_\infty}.$$ So
\begin{eqnarray}
\bigl(j_!\mathrm {inv}^\ast\mathrm {Hyp}(!,\psi; \lambda_1,\ldots,
\lambda_n;\emptyset)[-1]\bigr)|_{\eta_0}\cong
[n-1]_\ast\left({\mathscr L}_\psi\left(\frac{n-1}{t}\right)\otimes
\mathrm{inv}^\ast ({\mathscr
K}_{\lambda_1\cdots\lambda_{n-1}}\otimes{\mathscr
K}_{{\chi_2}^{n-2}})\right)|_{\eta_0}.
\end{eqnarray}
By the induction hypothesis, $\mathrm {inv}^\ast \mathrm
{Hyp}(!,\psi;\lambda_1,\ldots,\lambda_{n-1};\emptyset)[-1]$ is an
irreducible lisse sheaf of rank $n-1$ on $\mathbb G_{m,k}$. Combined
with (1), we get
$$j_!\mathrm {inv}^\ast \mathrm
{Hyp}(!,\psi;\lambda_1,\ldots,\lambda_{n-1};\emptyset)\cong
j_{!\ast} \mathrm {inv}^\ast \mathrm
{Hyp}(!,\psi;\lambda_1,\ldots,\lambda_{n-1};\emptyset),$$ and
$j_!\mathrm {inv}^\ast \mathrm
{Hyp}(!,\psi;\lambda_1,\ldots,\lambda_{n-1};\emptyset)$ is an
irreducible perverse sheaf of type $(\mathrm T_3)$ in \cite{L}
1.4.2. So ${\mathfrak F}(j_!\mathrm {inv}^\ast \mathrm
{Hyp}(!,\psi;\lambda_1,\ldots,\lambda_{n-1};\emptyset))$ is an
irreducible perverse sheaf on $\mathbb A_k^1$ of type $(\mathrm
T_3)$ by \cite{L} 1.4.2.1 (ii). By the induction hypothesis, we have
\begin{eqnarray}
\begin{array}{ccl}
\bigl(j_!\mathrm {inv}^\ast \mathrm
{Hyp}(!,\psi;\lambda_1,\ldots,\lambda_{n-1};\emptyset)[-1]\bigr)|_{\eta_{\infty}}&\cong&
\mathrm {inv}^\ast \left(\bigl(\mathrm
{Hyp}(!,\psi;\lambda_1,\ldots,\lambda_{n-1};\emptyset)[-1]\bigr)|_{\eta_0}\right)\\
&\cong& \mathrm {inv}^\ast \left(\bigoplus_{\lambda} \bigl({\mathscr
K}_\lambda\otimes U(\mathrm {mult}_0(\lambda))\bigr)|_{\eta_0}\right)\\
&\cong&  \bigoplus_{\lambda} \bigl(\mathrm {inv}^\ast ({\mathscr
K}_\lambda\otimes U(\mathrm
{mult}_0(\lambda)))\bigr)|_{\eta_\infty},
\end{array}
\end{eqnarray}
where $\mathrm{mult}_0(\lambda)$ denotes the number of times that
$\lambda$ appears in $\lambda_1,\ldots, \lambda_{n-1}$. By (2) and
\cite{L} 2.3.1.3 (i),  ${\mathfrak F}(j_!\mathrm {inv}^\ast \mathrm
{Hyp}(!,\psi;\lambda_1,\ldots,\lambda_{n-1};\emptyset))$ is lisse on
$\mathbb G_{m,k}$. It follows that $\mathrm
{Hyp}(!,\psi;\lambda_1,\ldots, \lambda_{n-1}, 1;\emptyset)[-1]$ is
an irreducible lisse sheaf on ${\mathbb G}_{m,k}$, and we have
$${\mathfrak F}(j_!\mathrm {inv}^\ast \mathrm
{Hyp}(!,\psi;\lambda_1,\ldots,\lambda_{n-1};\emptyset))=j_{!\ast}\mathrm
{Hyp}(!,\psi;\lambda_1,\ldots,\lambda_{n-1},1;\emptyset),$$ and the
long exact sequence in \cite{L} 2.3.2 gives rise to a short exact
sequence
\begin{eqnarray}
\begin{array}{cclcc}
0&\to& \bigl(\mathrm {Hyp}(!,\psi;\lambda_1,\ldots, \lambda_{n-1},
1;\emptyset)[-1]\bigr)^{\mathrm{Gal}(\bar\eta_0/\eta_0)}_{\bar\eta_0}&&\\
&\to& \bigl(\mathrm {Hyp}(!,\psi;\lambda_1,\ldots, \lambda_{n-1},
1;\emptyset)[-1]\bigr)_{\bar\eta_0}&&\\
&\to&{\mathfrak F}^{(\infty,0)}\big(\bigl(j_!\mathrm {inv}^\ast
\mathrm
{Hyp}(!,\psi;\lambda_1,\ldots,\lambda_{n-1};\emptyset)[-1]\bigr)|_{\eta_{\infty}}\big)
&\to& 0.
\end{array}
\end{eqnarray}
We have
\begin{eqnarray}
\begin{array}{cl}
&{\mathfrak F}^{(\infty,0)}\big(\bigl(j_!\mathrm {inv}^\ast \mathrm
{Hyp}(!,\psi;\lambda_1,\ldots,\lambda_{n-1};\emptyset)[-1]\bigr)|_{\eta_{\infty}}\big)\\
\cong& {\mathfrak F}^{(\infty,0)}\big(\bigoplus_{\lambda}
\bigl(\mathrm {inv}^\ast({\mathscr K}_\lambda\otimes U(\mathrm
{mult}_0(\lambda)))\bigr)|_{\eta_\infty}\big) \\
\cong& \bigoplus_{\lambda}\bigl({\mathscr K}_\lambda\otimes
U(\mathrm {mult}_0(\lambda))\bigr)|_{\eta_0}.
\end{array}
\end{eqnarray}
Using the formula in \cite{L} 2.3.1.1 (i) and (iii), one can check
\begin{eqnarray}
&&\mathrm {rank} \bigl(\mathrm {Hyp}(!,\psi;\lambda_1,\ldots,
\lambda_{n-1},
1;\emptyset)[-1]\bigr)^{\mathrm{Gal}(\bar\eta_0/\eta_0)}_{\bar\eta_0}=1,\\
&&\mathrm {rank} \bigl(\mathrm {Hyp}(!,\psi;\lambda_1,\ldots,
\lambda_{n-1}, 1;\emptyset)[-1]\bigr)=n.
\end{eqnarray}
It follows from (3)-(6) that
$$
\bigl(\mathrm {Hyp}(!,\psi;\lambda_1,\ldots, \lambda_{n-1},
1;\emptyset)[-1]\bigr)|_{\eta_0}\cong U(\mathrm
{mult}_0(1))|_{\eta_0}\bigoplus \bigoplus_{\lambda\not=1}
\bigl({\mathscr K}_\lambda\otimes U(\mathrm
{mult}_0(\lambda))\bigr)|_{\eta_0},$$ where $\mathrm{mult}_0(1)$ is
the number of trivial characters in $\lambda_1,\ldots,
\lambda_{n-1},1$.

By \cite{L} 2.3.3.1 (iii), we have
\begin{eqnarray*}
\bigl(\mathrm
{Hyp}(!,\psi;\lambda_1,\ldots,\lambda_{n-1},1;\emptyset)[-1]\bigr)|_{\eta_\infty}
&\cong&  {\mathfrak F}^{(0,\infty)}\bigl(\bigl(j_!\mathrm {inv}^\ast
\mathrm
{Hyp}(!,\psi;\lambda_1,\ldots,\lambda_{n-1};\emptyset)[-1]\bigr)|_{\eta_0}\bigr)\\
&&\qquad\bigoplus {\mathscr
F}^{(\infty,\infty)}\bigl(\bigl(j_!\mathrm {inv}^\ast \mathrm
{Hyp}(!,\psi;\lambda_1,\ldots,\lambda_{n-1};\emptyset)[-1]\bigr)|_{\eta_\infty}\bigr).
\end{eqnarray*}
By (2) and \cite{L} 2.4.3 (iii) b), we have
$${\mathfrak F}^{(\infty,\infty)}\bigl(\bigl(j_!\mathrm {inv}^\ast \mathrm
{Hyp}(!,\psi;\lambda_1,\ldots,\lambda_{n-1};\emptyset)[-1]\bigr)|_{\eta_\infty}\bigr)=0.$$
Using the formula in Theorem 0.2 and (1), we find
$${\mathfrak F}^{(0,\infty)}\bigl(\bigl(j_!\mathrm {inv}^\ast\mathrm {Hyp}(!,\psi; \lambda_1,\ldots,
\lambda_n;\emptyset)[-1]\bigr)|_{\eta_0}\bigr)\cong
[n]_\ast\big({\mathscr L}_\psi(nt)\otimes {\mathscr
K}_{\lambda_1\cdots\lambda_{n-1}}\otimes{\mathscr
K}_{{\chi_2}^{n-1}}\big)|_{\eta_\infty}.$$ It follows that
$$\bigl(\mathrm {Hyp}(!,\psi; \lambda_1,\ldots,
\lambda_{n-1},1;\emptyset)[-1]\bigr)|_{\eta_\infty}\cong
[n]_\ast\big({\mathscr L}_\psi(nt)\otimes {\mathscr
K}_{\lambda_1\cdots\lambda_{n-1}}\otimes{\mathscr
K}_{{\chi_2}^{n-1}}\big)|_{\eta_\infty}.$$
\end{proof}

Using Proposition 0.8 and the fact that
$$\mathrm {Hyp}(!,\psi;\emptyset; \rho_1,\ldots, \rho_m)
\cong \mathrm{inv}^\ast \mathrm {Hyp}(!,\psi^{-1};
\rho_1^{-1},\ldots, \rho_m^{-1};\emptyset),$$ we get the following.

\begin{corollary} Suppose $p$ is relatively prime to $2, 3, \ldots,
m$. Then $\mathrm {Hyp}(!,\psi;\emptyset; \rho_1,\ldots,
\rho_m)[-1]$ is an irreducible lisse sheaf of rank $m$ on ${\mathbb
G}_{m,k}$, and we have
\begin{eqnarray*}
\bigl(\mathrm {Hyp}(!,\psi; \emptyset; \rho_1,\ldots,
\rho_m)[-1]\bigr)|_{\eta_\infty}&\cong& \bigoplus_{\rho}
\bigl({\mathscr K}_\rho\otimes U(\mathrm
{mult}_\infty(\rho))\bigr)|_{\eta_\infty},\\
\bigl(\mathrm {Hyp}(!,\psi;\emptyset; \rho_1,\ldots,
\rho_m)[-1]\bigr)|_{\eta_0}&\cong& [m]_\ast\left({\mathscr
L}_\psi\left(-\frac{m}{t}\right)\otimes {\mathscr
K}_{\rho_1\cdots\rho_m}\otimes{\mathscr
K}_{{\chi_2}^{m-1}}\right)|_{\eta_0}.
\end{eqnarray*}
\end{corollary}

We now prove Proposition 0.7.

\begin{proof}[Proof of Proposition 0.7]
Using the fact that
$$\mathrm {Hyp}(!,\psi;\lambda_1,\ldots, \lambda_n; \rho_1,\ldots, \rho_m)
\cong \mathrm{inv}^\ast \mathrm {Hyp}(!,\psi^{-1};
\rho_1^{-1},\ldots, \rho_m^{-1};\lambda_1^{-1},\ldots,
\lambda_n^{-1}),$$ we can deduce (i) from (ii). Let's prove (ii).
When $n=0$, this follows from Corollary 3.1. Suppose $1\leq n<m$ and
suppose the lemma holds if there are $n-1$ characters $\lambda$'s.
To prove the lemma for $n$ characters $\lambda$'s, we use the fact
that
$$
\mathrm {Hyp}(!,\psi;\lambda_1,\ldots,
\lambda_n;\rho_1,\ldots,\rho_m)\cong {\mathscr K}_{\lambda_n}\otimes
\mathrm {Hyp}(!,\psi;\lambda_1\lambda_n^{-1},\ldots,
\lambda_{n-1}\lambda_n^{-1},
1;\rho_1\lambda_n^{-1},\ldots,\rho_m\lambda_n^{-1})$$ to reduce to
the case where $\lambda_n=1$. Then by the disjointness of
$\lambda$'s and $\rho$'s, we have $\rho_i\not =1$ for all $1\leq
i\leq m$. We have
$$\mathrm {Hyp}(!,\psi;\lambda_1,\ldots, \lambda_{n-1}, 1;\rho_1,\ldots, \rho_m)
=j^\ast{\mathfrak F}(j_!\mathrm {inv}^\ast \mathrm
{Hyp}(!,\psi;\lambda_1,\ldots,\lambda_{n-1};\rho_1,\ldots,
\rho_m)).$$ By the induction hypothesis, $\mathrm {inv}^\ast \mathrm
{Hyp}(!,\psi;\lambda_1,\ldots,\lambda_{n-1};\rho_1,\ldots,
\rho_m)[-1]$ is an irreducible lisse sheaf of rank $m$ on $\mathbb
G_{m,k}$, and
\begin{eqnarray*}
\bigl(\mathrm {inv}^\ast\mathrm
{Hyp}(!,\psi;\lambda_1,\ldots,\lambda_{n-1};\rho_1,\ldots,\rho_m)[-1]\bigr)|_{\eta_0}
\cong \bigoplus_{\rho} \bigl(\mathrm {inv}^\ast({\mathscr
K}_\rho\otimes U(\mathrm {mult}_\infty(\rho)))\bigr)|_{\eta_0}.
\end{eqnarray*}
Since $\rho_i\not=1$ for all $1\leq i\leq m$, this implies that
$$j_!\mathrm {inv}^\ast \mathrm
{Hyp}(!,\psi;\lambda_1,\ldots,\lambda_{n-1};\rho_1,\ldots,\rho_m)\cong
j_{!\ast} \mathrm {inv}^\ast \mathrm
{Hyp}(!,\psi;\lambda_1,\ldots,\lambda_{n-1};\rho_1,\ldots,\rho_m),$$
and  $j_!\mathrm {inv}^\ast \mathrm
{Hyp}(!,\psi;\lambda_1,\ldots,\lambda_{n-1};\rho_1,\ldots,\rho_m)$
is an irreducible perverse sheaf of type $(\mathrm T_3)$ in \cite{L}
1.4.2. So ${\mathfrak F}(j_!\mathrm {inv}^\ast \mathrm
{Hyp}(!,\psi;\lambda_1,\ldots,\lambda_{n-1};\rho_1,\ldots,\rho_m))$
is an irreducible perverse sheaf on $\mathbb A_k^1$ of type
$(\mathrm T_3)$ by \cite{L} 1.4.2.1 (ii). By the induction
hypothesis, we have
\begin{eqnarray}
\begin{array}{cl}
&\bigl(j_!\mathrm {inv}^\ast \mathrm
{Hyp}(!,\psi;\lambda_1,\ldots,\lambda_{n-1};\rho_1,\ldots,\rho_m)[-1]\bigr)|_{\eta_{\infty}}
\\
\cong& \mathrm {inv}^\ast\biggl([m-(n-1)]_\ast\biggl({\mathscr
L}_\psi\left(-\frac{m-(n-1)}{t}\right)\otimes {\mathscr
K}_{\lambda_1^{-1}\cdots\lambda_{n-1}^{-1}\rho_1\cdots
\rho_m}\otimes{\mathscr
K}_{{\chi_2}^{n+m-2}}\biggr)\\
&\qquad \qquad \bigoplus \bigoplus_{\lambda} \big({\mathscr
K}_\lambda\otimes U(\mathrm{mult}_0(\lambda))\big)\bigg)|_{\eta_\infty}\\
\cong & [m-(n-1)]_\ast\biggl({\mathscr L}_\psi(-(m-(n-1))t)\otimes
\mathrm{inv}^\ast\big( {\mathscr
K}_{\lambda_1^{-1}\cdots\lambda_{n-1}^{-1}\rho_1\cdots
\rho_m}\otimes{\mathscr
K}_{{\chi_2}^{n+m-2}}\big)\biggr)|_{\eta_\infty}\\&\qquad \bigoplus
\bigoplus_{\lambda} \bigl(\mathrm{inv}^\ast({\mathscr
K}_\lambda\otimes U(\mathrm
{mult}_0(\lambda)))\bigr)|_{\eta_\infty}.
\end{array}
\end{eqnarray} So all the
breaks of $\bigl(j_!\mathrm {inv}^\ast \mathrm
{Hyp}(!,\psi;\lambda_1,\ldots,\lambda_{n-1};\rho_1,\ldots,
\rho_m)[-1]\bigr)|_{\eta_\infty}$ are $<1$. By \cite{L} 2.3.1.3 (i),
${\mathfrak F}(j_!\mathrm {inv}^\ast \mathrm
{Hyp}(!,\psi;\lambda_1,\ldots,\lambda_{n-1};\rho_1,\ldots, \rho_m))$
is lisse on ${\mathbb G}_{m,k}$. It follows that $\mathrm
{Hyp}(!,\psi;\lambda_1,\ldots, \lambda_{n-1}, 1;\rho_1,\ldots,
\rho_m)[-1]$ is an irreducible lisse sheaf on ${\mathbb G}_{m,k}$,
we have
$${\mathfrak F}(j_!\mathrm {inv}^\ast \mathrm
{Hyp}(!,\psi;\lambda_1,\ldots,\lambda_{n-1};\rho_1,\ldots,
\rho_m))\cong j_{!\ast} \mathrm {Hyp}(!,\psi;\lambda_1,\ldots,
\lambda_{n-1}, 1;\rho_1,\ldots, \rho_m),$$ and the long exact
sequence in \cite{L} 2.3.2 gives rise to the short exact sequence
\begin{eqnarray}
\begin{array}{cclcc}
0&\to& \bigl(\mathrm {Hyp}(!,\psi;\lambda_1,\ldots, \lambda_{n-1},
1;\rho_1,\ldots,\rho_m)[-1]\bigr)^{\mathrm{Gal}(\bar\eta_0/\eta_0)}_{\bar\eta_0}&&
\\&\to& \bigl(\mathrm {Hyp}(!,\psi;\lambda_1,\ldots, \lambda_{n-1},
1;\rho_1,\ldots,\rho_m)[-1]\bigr)_{\bar\eta_0}&&
\\&\to&{\mathfrak F}^{(\infty,0)}\bigl(\bigl(j_!\mathrm {inv}^\ast \mathrm
{Hyp}(!,\psi;\lambda_1,\ldots,\lambda_{n-1};\rho_1,\ldots,\rho_m)[-1]
\bigr)|_{\eta_{\infty}}\bigr) &\to& 0.
\end{array}
\end{eqnarray}
Using the formulas in \cite{L} 2.3.1.1 (i) and (iii), one can check
\begin{eqnarray}
&&\mathrm {rank} \bigl(\mathrm {Hyp}(!,\psi;\lambda_1,\ldots,
\lambda_{n-1},
1;\rho_1,\ldots,\rho_m)[-1]\bigr)^{\mathrm{Gal}(\bar\eta_0/\eta_0)}_{\bar\eta_0}=1,\\
&& \mathrm {rank} \bigl(\mathrm {Hyp}(!,\psi;\lambda_1,\ldots,
\lambda_{n-1}, 1;\rho_1,\ldots,\rho_m)[-1]\bigr)=m.
\end{eqnarray}
Using the formula in Theorem 0.4 and (7), we find that
\begin{eqnarray}
\begin{array}{cl}
&{\mathfrak F}^{(\infty,0)}\bigl(\bigl(j_!\mathrm {inv}^\ast \mathrm
{Hyp}(!,\psi;\lambda_1,\ldots,\lambda_{n-1};\rho_1,\ldots,\rho_m)[-1]\bigr)|_{\eta_{\infty}}\bigr)\\
\cong& [m-n]_\ast\left({\mathscr
L}_\psi\left(-\frac{m-n}{t}\right)\otimes {\mathscr
K}_{\lambda_1^{-1}\cdots\lambda_n^{-1}\rho_1\cdots
\rho_m}\otimes{\mathscr
K}_{{\chi_2}^{n+m-1}}\right)|_{\eta_0}\bigoplus\bigoplus_{\lambda}
\bigl({\mathscr K}_\lambda\otimes U(\mathrm
{mult}_0(\lambda))\bigr)|_{\eta_0}.
\end{array}
\end{eqnarray}
It follows from (8)-(11) that
\begin{eqnarray*}
&&\bigl(\mathrm {Hyp}(!,\psi;\lambda_1,\ldots, \lambda_{n-1},
1;\rho_1,\ldots,\rho_m)[-1]\bigr)|_{\eta_0}\\&\cong&
[m-n]_\ast\left({\mathscr L}_\psi\left(-\frac{m-n}{t}\right)\otimes
{\mathscr K}_{\lambda_1^{-1}\cdots\lambda_n^{-1}\rho_1\cdots
\rho_m}\otimes{\mathscr K}_{{\chi_2}^{n+m-1}}\right)|_{\eta_0}
\\&& \qquad \bigoplus U(\mathrm {mult}_0(1))|_{\eta_0}\bigoplus
\bigoplus_{\lambda\not=1} \bigl({\mathscr K}_\lambda\otimes
U(\mathrm {mult}_0(\lambda))\bigr)|_{\eta_0}.
\end{eqnarray*}

By \cite{L} 2.3.3.1 (iii), we have
\begin{eqnarray*}
&&\bigl(\mathrm
{Hyp}(!,\psi;\lambda_1,\ldots,\lambda_{n-1},1;\rho_1,\ldots, \rho_m)[-1]\bigr)|_{\eta_\infty}\\
&\cong&  {\mathfrak F}^{(0,\infty)}\big(\bigl(j_!\mathrm {inv}^\ast
\mathrm
{Hyp}(!,\psi;\lambda_1,\ldots,\lambda_{n-1};\rho_1,\ldots,\rho_m)[-1]\bigr)|_{\eta_0}\big)\\
&& \quad \bigoplus {\mathscr
F}^{(\infty,\infty)}\bigl(\big(j_!\mathrm {inv}^\ast \mathrm
{Hyp}(!,\psi;\lambda_1,\ldots,\lambda_{n-1};\rho_n,\ldots,
\rho_m)[-1]\bigr)|_{\eta_\infty}\bigr).
\end{eqnarray*}
Since $\bigl(j_!\mathrm {inv}^\ast \mathrm
{Hyp}(!,\psi;\lambda_1,\ldots,\lambda_{n-1};\rho_1,\ldots,\rho_m)[-1]\bigr)|_{\eta_{\infty}}$
has breaks $<1$, we have
$${\mathfrak F}^{(\infty,\infty)}\bigl(\big(j_!\mathrm {inv}^\ast \mathrm
{Hyp}(!,\psi;\lambda_1,\ldots,\lambda_{n-1};\rho_1,\ldots,\rho_m)[-1]\big)|_{\eta_\infty}\bigr)=0$$
by \cite{L} 2.4.3 (iii) b). We have
\begin{eqnarray*}
&&{\mathfrak F}^{(0,\infty)}\bigl(\big(j_!\mathrm {inv}^\ast \mathrm
{Hyp}(!,\psi;\lambda_1,\ldots,\lambda_{n-1};\rho_1,\ldots,\rho_m)[-1]\bigr)|_{\eta_0}
\bigr)\\
&\cong& {\mathfrak F}^{(0,\infty)}\bigl(\bigoplus_{\rho}\bigl(
\mathrm {inv}^\ast ({\mathscr K}_\rho\otimes U(\mathrm
{mult}_\infty(\rho)))\bigr)|_{\eta_0}\bigr) \\
&\cong & \bigoplus_{\rho} \bigl(\mathrm {\mathscr K}_\rho\otimes
U(\mathrm {mult}_\infty(\rho))\bigr)|_{\eta_\infty}.
\end{eqnarray*}
It follows that
$$\bigl(\mathrm {Hyp}(!,\psi; \lambda_1,\ldots,
\lambda_n;\rho_1,\ldots, \rho_m)[-1]\bigr)|_{\eta_\infty}\cong
\bigoplus_{\rho} \bigl({\mathscr K}_\rho\otimes U(\mathrm
{mult}_\infty(\rho))\bigr)|_{\eta_\infty}.$$

Finally, we prove (iii). First consider the case $\lambda_n=1$. By
the disjointness of $\lambda$'s and $\rho$'s, we have $\rho_i\not=1$
for all $i$. We have
$$\mathrm
{Hyp}(!,\psi;\lambda_1,\ldots,
\lambda_{n-1},1;\rho_1,\ldots,\rho_n)\cong j^\ast{\mathfrak
F}(j_!\mathrm {inv}^\ast \mathrm
{Hyp}(!,\psi;\lambda_1,\ldots,\lambda_{n-1};\rho_1,\ldots,
\rho_n)).$$ By (ii), $\mathrm {inv}^\ast \mathrm
{Hyp}(!,\psi;\lambda_1,\ldots,\lambda_{n-1};\rho_1,\ldots,
\rho_n)[-1]$ is an irreducible lisse sheaf on $\mathbb G_{m,k}$ and
\begin{eqnarray}
\bigl(\mathrm {inv}^\ast \mathrm
{Hyp}(!,\psi;\lambda_1,\ldots,\lambda_{n-1};\rho_1,\ldots,
\rho_n)[-1]\bigr)|_{\eta_0}\cong \bigoplus_{\rho} \bigl(\mathrm
{inv}^\ast ({\mathscr K}_\rho\otimes U(\mathrm
{mult}_\infty(\rho)))\bigr)|_{\eta_0}.
\end{eqnarray}
Since $\rho_i\not=1$ for all $i$, this implies that
$$j_!\mathrm {inv}^\ast \mathrm
{Hyp}(!,\psi;\lambda_1,\ldots,\lambda_{n-1};\rho_1,\ldots,
\rho_n)\cong j_{!\ast}\mathrm {inv}^\ast \mathrm
{Hyp}(!,\psi;\lambda_1,\ldots,\lambda_{n-1};\rho_1,\ldots,
\rho_n),$$ and $j_!\mathrm {inv}^\ast \mathrm
{Hyp}(!,\psi;\lambda_1,\ldots,\lambda_{n-1};\rho_1,\ldots, \rho_n)$
is an irreducible perverse sheaf of type $(\mathrm T_3)$ in \cite{L}
1.4.2. So ${\mathfrak F}(j_!\mathrm {inv}^\ast \mathrm
{Hyp}(!,\psi;\lambda_1,\ldots,\lambda_{n-1};\rho_1,\ldots, \rho_n))$
is an irreducible perverse sheaf on $\mathbb A_k^1$ of type
$(\mathrm T_3)$ by \cite{L} 1.4.2.1 (ii). By (ii), we have
\begin{eqnarray}
\begin{array}{cl}
&\bigl(j_!\mathrm {inv}^\ast \mathrm {Hyp}(!,\psi; \lambda_1,\ldots,
\lambda_{n-1};\rho_1,\ldots, \rho_n)[-1]\bigr)|_{\eta_\infty} \\
\cong & \big({\mathscr L}_\psi(-t)\otimes \mathrm{inv}^\ast
{\mathscr K}_{\lambda_1^{-1}\cdots\lambda_{n-1}^{-1}\rho_1\cdots
\rho_n}\big)|_{\eta_\infty}\bigoplus \bigoplus_{\lambda}\bigl(
\mathrm{inv}^\ast({\mathscr K}_\lambda\otimes U(\mathrm
{mult}_0(\lambda)))\bigr)|_{\eta_\infty}.
\end{array}
\end{eqnarray} By (13) and
\cite{L} 2.3.1.2, ${\mathfrak F}(j_!\mathrm {inv}^\ast \mathrm
{Hyp}(!,\psi;\lambda_1,\ldots,\lambda_{n-1};\rho_1,\ldots, \rho_n))$
is lisse on $\mathbb G_{m,k}-\{1\}$. Therefore $\mathrm
{Hyp}(!,\psi; \lambda_1,\ldots, \lambda_{n-1},1;\rho_1,\ldots,
\rho_n)[-1]$ is an irreducible lisse sheaf when restricted to
${\mathbb G}_{m,k}-\{1\}$. Using \cite{L} 2.3.1.1 (i), one can show
its rank is $n$.

Let $\mathrm {tran}:\mathbb A_k^1\to \mathbb A_k^1$ be the
translation by $1$. By the definition of the Fourier transformation,
we have
\begin{eqnarray}
\begin{array}{cl}
& {\mathfrak F}(j_!\mathrm {inv}^\ast \mathrm
{Hyp}(!,\psi;\lambda_1,\ldots,\lambda_{n-1};\rho_1,\ldots,
\rho_n))\\
\cong& \mathrm {tran}^\ast {\mathfrak F}({\mathscr L}_\psi\otimes
j_!\mathrm {inv}^\ast \mathrm
{Hyp}(!,\psi;\lambda_1,\ldots,\lambda_{n-1};\rho_1,\ldots, \rho_n)).
\end{array}
\end{eqnarray}
Note that ${\mathscr L}_\psi\otimes j_!\mathrm {inv}^\ast \mathrm
{Hyp}(!,\psi;\lambda_1,\ldots,\lambda_{n-1};\rho_1,\ldots, \rho_n)$
is also an irreducible perverse sheaf on $\mathbb A_k^1$ of type
$(\mathrm T_3)$. Hence ${\mathfrak F}({\mathscr L}_\psi\otimes
j_!\mathrm {inv}^\ast \mathrm
{Hyp}(!,\psi;\lambda_1,\ldots,\lambda_{n-1};\rho_1,\ldots, \rho_n))$
is an irreducible perverse sheaf  on $\mathbb A_k^1$ of type
$(\mathrm T_3)$. So we have
\begin{eqnarray}
\begin{array}{cl}
&\mathscr H^{-1} \bigl({\mathfrak F}(\mathscr L_\psi\otimes
j_!\mathrm {inv}^\ast \mathrm
{Hyp}(!,\psi;\lambda_1,\ldots,\lambda_{n-1};\rho_1,\ldots,
\rho_n))\bigr)_{\bar 0} \\
\cong&\mathscr  H^{-1} \bigl( {\mathfrak F}(\mathscr L_\psi\otimes
j_!\mathrm {inv}^\ast \mathrm
{Hyp}(!,\psi;\lambda_1,\ldots,\lambda_{n-1};\rho_1,\ldots,
\rho_n))\bigr)_{\bar \eta_0}^{\mathrm {Gal}(\bar\eta_0/\eta_0)},
\end{array}
\end{eqnarray}
and the long exact sequence in \cite{L} 2.3.2 gives rise to the
short exact sequence
\begin{eqnarray}
\begin{array}{ccccc}
0&\to& \mathscr H^{-1} \bigl({\mathfrak F}(\mathscr L_\psi\otimes
j_!\mathrm {inv}^\ast \mathrm
{Hyp}(!,\psi;\lambda_1,\ldots,\lambda_{n-1};\rho_1,\ldots,
\rho_n))\bigr)_{\bar 0}\\ &\to&\mathscr H^{-1} \bigl( {\mathscr
F}(\mathscr L_\psi\otimes j_!\mathrm {inv}^\ast \mathrm
{Hyp}(!,\psi;\lambda_1,\ldots,\lambda_{n-1};\rho_1,\ldots,
\rho_n))\bigr)_{\bar \eta_0}&& \\
&\to& {\mathfrak F}^{(\infty,0)}\bigl(\big(\mathscr L_\psi\otimes
j_!\mathrm {inv}^\ast \mathrm
{Hyp}(!,\psi;\lambda_1,\ldots,\lambda_{n-1};\rho_1,\ldots,
\rho_n)[-1]\big)|_{\eta_\infty}\bigr) &\to & 0.
\end{array}
\end{eqnarray}
We have
\begin{eqnarray}
\begin{array}{cl}
&{\mathfrak F}^{(\infty,0)}\bigl(\bigl({\mathscr L}_\psi\otimes
j_!\mathrm {inv}^\ast \mathrm
{Hyp}(!,\psi;\lambda_1,\ldots,\lambda_{n-1};\rho_1,\ldots,
\rho_n)[-1]\bigr)|_{\eta_\infty}\big)\\
\cong& {\mathfrak F}^{(\infty,0)}\bigl(\bigl(\mathrm {inv}^\ast
{\mathscr K}_{\lambda_1^{-1}\cdots\lambda_{n-1}^{-1}\rho_1\cdots
\rho_n}\bigr)|_{\eta_\infty}\bigoplus \bigoplus_{\lambda}
\bigl({\mathscr L}_\psi\otimes \mathrm{inv}^\ast({\mathscr
K}_\lambda\otimes U(\mathrm
{mult}_0(\lambda)))\bigr)|_{\eta_\infty}\bigr)\\
\cong& {\mathscr
K}_{\lambda_1^{-1}\cdots\lambda_{n-1}^{-1}\rho_1\cdots
\rho_n}|_{\eta_0}
\end{array}
\end{eqnarray}
by \cite{L} 2.4.3 (ii) b) and 2.5.3.1. Applying $\mathrm
{tran}^\ast$ to (15)-(17) and using (14), we see that
$$\bigl(\mathrm {Hyp}(!,\psi;\lambda_1,\ldots,
\lambda_{n-1},1;\rho_1,\ldots,\rho_n)[-1]\big)_{\bar 1}\cong
\big(\mathrm {Hyp}(!,\psi;\lambda_1,\ldots,
\lambda_{n-1},1;\rho_1,\ldots,\rho_n)[-1]\big)_{\bar
\eta_1}^{\mathrm{Gal}(\bar\eta_1/\eta_1)},$$ and we have a short
exact sequence
\begin{eqnarray*}
0&\to& \big(\mathrm {Hyp}(!,\psi;\lambda_1,\ldots,
\lambda_{n-1},1;\rho_1,\ldots,\rho_n)[-1]\big)_{\bar
\eta_1}^{\mathrm{Gal}(\bar\eta_1/\eta_1)} \\
&\to& \big(\mathrm {Hyp}(!,\psi;\lambda_1,\ldots,
\lambda_{n-1},1;\rho_1,\ldots,\rho_n)[-1]\big)_{\bar \eta_1}\\
&\to& (\mathrm{tran}^\ast{\mathscr
K}_{\lambda_1^{-1}\cdots\lambda_{n-1}^{-1}\rho_1\cdots
\rho_n})|_{\bar\eta_1} \to 0.
\end{eqnarray*}

By (13), all the breaks of $(j_!\mathrm {inv}^\ast \mathrm
{Hyp}(!,\psi;\lambda_1,\ldots,\lambda_{n-1};\rho_1,\ldots,
\rho_n)[-1])|_{\eta_\infty}$ are $\leq 1$. It follows that
$${\mathfrak F}^{(\infty,\infty)}\bigl(\bigl(j_!\mathrm {inv}^\ast
\mathrm {Hyp}(!,\psi;\lambda_1,\ldots,\lambda_{n-1};\rho_1,\ldots,
\rho_n)[-1]\bigr)|_{\eta_\infty}\big)=  0.$$ By (12) and \cite{L}
2.5.3.1, we have
$${\mathfrak F}^{(0,\infty)}\bigl(\bigl(j_!\mathrm {inv}^\ast
\mathrm {Hyp}(!,\psi;\lambda_1,\ldots,\lambda_{n-1};\rho_1,\ldots,
\rho_n)[-1]\bigr)|_{\eta_0}\big)\cong \bigoplus_{\rho}
\bigl({\mathscr K}_\rho\otimes U(\mathrm
{mult}_\infty(\rho))\bigr)|_{\eta_\infty}.$$ By \cite{L} 2.3.3.1
(iii), we have
\begin{eqnarray*}
&&\bigl(\mathrm
{Hyp}(!,\psi;\lambda_1,\ldots,\lambda_{n-1},1;\rho_1,\ldots, \rho_n)[-1]\bigr)|_{\eta_\infty}\\
&\cong&  {\mathfrak F}^{(0,\infty)}\bigl(\bigl(j_!\mathrm {inv}^\ast
\mathrm
{Hyp}(!,\psi;\lambda_1,\ldots,\lambda_{n-1};\rho_1,\ldots,\rho_n)[-1]\bigr)|_{\eta_0}\bigr)\\
&&\qquad \qquad \bigoplus {\mathscr
F}^{(\infty,\infty)}\bigl(\bigl(j_!\mathrm {inv}^\ast \mathrm
{Hyp}(!,\psi;\lambda_1,\ldots,\lambda_{n-1};\rho_n,\ldots,
\rho_n)[-1]\bigr)|_{\eta_\infty}\bigr)\\
&\cong& \bigoplus_{\rho} \big({\mathscr K}_\rho\otimes U(\mathrm
{mult}_\infty(\rho))\big)|_{\eta_\infty}.
\end{eqnarray*}
This proves those statements in Proposition 0.7 (iii) about the
behavior at $1$ and at $\infty$ for $\lambda_n=1$. For general
$\lambda_n$, we use the fact that $$\mathrm
{Hyp}(!,\psi;\lambda_1,\ldots, \lambda_n;\rho_1,\ldots,\rho_n)\cong
{\mathscr K}_{\lambda_n}\otimes \mathrm
{Hyp}(!,\psi;\lambda_1\lambda_n^{-1},\ldots,
\lambda_{n-1}\lambda_n^{-1},
1;\rho_1\lambda_n^{-1},\ldots,\rho_n\lambda_n^{-1}).$$ For the
statement in Proposition 0.7 (iii) about the behavior at $0$, we use
the fact that
$$\mathrm
{Hyp}(!,\psi;\lambda_1,\ldots,\lambda_n;\rho_1,\ldots, \rho_n) \cong
\mathrm {inv}^\ast \mathrm {Hyp}(!,\psi^{-1};\rho_1^{-1},\ldots,
\rho_n^{-1};\lambda_1^{-1},\ldots,\lambda_n^{-1}).$$
\end{proof}

\end{document}